\documentclass{osa-article}

\articletype{Research Article}

\usepackage{tikz}
\usetikzlibrary{quotes,angles}
\usepackage{graphicx}
\usepackage{subcaption}
\usepackage[shortlabels]{enumitem}
\usepackage{hyperref}
\usetikzlibrary{arrows}
\usepackage{float}
\usepackage{graphicx}
\usepackage[section]{placeins}
\graphicspath{{NewBraggFigs/}}
\usetikzlibrary{decorations.pathmorphing}
\tikzset{snake it/.style={decorate, decoration=snake}}

\usepackage{sansmath}
\usetikzlibrary{shadings,intersections}

\makeatletter
\newcommand*{\rom}[1]{\expandafter\@slowromancap\romannumeral #1@}

\makeatother

\numberwithin{equation}{section}

\numberwithin{theorem}{section}

\DeclareMathOperator*{\argmin}{arg\,min}

\usepackage{fullpage}

\newcommand{\paren}[1]{\left(#1\right)}               

\newcommand{\va}{\mathbf{a}}
\newcommand{\vb}{\mathbf{b}}

\newcommand{\vx}{\mathbf{x}}
\newcommand{\vy}{\mathbf{y}}
\newcommand{\vz}{\mathbf{z}}
\newcommand{\vd}{\mathbf{d}}
\newcommand{\vs}{\mathbf{s}}

\newcommand{\be}{\begin{equation}}
\newcommand{\ee}{\end{equation}}
\newcommand{\bea}{\begin{eqnarray}}
\newcommand{\eea}{\end{eqnarray}}
\newcommand{\bean}{\begin{eqnarray*}}
\newcommand{\eean}{\end{eqnarray*}}

\newcommand{\bel}[1]{\begin{equation}\label{#1}}

\newcommand{\eel}[1]{{\label{#1}\end{equation}}}

\begin{document}

\title{A novel reconstruction technique for two-dimensional Bragg scatter imaging}

\author{James W. Webber,\authormark{1,$\dagger$} Eric L. Miller\authormark{1,$\ddagger$}}

\address{\authormark{1}The Department of Electrical and Computer Engineering, Tufts University, 161 College Ave, Medford, MA 02155, USA\\
}

\email{\authormark{$\dagger$}james.webber@tufts.edu} 
\email{\authormark{$\ddagger$}eric.miller@tufts.edu} 



\begin{abstract}
Here we introduce a new reconstruction technique for two-dimensional Bragg Scattering Tomography (BST), based on the Radon transform models of [arXiv preprint, arXiv:2004.10961 (2020)]. Our method uses a combination of ideas from multibang control and microlocal analysis to construct an objective function which can regularize the BST artifacts; specifically the boundary artifacts due to sharp cutoff in sinogram space (as observed in [arXiv preprint, arXiv:2007.00208 (2020)]), and artifacts arising from approximations made in constructing the model used for inversion. We then test our algorithm in a variety of Monte Carlo (MC) simulated examples of practical interest in airport baggage screening and threat detection. The data used in our studies is generated with a novel Monte-Carlo code presented here.  The model, which is available from the authors upon request, captures both the Bragg scatter effects described by BST as well as beam attenuation and Compton scatter. 
\end{abstract}

\section{Introduction} 
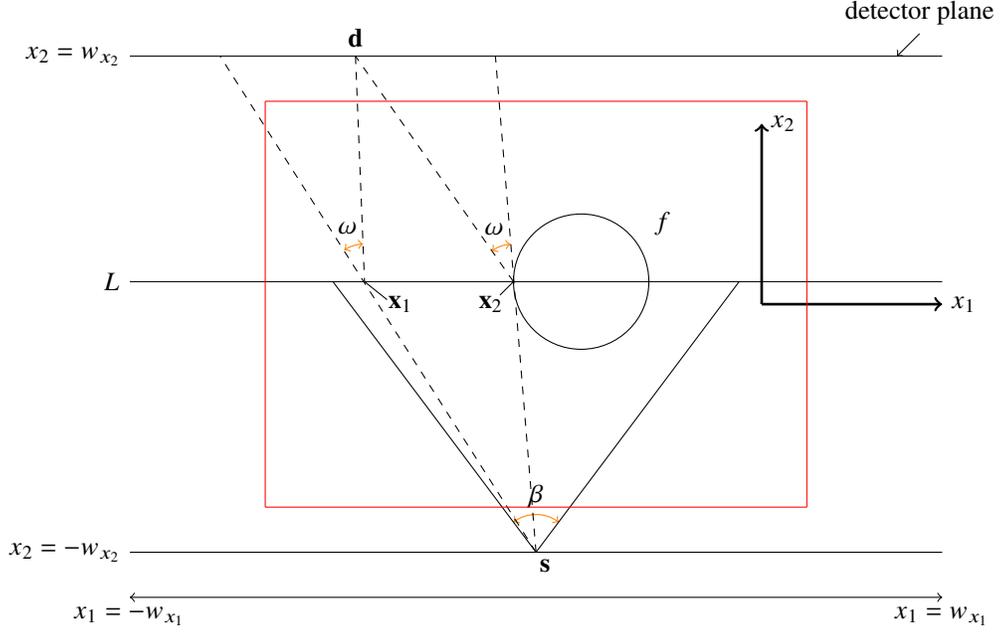
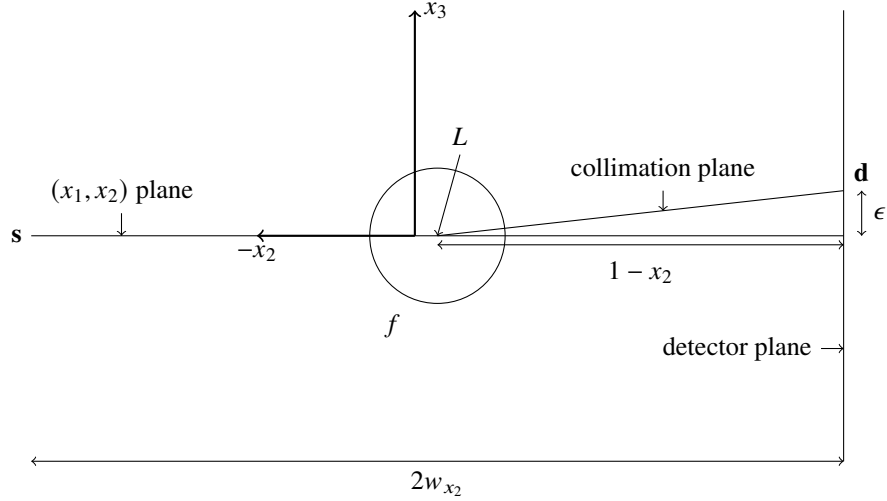
\begin{figure}[!h]
\centering
\begin{subfigure}{1\linewidth}
\centering
\begin{tikzpicture}[scale=6]
\draw  (-0.9,0.5)node[left] {$L$}--(0.9,0.5);
\draw  (-0.9,1)node[left] {$x_2=w_{x_2}$}--(0.9,1);
\draw  [->,line width=1pt] (0.5,0.45)--(0.9,0.45)node[right] {$x_1$};
\draw [->,line width=1pt]  (0.5,0.45)--(0.5,0.85)node[right] {$x_2$};
\coordinate (origo) at (0,-0.1);
\coordinate (pivot) at (-0.45,0.5);
\coordinate (bob) at (0.45,0.5);
\draw pic[draw=orange, <->,"$\beta$", angle eccentricity=1.5] {angle = bob--origo--pivot};
\draw [dashed] (0,-0.1)--(-0.7,1);
\draw [dashed] (-0.4,1)--(-0.38,0.5);
\coordinate (origo2) at (-0.38,0.5);
\coordinate (pivot2) at (-0.7,1);
\coordinate (bob2) at (-0.4,1);
\draw pic[draw=orange, <->,"$\omega$", angle eccentricity=1.5] {angle = bob2--origo2--pivot2};
\node at (-0.3,0.45) {$\vx_1$};
\draw (-0.38,0.5)--(-0.33,0.46);
\draw (0,-0.1)--(-0.45,0.5);
\draw [dashed] (0,-0.1)--(-0.09,1); 
\draw (0,-0.1)--(0.45,0.5);
\draw (-0.9,-0.1) node[left] {$x_2=-w_{x_2}$} --(0.9,-0.1);
\draw (0.1,0.5) circle (0.15);
\node at (0.28,0.63) {$f$};
\node at (-0.1,0.45) {$\vx_2$};
\draw  (-0.08,0.47)->(-0.05,0.5);
\draw [dashed] (-0.05,0.5)--(-0.4,1)node[above] {$\vd$};
\draw [->] (0.85,1.05)node[above] {detector plane}--(0.8,1);
\node at (0.02,-0.13) {$\vs$};
\draw [red] (-0.9+0.3,0)--(-0.9+0.3,0.9);
\draw [red] (0.9-0.3,0)--(0.9-0.3,0.9);
\draw [red] (-0.9+0.3,0.9)--(0.9-0.3,0.9);
\draw [red] (0.9-0.3,0)--(-0.9+0.3,0);
\draw [<->] (-0.9,-0.2)node[below]{$x_1=-w_{x_1}$}--(0.9,-0.2)node[below]{$x_1=w_{x_1}$};
\coordinate (origo1) at (-0.05,0.5);
\coordinate (pivot1) at (-0.4,1);
\coordinate (bob1) at (-0.09,1);
\draw pic[draw=orange, <->,"$\omega$", angle eccentricity=1.5] {angle = bob1--origo1--pivot1};
\end{tikzpicture}
\caption{$(x_1,x_2)$ (source fan-beam) plane cross-section. The source ($\vs$) opening angle is $\beta$ and we have shown two scattering locations at $\vx_1, \vx_2\in L$ with scattering angle $\omega$. The scanning tunnel $[-\frac{2}{3}w_{x_1},\frac{2}{3}w_{x_1}]\times [-\frac{3}{4}w_{x_2},\frac{3}{4}w_{x_2}]$ is shown as a red rectangle.}
\label{fig1}
\end{subfigure}
\begin{subfigure}{1\linewidth}
\centering
\hspace*{0.6cm}
\begin{tikzpicture}[scale=6]
\begin{scope}
\draw  [->,line width=0.8pt] (-0.15,0)--(-0.5,0)node[below] {$-x_2$};
\draw  [->,line width=0.8pt] (-0.15,0)--(-0.15,0.5)node[right] {$x_3$};
\draw  (-1,0) node[left] {$\vs$}--(0.8,0);
\draw  (0.8,-0.5)--(0.8,0.5);
\draw [<->] (0.84,0)--(0.84,0.1)node[above] {$\vd$};
\node at (0.88,0.05) {$\epsilon$};
\draw  (0.8,0.1)--(-0.1,0);
\draw (-0.1,0) circle (0.15);
\node at (-0.2,-0.2) {$f$};
\draw [->] (0.4,0.1) node[above] {collimation plane}--(0.4,0.055);
\draw [->] (-0.8,0.05) node[above] {$(x_1,x_2)$ plane}--(-0.8,0);
\draw [->] (0.75,-0.25) node[left] {detector plane}--(0.8,-0.25);
\draw [<->] (-1,-0.5)--(0.8,-0.5);
\node at (-0.1,-0.55) {$2w_{x_2}$};
\draw [<->] (-0.1,-0.02)--(0.8,-0.02);
\node at (0.35,-0.08) {$1-x_2$};
\draw [->](-0.05,0.18)node[above]{$L$}--(-0.1,0);
\end{scope}
\end{tikzpicture}
\caption{$(x_2,x_3)$ plane cross-section. Note that $L$ is now orthogonal to the page (parallel to $x_1$).}
\label{figyz}
\end{subfigure}
\caption{The X-ray scanner geometry. The scanned object is labelled as $f$, with $\text{supp}(f)\subset [-\frac{2}{3}w_{x_1},\frac{2}{3}w_{x_1}]\times [-\frac{3}{4}w_{x_2},\frac{3}{4}w_{x_2}]$ (the red rectangle in figure \ref{fig1}). The detectors are collimated to planes, and the scattering events occur along lines $L=\{x_2=x',x_3=0\}$, for some $-w_{x_2}<x'<w_{x_2}$. The scatter from $L$ is measured by detectors $\vd\in\{x_2=1,x_3=\epsilon\}$, for some $\epsilon>0$.}
\label{figmain}
\end{figure}

In this paper we present a novel reconstruction technique for two-dimensional BST, based on the generalized Radon transform models of \cite{webber2020bragg}. Our method employs a combination of ideas from multibang control and microlocal analysis to derive new regularization penalties, which prove to be effective in combatting the high level of noise and systematic error in BST data  (e.g. error due to beam attenuation and Compton scatter, the two primary physical mechanisms not accounted for under the BST model). The BST model arises from a scanning geometry (first introduced in \cite{webber2020bragg}) which uses translating sources to inspect what is well approximated as a line in image space. See figure \ref{figmain}. A more detailed description of the sensing geometry is given later in section \ref{geometry}. Here BST refers to the imaging of the Bragg differential cross section function (denoted by $f$ in this paper) from Bragg scatter data, and is not exclusive to the sensing geometry of \cite{webber2020bragg}.

The literature considers a variety of reconstruction techniques and experimental methods in BST \cite{hassan2016snapshot,greenberg2013snapshot,maccabe2012pencil,greenberg2014structured,greenberg2013coding,greenberg2015optimization}.
In \cite{greenberg2013snapshot} the authors consider snap shot, pencil beam, coded aperture imaging of crystalline powders, such as NaCl and Al. The technique, referred to as Coded Aperture X-ray Scatter Imaging (CAXSI), uses a pencil beam source to illuminate a small sample of crystalline powder. The scattered rays are then passed through a coded aperture mesh, and the resulting intensity is recorded by a linear (1-D) array of energy-resolved detectors. The use of coded aperture offers information about the scattering directions of the incoming photons, and thus improves the problem stability. Mathematically, the coded aperture modelling is represented by a kernel weighting ($t$ in \cite[equation (3)]{greenberg2013snapshot}). The physical modeling then leads the authors to a linear inverse problem to solve for the Bragg differential cross section function $f$. To obtain a solution the authors apply Total Variation (TV) regularization and minimize the Poisson log-likelihood function. An iterative, Generalized Expectation-Maximization (GEM) algorithm is then implemented to minimize the objective, with good results on experimental data. 

In \cite{maccabe2012pencil} CAXSI is considered, using an experimental setup much like that of \cite{greenberg2013snapshot} with a planar (2-D) array of detectors (the detectors are not energy-resolved). A Generalized Maximum Likelihood (GML) estimator is then applied to estimate $f$. The GML algorithm is a multiplicative, iterative update intended to match the BST model to the data with a mean background and Poisson statistics applied. GML is similar to GEM, except GEM applies an additional maximization step \cite[equation (10)]{greenberg2013snapshot} after each iteration, which amounts to Poisson denoising with TV.

In \cite{hassan2016snapshot,greenberg2014structured,greenberg2013coding,greenberg2015optimization} a number of CAXSI variations are considered, for example, using fan beam sources in \cite{hassan2016snapshot,greenberg2014structured}. In these works, a set of ``reference" differential cross section curves is used in the reconstruction algorithm, whereby each pixel is assigned to a material from the reference library based on the normalized correlation between the reconstructed and reference cross section values at the pixel. In this paper we do not assume knowledge of a cross section reference library, and wish to keep the material properties general. We assume only that the form factor curves are $L^2$ functions.

We introduce a new regularization scheme for small sample ($<3\text{cm}$ dimension), low effective atomic number (e.g. H-C-N-O compounds) BST. The technique we propose uses a collection of ideas from Compressed Sensing (CS), multibang control and microlocal analysis. The application of such ideas has not yet been investigated in the BST literature. Further, in the experiments conducted in \cite{hassan2016snapshot,greenberg2013snapshot,maccabe2012pencil,greenberg2014structured,greenberg2013coding,greenberg2015optimization}, typically only point-like samples are considered (i.e. $<1\text{cm}$ dimension). We consider objects up to 3cm in width, length and depth. We also make no use of coded apertures in our experimental setup to restrict the scattering directions and increase problem stability. The linear collimation technology of the scanner of figure \ref{figmain} is used to the same effect, to restrict the scattering sites to lines parallel to $x_1$ in the $(x_1,x_2)$ plane.

The reconstruction target $f$ (illustrated as a sphere in figure \ref{figmain}) is three-dimensional $f=f(q,x_1,x_2)$, where $q$ denotes the momentum transfer of the scattering interaction (see equation \eqref{qdef}). As is done in \cite{webber2020bragg}, the recovery of $f$ is obtained slice by slice on planes orthogonal to the $x_2$ direction. We consider only 2-D reconstruction here (in $(q,x_1)$ space), for a variety of $x_2$ values, in our simulations, and leave the piecing together of the 2-D slices (to form a full 3-D image) to future work.  A major part of the regularization idea we propose is to assume that $f$ is separable as a function of $q$ and $(x_1,x_2)$ (in the spatial domain). We model the spatial component of $f$ using an overcomplete piecewise constant dictionary. This is a standard idea in CS \cite{8834789,li2011compressed,Herrholz2014}, although typically in CS the library (basis) functions (or atoms \cite{Herrholz2014}) are chosen to span the whole imaging domain. The basis functions we use only cover the $(x_1,x_2)$ domain, with $f$ having more general $L^2$ properties in the $q$ domain. The assumptions made here regarding the piecewise constant model for $f$ discussed above are consistent with the BST literature \cite{hassan2016snapshot,greenberg2013snapshot,maccabe2012pencil,greenberg2014structured,greenberg2013coding,greenberg2015optimization}, and what is expected in practice. That is, we expect $f$ to be expressible as a finite sum of characteristics in the spatial domain. For example, $f$ could be a block of explosives in airport baggage or a sample of narcotic powder (e.g. fentanyl) in mail at customs. 

Another major component of our regularization penalty is the use of ideas from multibang control \cite{clason2014multi,clason2016convex,clason2018total}. The multibang penalty \cite{clason2014multi} is used to enforce solutions whereby the function outputs are constrained to a finite set. We aim to apply multibang ideas here to enforce the piecewise constant structure of $f$ in $(x_1,x_2)$ space, as discussed in the last paragraph. We do this by defining a set of binary switches $a_j$ to either activate or deactivate a characteristic function from our library. The multibang penalties are applied to enforce binary solutions for the $a_j$. So the finite set of solutions for the fitted $a_j$, in our case, is $\{0,1\}$. As our proposed objective function has smoothly defined gradients, we seek a relaxed alternative to the multibang penalty (the multibang penalty of \cite{clason2014multi} is not smooth), which we introduce later in section \ref{reconmethod}. 

In addition to CS and multibang techniques, we also employ filtering ideas from microlocal analysis \cite{frikel2013characterization,borg2018analyzing}, to suppress the boundary artifacts typically observed in BST reconstruction, e.g., as observed in \cite{webber2020microlocal} in reconstructions from Bragg integral data. The filtering techniques from the literature are shown in section \ref{microlocal} to offer significantly improved image quality and artifact reduction in reconstructions from Monte Carlo data.

The remainder of the this paper is organized as follows. In section \ref{methodology} we explain our methodology. This includes a review of the physical model and Bragg transform from \cite{webber2020bragg} in section \ref{braggtrans}, before moving on to explain our new reconstruction method in section \ref{reconmethod}. The reconstruction technique we propose is formalized as an algorithm in section \ref{algorithm1}, and the pre-processing using microlocal filters is explained in section \ref{microlocal}. In section \ref{results} we present our results on a wide variety of Monte Carlo and analytic simulations of interest in threat detection, and give comparison to a TV regularized solution.

\section{Methodology}
\label{methodology}
In this section we introduce our reconstruction technique and explain the filtering techniques used as data pre-processing from microlocal analysis. First we review the sensing geometry and Bragg transform of \cite{webber2020bragg}, and introduce some notation.
\subsection{The sensing geometry}
\label{geometry}
The scanner of figure \ref{figmain} is equipped with linear detector collimation technology, which we will refer to as ``Venetian blind" type collimation. The scanned object $f$ travels through the scanning tunnel (the red rectangle in figure \ref{fig1}) in the $x_3$ direction on a conveyor belt, and is illuminated by a line of X-ray sources, located opposite a plane of detectors. The scanner sources (with coordinate $\vs$) are fixed and switched along $\{x_2=-w_{x_2},x_3=0\}$, and are assumed to be polychromatic 2-D fan-beam (in the $(x_1,x_2)$ plane) with opening angle $\beta$. The detectors (with coordinate $\vd$) are assumed to be energy-resolved and lie on the $\{x_2=w_{x_2}\}$ plane, with small (relative to the scanning tunnel size) offset $\epsilon$ in the $x_3$ direction. The detectors are collimated to record photons which scatter on planes in $\mathbb{R}^3$, and the planes of collimation are orientated to intersect the source $(x_1,x_2)$ plane along horizontal lines (parallel to $x_1$). Hence the photon arrivals measured by the scanner detectors are scattered from horizontal lines embedded in the $(x_1,x_2)$ plane. An example $\epsilon$ is shown in figure \ref{figyz}, which maps to the line $L=\{x_2=0,x_3=0\}$ at the half way point.

\subsection{The Bragg transform}
\label{braggtrans}
Let $\mathcal{I}=[-w_{x_2},w_{x_2}]$, let $\mathfrak{E}=[E_m,E_M]$ be the energy range, and let $\Phi : x_2\to \epsilon$ be a diffeomorphic map from the scanned line profile $x_2$ to the detector array position $\epsilon$. Let $\vx=(x_1,x_2)$. Then the Bragg transform $\mathfrak{B}_{a} : L^2_0(\mathfrak{E}\times[-w_{x_1},w_{x_1}]\times \mathcal{I})\to L^2((0,\infty)\times[-w_{x_1},w_{x_1}]^2\times\Phi(\mathcal{I}))$ defines a mapping from the target $f$ to the Bragg scatter measured by the scanner of figure \ref{figmain} \cite[page 7]{webber2020bragg}
\begin{equation}
\label{equBG}
\begin{split}
\mathfrak{B}_af(E,s_1,d_1,\Phi(x_2))=\int_{\mathbb{R}}\chi_{[-w,w]}(x_1-s_1)&I_0(E,\vx)P(\theta(\vd,\vs,\vx))\mathrm{d}\Omega_{\vx,\vd}\\
&\times f\left(\frac{E\sin\theta(\vd,\vs,\vx)}{hc},\vx\right)\mathrm{d}x_1,
\end{split}
\end{equation}
where $\vs=(s_1,-w_{x_2},0)$, $\vd=(d_1,w_{x_2},\Phi(x_2))$, $\chi_{S}$ denotes the characteristic function on a set $S$ and $f(q,\vx)=n_c(\vx)F(q,\vx)$ is the number of cells per unit volume ($n_c$) multiplied by the Bragg differential cross section ($F$). Here $h$ is Planck's constant and $c$ is the speed of light in a vacuum. 

We consider the recovery of the 2-D functions $f(\cdot,\cdot,x_2)$ from $\mathfrak{B}_af(\cdot,\cdot,\cdot,\Phi(x_2))$, for each $x_2\in\mathcal{I}$. That is we consider the slice-by-slice reconstruction of $f$ from the 4-D data $\mathfrak{B}_af$. We focus exclusively on 2-D reconstruction here. This is to say that we do not consider the piecing together of the 2-D slices to form a full 3-D image. This is left to future work. With this in mind we adopt the short-hand notation $f(q,x_1)=f(\cdot,\cdot,x_2)$, for some fixed $x_2\in\mathcal{I}$.  

The remaining terms are defined as follows. The source width $w$ is determined by the source opening angle $\beta$ (see figure \ref{fig1})
\begin{equation}
w(x_2)=(1+x_2)\tan\frac{\beta}{2}.
\end{equation}
The solid angle is
\begin{equation}
\label{sangle}
\mathrm{d}\Omega_{\vx,\vd}=D_A\times \frac{((\vx,0)-\vd)\cdot (0,-1,0)^T}{|(\vx,0)-\vd|^3},
\end{equation}
where $D_A$ is the detector area, and the Bragg angle ($\theta=\frac{\omega}{2}$) is determined by
\begin{equation}
\label{bangle}
\cos 2\theta(\vd,\vs,\vx)=\frac{((\vx,0)-\vs)\cdot(\vd-(\vx,0))}{|((\vx,0)-\vs)||(\vd-(\vx,0))|}.
\end{equation}
The polarisation factor $P(\theta)$ is given by
\begin{equation}
P(\theta)=\frac{1+\cos^22\theta}{2}
\end{equation}
and the initial source intensity is
\begin{equation}
\label{I0}
I_0(E,\vx)=\frac{I_0(E)}{|\vs-\vx|^2}=\frac{I_0(E)}{x_1^2+(x_2+1)^2},
\end{equation}
where $I_0>0$ is the initial energy spectrum (e.g. a Tungsten target X-ray tube). The momentum transfer is defined 
\begin{equation}
\label{qdef}
q=\frac{E}{hc}\sin\theta,
\end{equation}
where $E$ is given in units of kilo-electron-volts (keV) and $q$ is given in units of inverse Angstroms ($\AA^{-1}$). So $hc$ is the conversion factor from $\AA^{-1}$ to keV. Equation \eqref{qdef} is derived from the Bragg equation \cite{bragg1913reflection}
\begin{equation}
\label{braggequ}
 \frac{hc}{E}=\lambda=2d_{H}\sin\theta,
\end{equation}
where $\lambda$ is the photon wavelength, and $d_H$ is the spacing between the reflection planes within the crystal. For example, for cubic structures
\begin{equation}
\label{dH}
d_{H}=\frac{a_0}{\sqrt{h^2+k^2+l^2}},
\end{equation}
where $H=(h,k,l)$ is the Miller index and $a_0$ is the uniform lattice spacing of the crystal.

The operator $\mathfrak{B}_a$ is the same as considered in \cite[equation (3.13)]{webber2020bragg}, but with the attenuation terms $A_1$ and $A_2$ removed from the modelling (using the notation of \cite{webber2020bragg}). We do this as later in our simulations we assume no prior knowledge of the attenuation map. Further, as discussed in \cite{webber2020bragg}, the neglection of attenuation effects from the modelling is needed to prove linear invertiblity.

\begin{figure}[!h]
\centering
\begin{subfigure}{0.32\textwidth}
\includegraphics[width=0.9\linewidth, height=4cm, keepaspectratio]{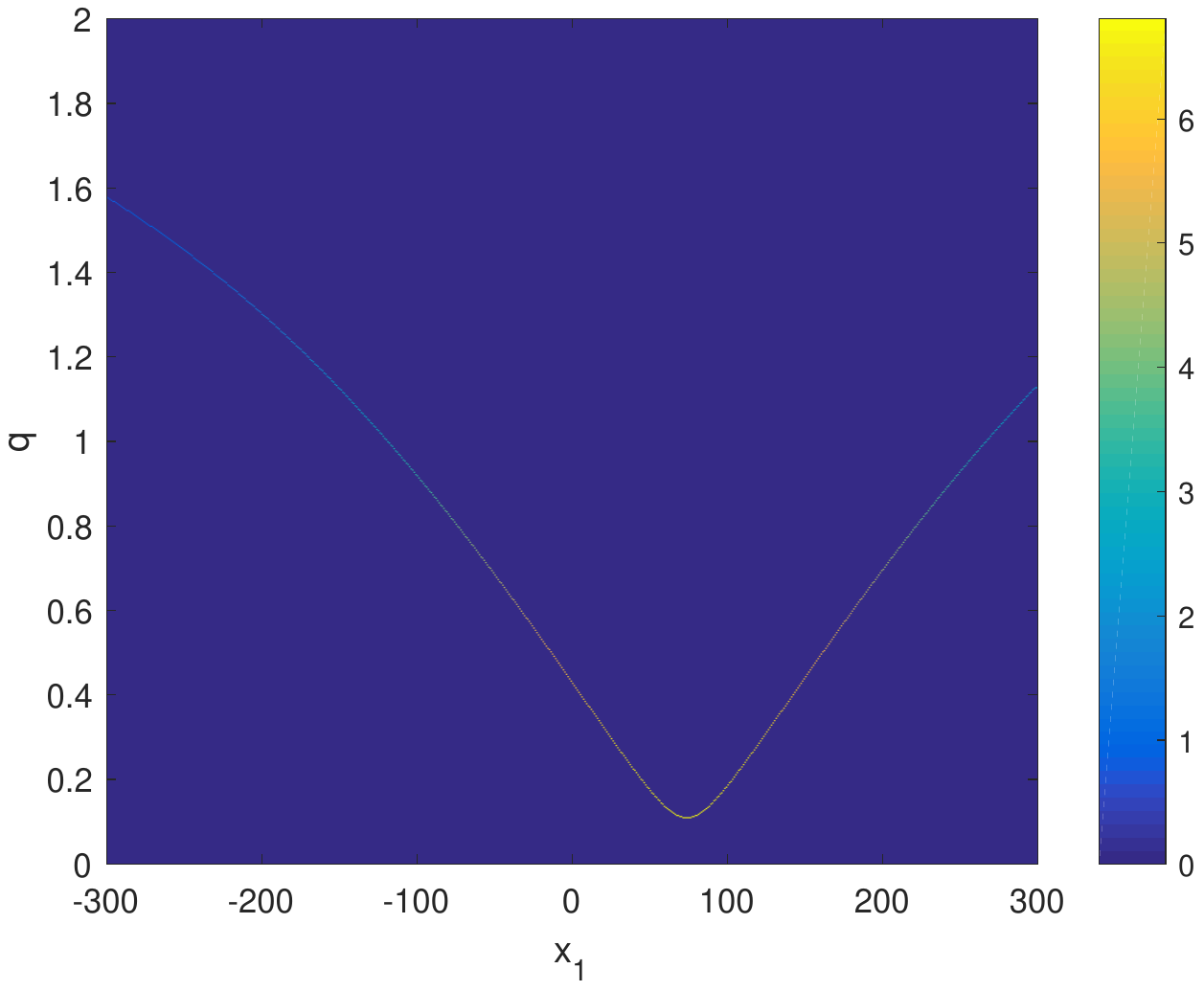}
\end{subfigure}
\begin{subfigure}{0.32\textwidth}
\includegraphics[width=0.9\linewidth, height=4cm, keepaspectratio]{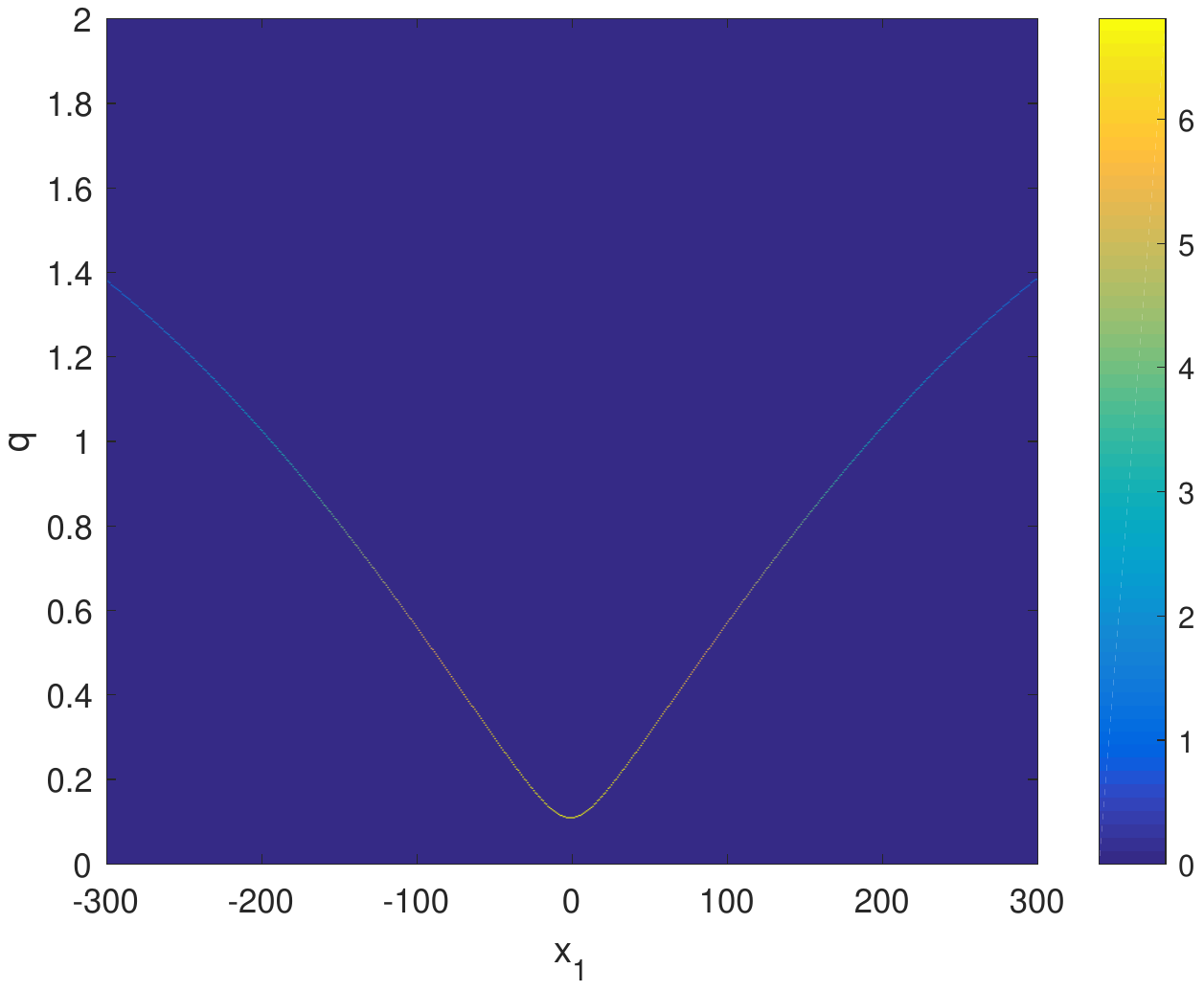} 
\end{subfigure}
\begin{subfigure}{0.32\textwidth}
\includegraphics[width=0.9\linewidth, height=4cm, keepaspectratio]{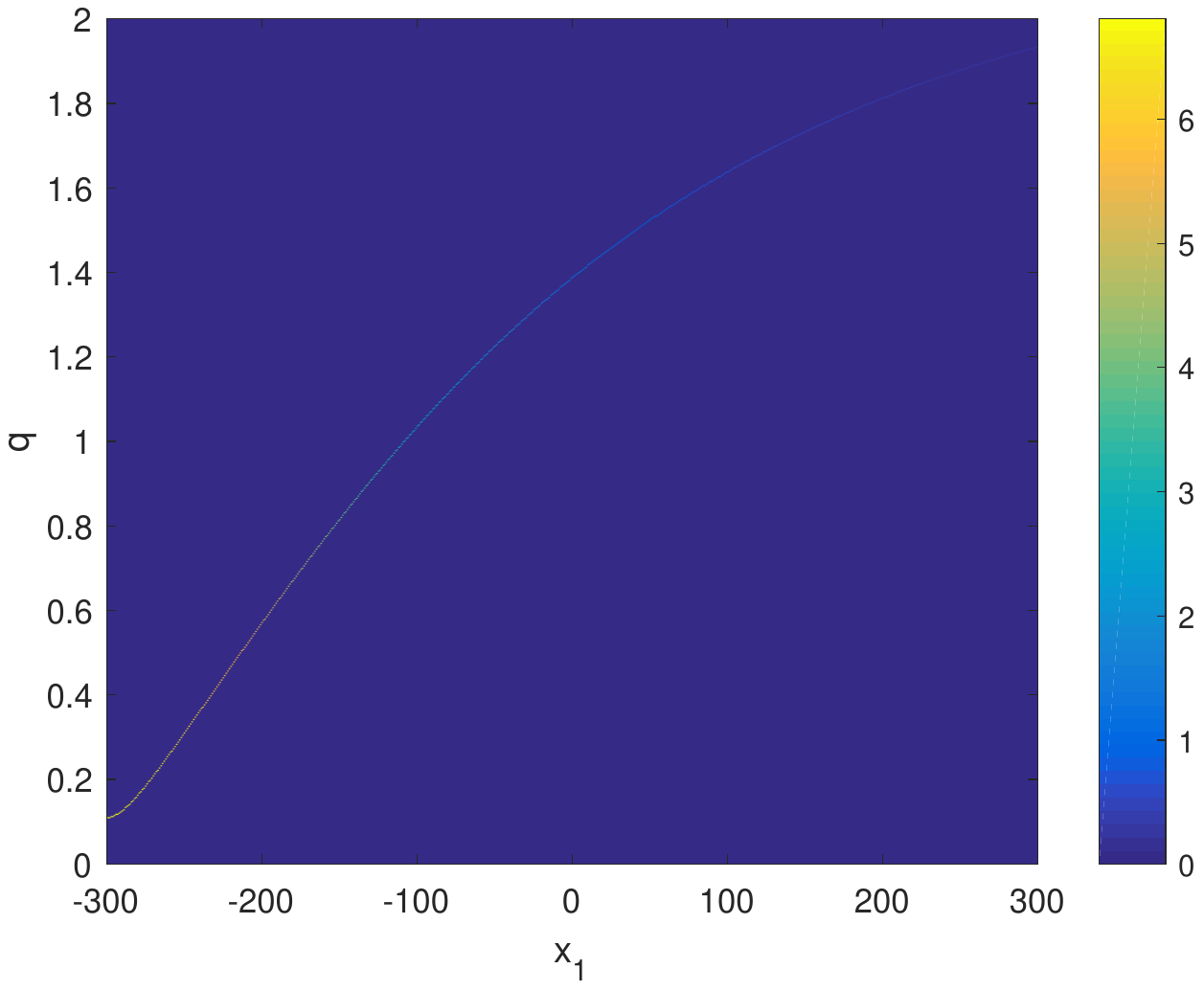}
\end{subfigure}
\begin{subfigure}{0.32\textwidth}
\includegraphics[width=0.9\linewidth, height=4cm, keepaspectratio]{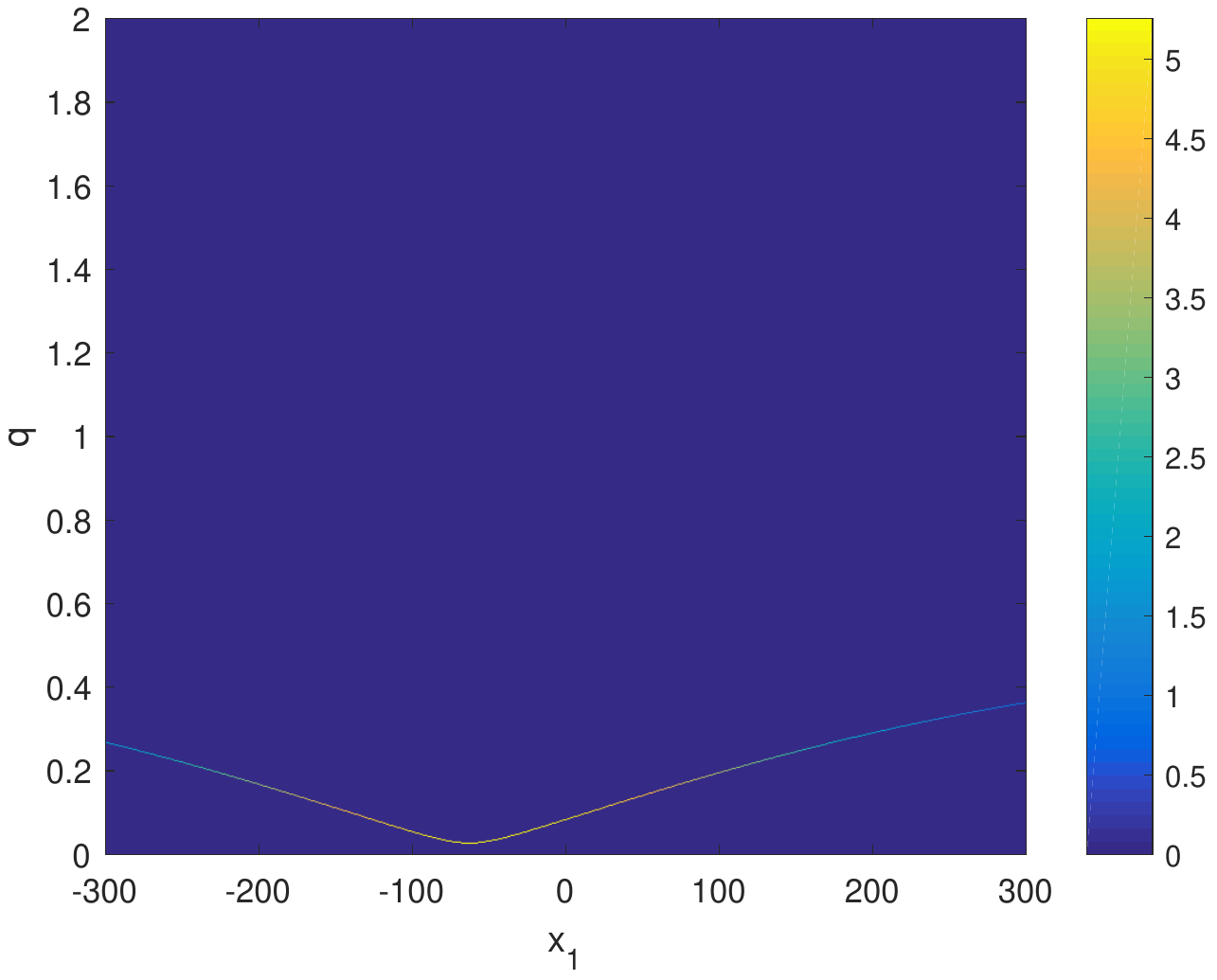}
\end{subfigure}
\begin{subfigure}{0.32\textwidth}
\includegraphics[width=0.9\linewidth, height=4cm, keepaspectratio]{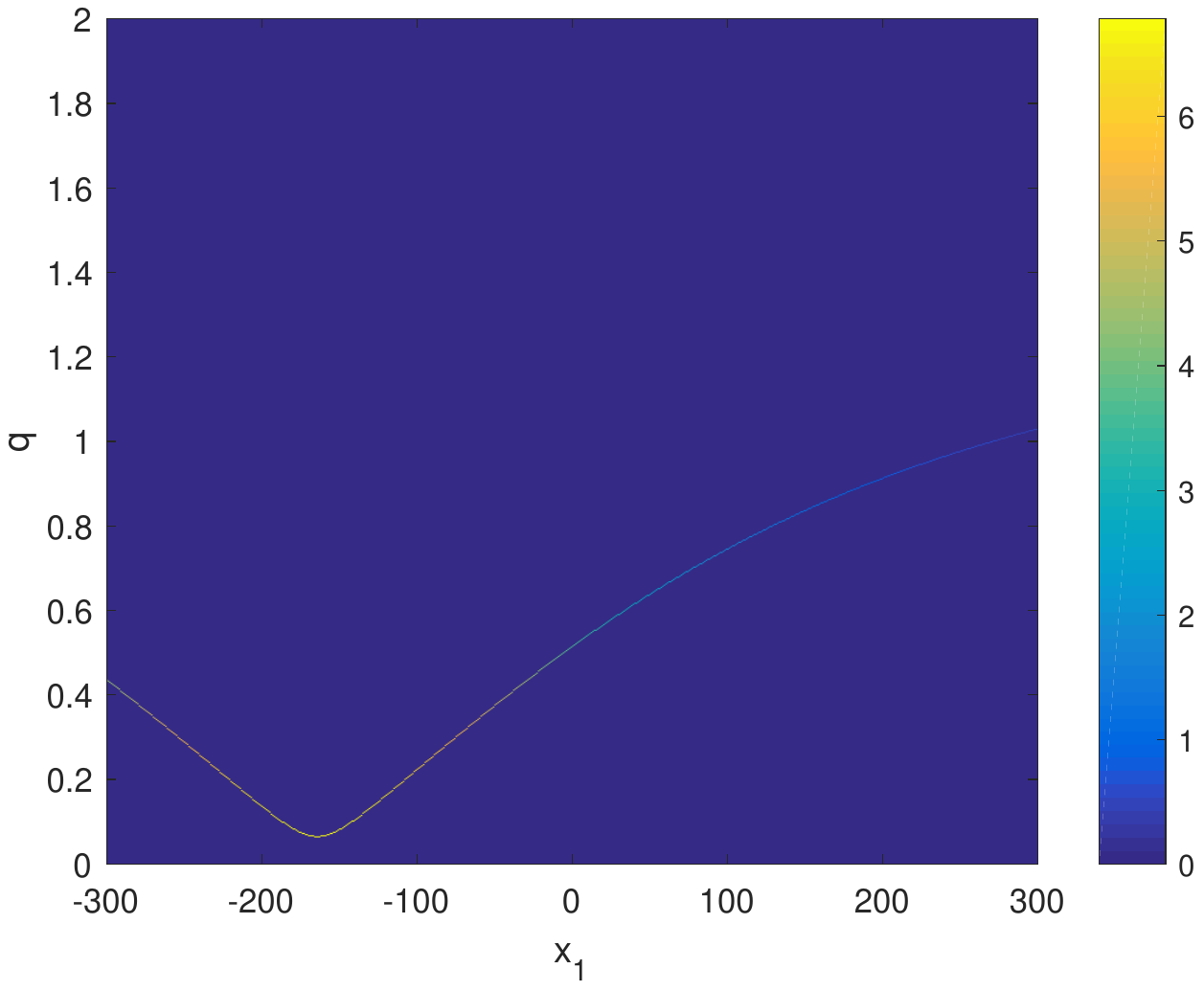} 
\end{subfigure}
\begin{subfigure}{0.32\textwidth}
\includegraphics[width=0.9\linewidth, height=4cm, keepaspectratio]{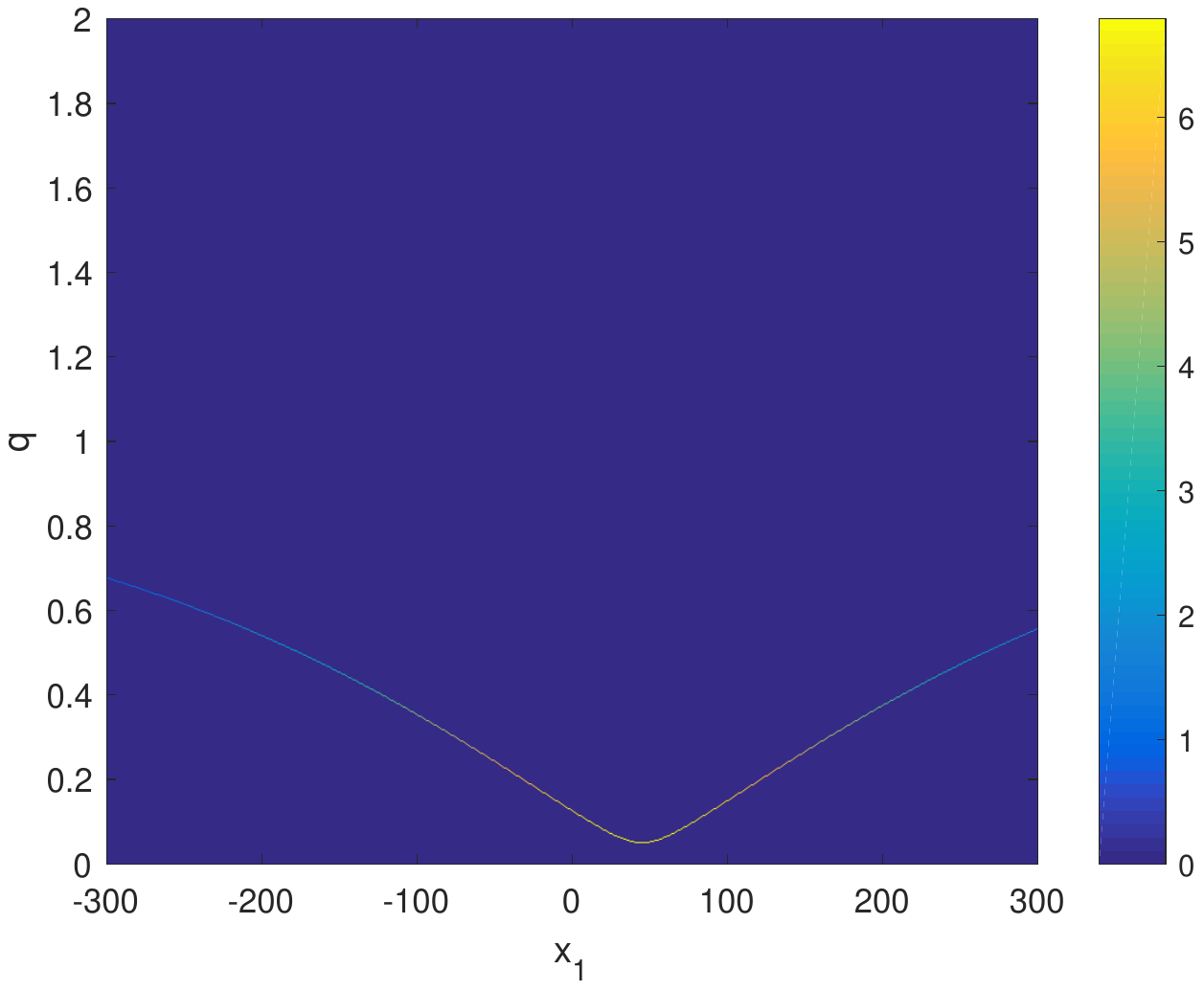}
\end{subfigure}
\caption{Weighted Bragg curve examples for varying $E$, $s_1$ and $d_1$. $x_2=0$ is fixed. The figures were formed by reshaping the rows of $A$ into 2-D images. The curves on the top row are chosen so that $s_1=d_1$ as is considered in \cite{webber2020bragg}. We see that the curves on the top row have the same shape as those shown in \cite[Figure 7]{webber2020bragg}.}
\label{BC1}
\end{figure}
\subsection{The reconstruction method}
\label{reconmethod}
Throughout this paper, $f(q,x_1)$ (when discretized) will be represented as an $n\times m$ image, with $n$ the number of $q$ samples, and $m$ the number of $x_1$ samples. Let us fix $x_2\in \mathcal{I}$, and let $A\in \mathbb{R}^{p\times (mn)}$ denote the discretized form of the linear operator $\mathfrak{B}_a$. See figure \ref{BC1} for an illustration of a discretized Bragg operator in $(q,x_1)$ space. 

The regularization method we propose is derived from a set of a-priori assumptions regarding the target function $f$. We assume that $f$ is of the form
\begin{equation}
\label{fform1}
f(q,x_1)=\sum_{j=1}^la_jf_j(q)\chi_j(x_1),
\end{equation}
where $f_j\in L^2(\mathfrak{E})$, and $\chi_j=\chi_{I_j}$ with $I_j=[-w_j+x^c_j,w_j+x^c_j]$ an interval with width $w_j$ and center $x^c_j$. The $f_j$ have the form of a delta-comb (see \cite[equation (3.10)]{webber2020bragg}). See figure \ref{Fq1} for some example $f_j$ curves. As discussed in the introduction, the form \ref{fform1} for $f$ is consistent with the BST literature and what is expected in practice, and thus it is reasonable to assume that $f$ can be expressed by the expansion \ref{fform1}.

We introduce a library (finite set) of characteristic functions $\chi_1,\ldots,\chi_l$ from which $f$ can be formed. The $a_j\in\{0,1\}$ act as binary switches to either activate or deactivate characteristic $j$ from the library. The $\chi_j$ library is chosen to comprehensively cover the support of $f$ in $x_1$, but to also be restrictive enough to sufficiently regularize the solution. That is, we choose the lowest cardinality $x^c_j, w_j$ set such that the characteristic centers and widths of interest are accurately represented. The $x^c_j, w_j$ set used later in section \ref{results} in our simulations is shown to be sufficient to cover sufficiently the support of $f$ in $x_1$, with good results, for the most part, on phantoms which are comprised of characteristics lying outside of the chosen library.

Let $\vz_j$ denote the vectorized form of $\chi_j$ and let $C_j=\paren{\vz_{j1}I_n,\ldots,\vz_{jm}I_n}^T$, where $\vz_{ji}$ is the $i^{\text{th}}$ entry of $\vz_j$ and $I_n$ is the $n\times n$ identity matrix. Then we define $A_j=AC_j$ as the restriction of  $A$ to characteristic $\chi_j$. Let $\mathcal{Z}=(\vz_1,\ldots,\vz_l)$ be the matrix with $\vz_j$ as columns. We define the Gram matrix
\begin{equation}
\label{Gram}
G=\mathcal{Z}^T\mathcal{Z}=(\vz_1,\ldots,\vz_l)^T(\vz_1,\ldots,\vz_l).
\end{equation}
Then $G$ is such that $G_{ij}\neq 0$ if $\chi_i$ and $\chi_j$ intersect and $G_{ij}= 0$ otherwise. We are now ready to define our objective function.


We propose to minimize the functional
\begin{equation}
\begin{split}
\label{equ1}
\mathcal{C}(\va,Y)=\sum_{k=1}^p\Bigg[\paren{\sum_{j=1}^{l}a_jA_j\vy_j}_k -b_k\log&\paren{\sum_{j=1}^{l}a_jA_j\vy_j}_k\Bigg]+\lambda\sum_{j=1}^l\|\vy_j\|_1\\
&+\alpha\sum_{j=1}^la_j(1-a_j)+\gamma\sum_{i<j}G_{ij}a_ia_j,
\end{split}
\end{equation}
where $\va=(a_1,\ldots,a_l)$, $a_j\in[0,1]$, $Y=(\vy_1,\ldots,\vy_l)\in\mathbb{R}^{n\times l}_+$ and $\vy_j$ is the discretized form of $f_j$. The negative Poisson log-likelihood function in the first term of \eqref{equ1} is included since we expect the photon arrivals to follow a Poisson noise model, as is, for example, used in the BST literature \cite{hassan2016snapshot,greenberg2013snapshot,maccabe2012pencil,greenberg2014structured,greenberg2013coding,greenberg2015optimization}, and also in Positron Emission Tomography (PET) \cite[page 5]{ehrhardt2014joint}. Here the notation $\paren{\textbf{b}}_k=b_k$ of \eqref{equ1} denotes the $k^{\text{th}}$ entry of $\textbf{b}$. The penalty term
\begin{equation}
\label{MB2}
\text{MB}_{1}(\va)=\sum_{j=1}^la_j(1-a_j)
\end{equation}
is included to enforce binary solutions for $a$. The idea follows a similar intuition to that of multibang control \cite{clason2014multi,clason2016convex,clason2018total}, where the authors seek solutions with a finite set of values. The multibang penalty \cite{clason2014multi} is defined (using the notation of \cite{clason2014multi}) as 
\begin{equation}
\label{MB}
\text{MB}_0(u)=\int_{\Omega}\prod_{j=1}^d|u(x)-u_j|_0,
\end{equation}
where $\Omega\in\mathbb{R}^n$ is the region of interest, $\{u_1,\ldots,u_d\}$ is the finite set of solutions considered, and $|u|_0=1$ if $u=0$, $|u|_0=0$ otherwise. In our case $\{u_1,u_2\}=\{0,1\}$ and $d=2$. The objectives of \cite{clason2014multi,clason2016convex,clason2018total} are then minimized using a semismooth Newton approach, as proposed in \cite{clason2014multi}. As the sum of the remaining terms of \eqref{equ1} (i.e. with $\alpha\text{MB}_1(\va)$ removed from the summation) has a trivially computable gradient, we seek a multibang type penalty with smooth gradients so that the solution to \eqref{equ1} can be obtained using the quasi-Newton methods of \cite{byrd1995limited}. With this in mind we introduce the ``relaxed" multibang regularizer $\text{MB}_1$ here, with smoothing parameter $\alpha$  (increasing $\alpha$ more harshly enforces binary solutions and vice-versa). The term
\begin{equation}
\text{GM}(\va)=\sum_{i<j}G_{ij}\paren{\va\va^T}_{ij}=\sum_{i<j}G_{ij}a_ia_j
\end{equation}
is included as a soft constraint in \eqref{equ1} to enforce no overlap between the activated characteristics. We do this, as two objects cannot occupy the same space (i.e. $f$ is well-defined). In the simulations conducted in section \ref{results}, $\gamma$ is set to some large value (orders of magnitude greater than $\alpha$ and $\lambda$) to ensure there is no overlap in the $\chi_j$. Finally, the $L^1$ regularization penalty (with parameter $\lambda$) is included to enforce sparsity in the $q$ domain, as we expect the $\vy_j$ to be sparse given the delta-comb structure of the $f_j$ for crystalline material (see figure \ref{Fq1}).

To minimize $\mathcal{C}$ we use the L-BFGS-B method of \cite{byrd1995limited} with box contraints, setting non-negativity constraints on the $\vy_j$ and constraining $a_j\in [0,1]$. A solution is obtained by first fixing $\va$, and solving \eqref{equ1} for $\vy$ for a number of inner iterates $n_2$, then fixing the $\vy$ output and solving for $\va$ over $n_2$ iterates. The entire process is then repeated for $n_1$ outer iterations until convergence is reached. A more detailed explanation of the reconstruction algorithm is given below.

\subsection{Algorithm: 2-D Bragg Scatter Reconstruction (2DBSR)}
\label{algorithm1}
Initialize $b$, $\alpha$, $\lambda$, $\gamma$, $n_1$ and $n_2$. Initialize $\va^0$, $\vy^0$, and characteristic library $\mathcal{Z}=(\vz_1,\ldots,\vz_l)$. Define $G=\mathcal{Z}^T\mathcal{Z}$. Then for $n_{\text{outiter}}=1$ to $n_1$ do:\\
\\
{\textbf{Stage 1}}:
\begin{itemize}
\item Fix $\va^0$ and define $\mathcal{A}=[a^0_{1}A_1,\ldots,a^0_{l}A_l]$,
$$\mathcal{F}(\vy)=\sum_k\paren{\mathcal{A}\vy}_k -b_k\log\paren{\mathcal{A}\vy}_k+\lambda\|\vy\|_1$$
\item Calculate the gradients $\nabla_{\vy}\mathcal{F}$.
\item Run $n_2$ iterations of L-BFGS-B with $\vy^0$, $\mathcal{F}$ and $\nabla_{\vy}\mathcal{F}$ as input, constraining the solution $\vy\in [0,\infty)^{n\times l}$.
\item Set $\vy^0=\vy$.
\end{itemize}
{\textbf{Stage 2}}:
\begin{itemize}
\item Fix $\vy^0$ and define $\mathcal{Y}=[A_1\vy^0_1,\ldots,A_l\vy^0_l]$,
$$\mathcal{G}(\va)=\sum_k\paren{\mathcal{Y}\va}_k -b_k\log\paren{\mathcal{Y}\va}_k+\alpha\sum_{j=1}^la_j(1-a_j)+\gamma\sum_{i<j}G_{ij}a_ia_j$$
\item Calculate the gradients $\nabla_{\va}\mathcal{G}$.
\item Run $n_2$ iterations of L-BFGS-B with $\va^0$, $\mathcal{G}$ and $\nabla_{\va}\mathcal{G}$ as input, constraining the solution $\va\in [0,1]^{l}$.
\item Set $\va^0=\va$.
\end{itemize}
Equivalently, initialize the input parameters, then repeat stages 1 and 2 above $n_1$ times until convergence in reached. Note that in stage 1, since the entries of $\vy$ are constrained to be positive, the gradient of $\|\vy\|_1$ is defined and smooth at $\vy=0$, i.e., 
$$\nabla_{\vy}\|\vy\|_1=\underbrace{(1,\ldots,1)^T}_{n\times l\ \text{times}}$$
on $[0,\infty)^{n\times l}$.

\subsection{Artifact reduction and data pre-processing using microlocal filters}
\label{microlocal}
Here we discuss the idea to apply filtering techniques from microlocal analysis to the Bragg data, with the aim to suppress boundary type artifacts in the reconstruction, e.g., as are discovered in \cite{webber2020microlocal} in reconstructions from Bragg curve integral data.

In the literature on microlocal analysis in 2-D limited data X-ray CT \cite{frikel2013characterization,borg2018analyzing}, the authors make use of smoothing filters to reduce streaking artifacts which appear along straight lines at the boundary of the data set. Let $R$ denote the classical Radon transform and let $S\subset S^1\times \mathbb{R}$ be the subset of sinogram space for which $Rf\in L^2(S^1\times \mathbb{R})$ is known, where $f$ is the reconstruction target. Then the reconstruction in limited data CT is typically obtained by direct filtered backprojection
\begin{equation}
\label{FBP}
\tilde{f}=R^*\Lambda \chi_S Rf,
\end{equation}
where $\Lambda$ is the standard ramp filter \cite[chapter 3]{kak1988principles}. That is, the missing data is set to zero and then the inverse Radon transform is applied to recover $\tilde{f}\approx f$. Using a multiplicative filter with such a sharp cutoff as $\chi_S$ has been shown to produce heavy streaking artifacts in the reconstruction, which appear along lines corresponding to data points at the boundary of $S^1\times \mathbb{R}\backslash S$ (i.e. where $Rf$ is not known). In \cite{frikel2013characterization,borg2018analyzing} it is proposed to suppress such artifacts by replacing $\chi_S$ in \eqref{FBP} with $\psi\in C^{\infty}\paren{S^1\times \mathbb{R}}$, such that $\psi(s,\theta)=0$ for $(s,\theta)\in \paren{S^1\times \mathbb{R}}\backslash S$ and $\psi(s,\theta)=1$ on most of the interior of $S$ (away from the boundary), with $\psi$ smoothly tending to zero near the boundary of $S$. We use a similar idea here for data pre-processing. 
\begin{figure}[h!]
\centering
\includegraphics[ width=0.4\linewidth, height=0.4\linewidth, keepaspectratio]{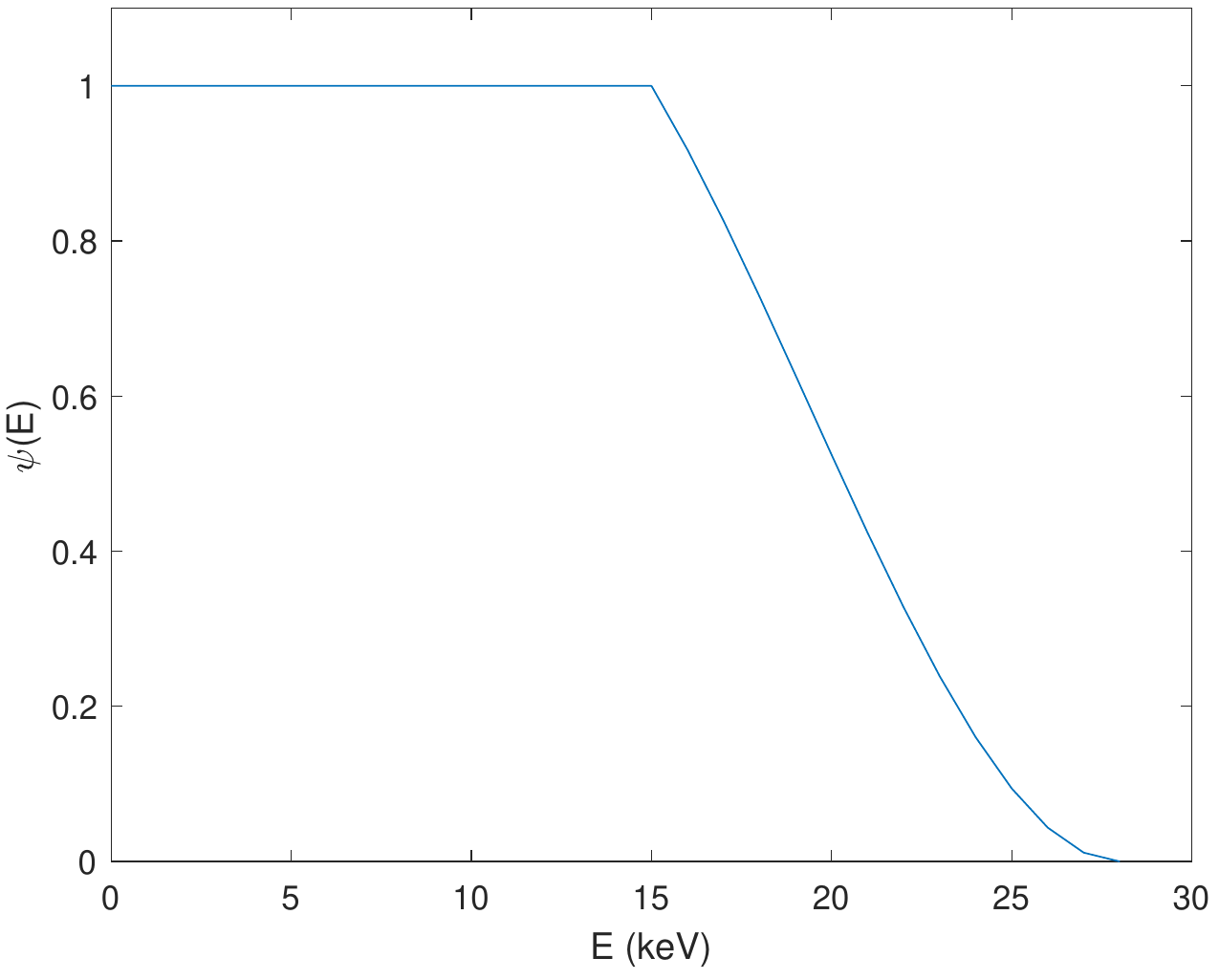} 
\caption{Plot of $\psi$, as in equation \eqref{psi}, with $E_M=29$keV and $E_m=1$keV.} \label{psiplot}
\end{figure}

In our case the Bragg data $\mathfrak{B}_af\in L^2([E_m,E_M]\times\Omega\times\Phi(\mathcal{I}))$, where $\Omega=[-w_{x_1},w_{x_1}]^2$. Thus, the boundary of the sinogram in our case is a 3-D cube in energy, source position and detector position, for each $x_2$ considered. The full data consists of all $(E,s_1,d_1)\in\mathbb{R}_+\times \mathbb{R}^2$. Due to finite scanner width and limitations on the energy range however, we consider the limited data on $[E_m,E_M]\times\Omega$. To deal with the cutoff at the boundary of $[E_m,E_M]\times\Omega$ we construct a filtering function $\psi\in C^{\infty}\paren{[E_m,E_M]\times\mathbb{R}^2\times\Phi(\mathcal{I})}$ which goes to zero smoothly in the energy variable $E$, as $E\to E_M$. In the remaining variables $s_1$ and $d_1$, we see a natural smoothing effect in the data towards zero for $s_1, d_1$ near the $x_1$ limit of the portal scanner (i.e. at $x_1=\pm w_{x_1}$) due to solid angle effects and the reduction in source intensity with increasing $x_1$ (see \eqref{I0}). That is, with a wide enough scanner (or large enough $w_{x_1}$), the detectors with $d_1\to\pm w_{x_1}$ have a more restrictive view of the scatter than those $d_1\approx 0$, and thus measure less scatter due to lower solid angle. The source projections with $s_1\to\pm w_{x_1}$ produce less counts due to reduced illumination of the object (i.e. at the scanner edge much of the source fan-beam does not intersect the scanning region), and decreased input intensity at further distances from the source. Also, as $E\to E_m$, we see a natural smoothing of the signal towards zero due to increased effects due to attenuation and self-absorption at low $E$. The self-absorbed photons are those which are absorbed by a material (due to the photoelectric effect) after being scattered initially from the same material. While the natural smoothing does not reach zero exactly at the boundary (e.g. at $s_1, d_1= \pm w_{x_1}$), we found the smoothing sufficient to reduce significantly the presence of boundary artifacts in the reconstruction. Therefore, we require no additional smoothing in the data in the $s_1$ and $d_1$ variables, and in $E$ as $E\to E_m$. 

With this in mind, we define $\psi$ as the third order polynomial in $E$
\begin{equation}
\label{psi}
\psi(E,s_1,d_1,\epsilon)=-\paren{\frac{2}{E_M}\paren{E_M-E}}^3+2\paren{\frac{2}{E_M}\paren{E_M-E}}^2
\end{equation}
for $E\in [\frac{E_M}{2},E_M]$, and $\psi(E)=1$ on $[E_m,\frac{E_M}{2}]$. See figure \ref{psiplot} for a plot of $\psi$ with $E_M=29$keV and $E_m=1$keV. We assume that $E_m<\frac{E_M}{2}$, and this assumption is satisfied in the simulations conducted in section \ref{results}. We choose to smooth $50\%$ of the data here (i.e. for $E\in [\frac{E_M}{2},E_M]$), as $50\%$ sinogram smoothing proved effective as a cutoff in \cite{frikel2013characterization,borg2018analyzing}. \begin{figure}[!h]
\centering
\begin{subfigure}{0.35\textwidth}
\includegraphics[width=1\linewidth, height=4.8cm, keepaspectratio]{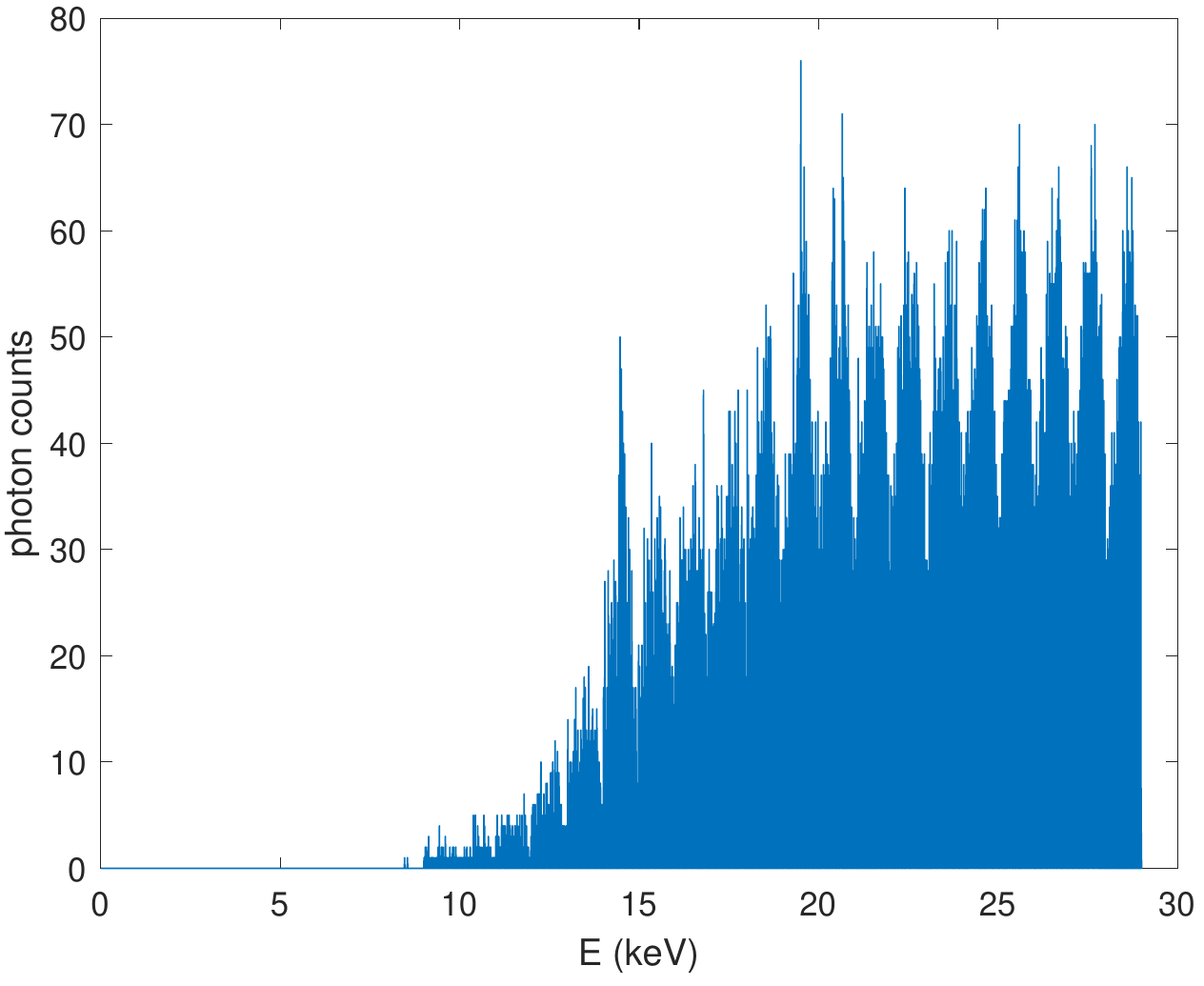}
\subcaption{$\mathfrak{B}_af$}\label{MFa}
\end{subfigure}
\begin{subfigure}{0.35\textwidth}
\includegraphics[width=1\linewidth, height=4.8cm, keepaspectratio]{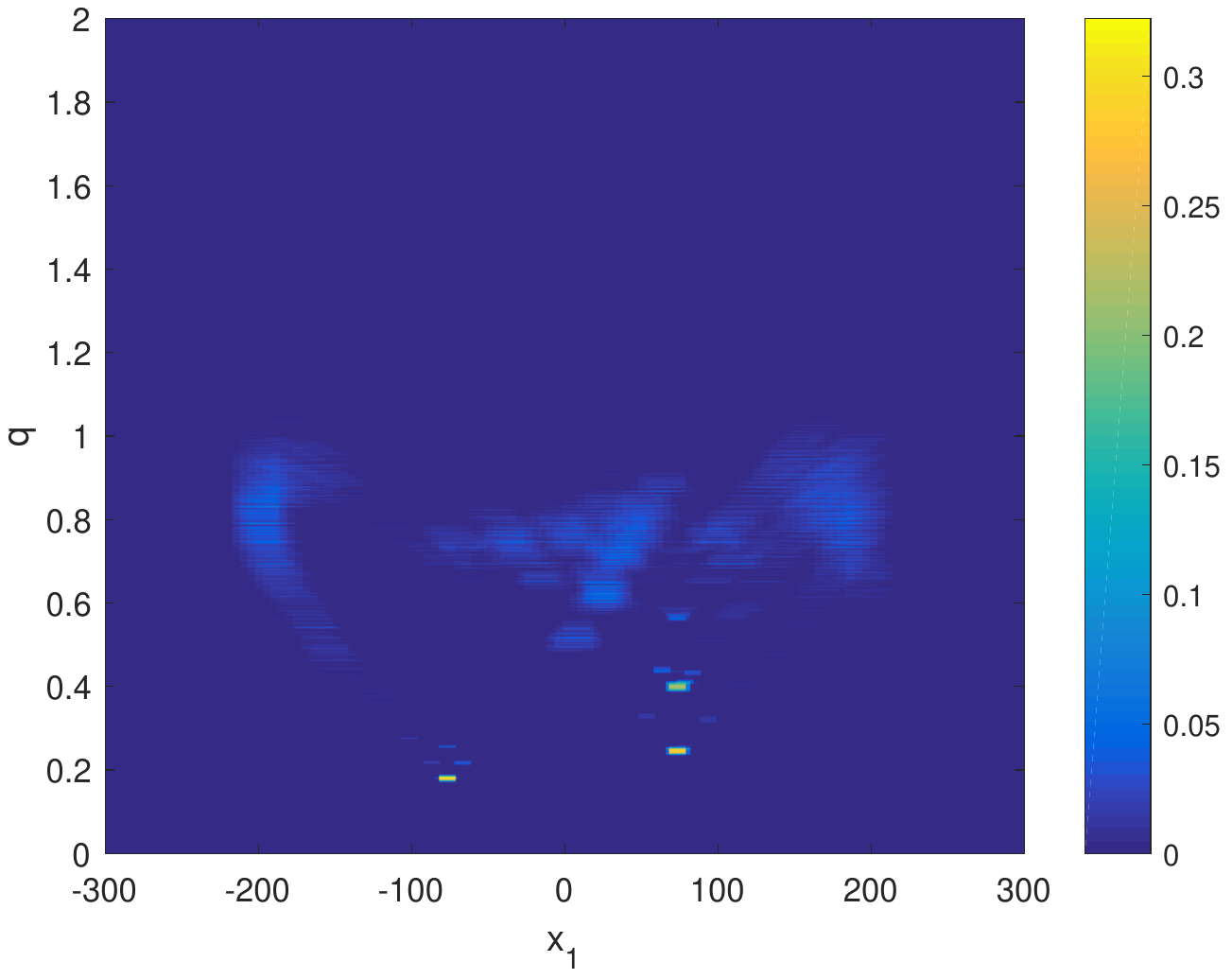}
\subcaption{$\mathfrak{B}_af$ reconstruction}\label{MFc}
\end{subfigure}
\begin{subfigure}{0.35\textwidth}
\includegraphics[width=1\linewidth, height=4.8cm, keepaspectratio]{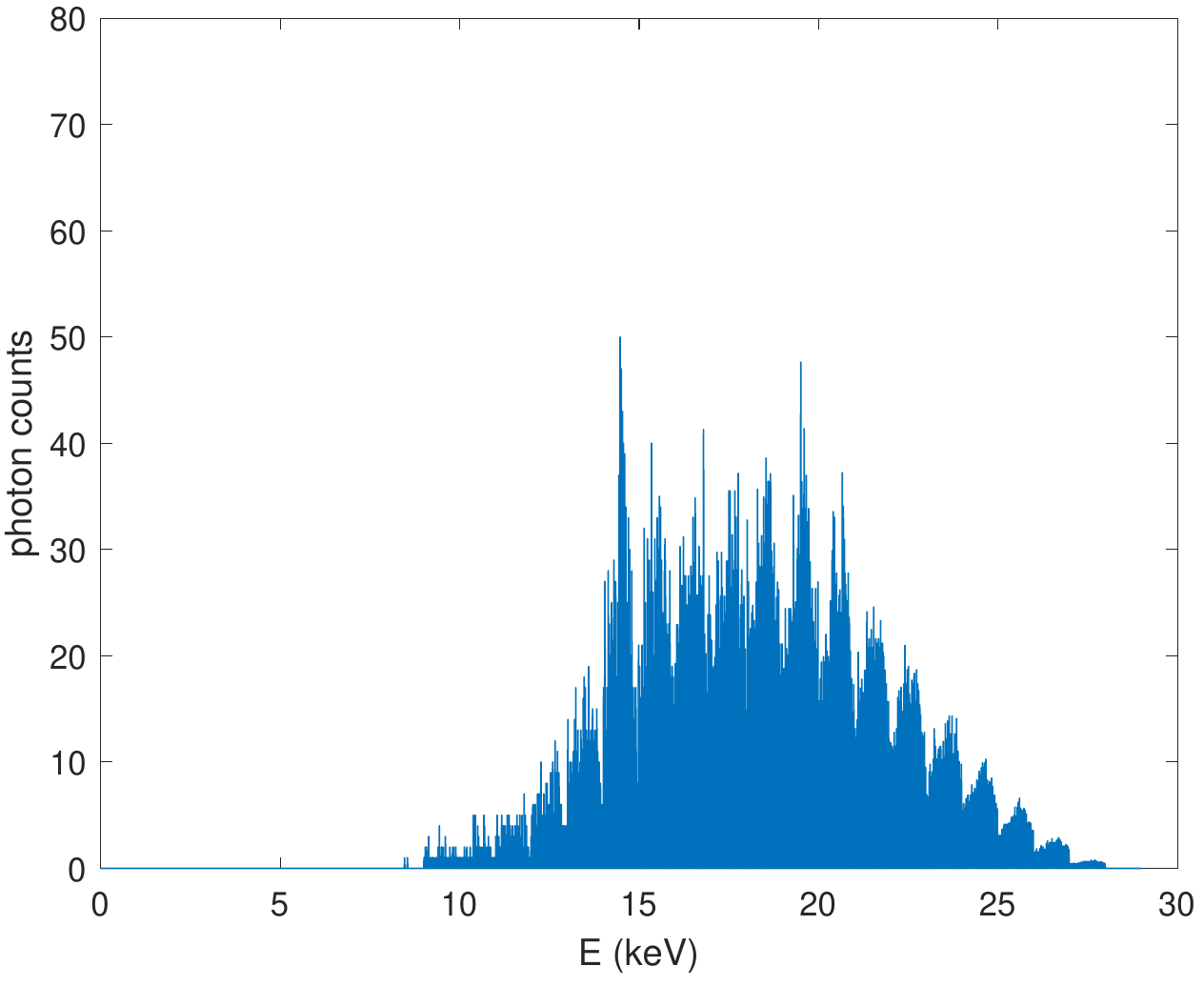} 
\subcaption{$\psi\mathfrak{B}_af$}\label{MFb}
\end{subfigure}
\begin{subfigure}{0.35\textwidth}
\includegraphics[width=1\linewidth, height=4.8cm, keepaspectratio]{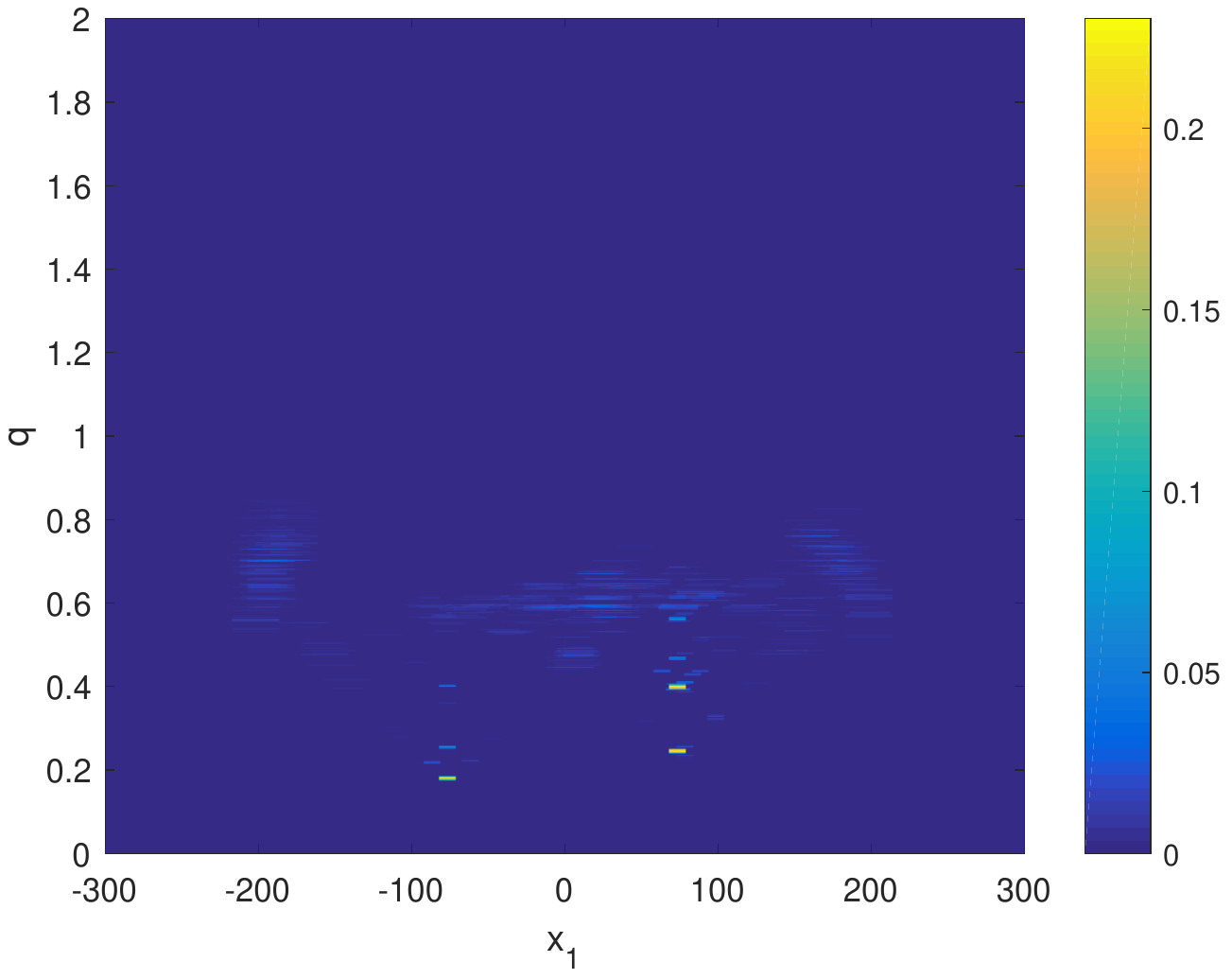}
\subcaption{$\psi\mathfrak{B}_af$ reconstruction}\label{MFd}
\end{subfigure}
\caption{Microlocal filtering examples.}
\label{MF1}
\end{figure}
Finally, to obtain our solution, we discretize $\psi\mathfrak{B}_af$ into a data vector $\vb$, and use $\vb$ as input to Algorithm \ref{algorithm1}. 

To show the effectiveness of this idea, and how it works, refer to figure \ref{MF1}. Figure \ref{MFa} shows a plot of Monte Carlo data $\mathfrak{B}_af$ of the scatter measured by the portal scanner from NaCl and Carbon (diamond structure) spheres. In this case $x_2=0$ and the spheres are centered along the line $\{x_3=0,x_2=0\}$. Figure \ref{MFb} shows the same data filtered by $\psi$. To clarify, the plots in figures \ref{MFa} and \ref{MFb} are the vectorized forms of the full 3-D data cube (in $(E,s_1,d_1)$), where each energy bin on the plot contains the corresponding vectorized 2-D datasets in $(s_1,d_1)$. The $s_1$ and $d_1$ positions used are specified later in section \ref{IP}. We present the plots in this way purely to show the effects of the $\psi$ filtering. We can see that the plot of figure \ref{MFb} is that of figure \ref{MFa} point multiplied by $\psi$, as in figure \ref{psiplot}. 

The reconstructions of $f(q,x_1,0)$ corresponding to figures \ref{MFa} and \ref{MFb}are shown in figures \ref{MFc} and \ref{MFd}, using $\mathfrak{B}_af$ and $\psi\mathfrak{B}_af$ as input to Algorithm \ref{algorithm1}, respectively, with $n_{\text{outiter}}=1$ outer iterations performed. Later in section \ref{results} we will show the results with higher $n_{\text{outiter}}$, and explain the imaging parameters in more detail. This example is included purely to illustrate the effectiveness of microlocal filters for our problem. In figure \ref{MFc} we see significant image artifacts, which roughly follow the shape of the Bragg integration curves of figure \ref{BC1}, in particular those along the top row. The curves along the top row of figure \ref{BC1} correspond to the case when $E=E_M$, that is, at the sinogram boundary in $E$. Thus we see artifacts in the reconstruction due to the sharp cutoff in sinogram at $E=E_M$. When the cutoff is smoothed by $\psi$, using $\psi\mathfrak{B}_af$ as input to Algorithm \ref{algorithm1} in figure \ref{MFd}, we notice a significant reduction in the artifacts, and a marked improvement in image quality.

\section{Testing and results}
\label{results}
In this section we perform the analytic and Monte Carlo data testing of Algorithm \ref{algorithm1} in a variety of imaging scenarios of interest in airport baggage screening and threat detection. First we discuss the imaging specifics (i.e. source/detector positions, source energies etc.) and then give simulated reconstructions from analytic and Monte Carlo data.

\subsection{The imaging parameters}
\label{IP}
Let us consider the scanning geometry depicted in figure \ref{figmain}. Throughout the simulations conducted in this section we convert to a more practical measurement system, using units of millimeters (mm), and set $w_{x_1}=300$mm and $w_{x_2}=410$mm. $\Phi$ is given explicitly by the linear map
\begin{equation}
\label{phi}
\Phi(x_2)=\frac{75}{820}(410-x_2),
\end{equation}
with $\Phi(-w_{x_2})=75$mm, $\Phi(w_{x_2})=0$mm and $\Phi(0)=37.5$mm. The chosen measurement system, $w_{x_1}$, $w_{x_2}$, and $\Phi$ are based on preliminary design specifics for the scanner of figure \ref{figmain}. For this study we use 31 source positions $s_1\in\{-300+20(j-1) : 1\leq j\leq 31\}$mm spaced at 20mm intervals, with source opening angle $\beta=120^{\circ}$. The sources are polychromatic (e.g. an X-ray tube) fan-beam, and we consider 29 spectrum energies $E\in\{1,\ldots,29\}$keV, with energy bin width 1keV. That is, the source spectrum is modeled as the sum of 29 delta distributions with centers $E\in\{1,\ldots,29\}$keV. A more accurate model would employ a fully continuous spectra and take averages over each energy bin. However, we leave the considerations of such effects to future work, and focus on the errors due to attenuation and Compton scatter (which are far more significant sources of error than energy bin averaging) here specifically.

As discussed in section \ref{reconmethod}, we aim to recover $f(q,x_1)$ on the range $[0,2]\AA^{-1}\times [-300,300]\text{mm}$. The maximum energy $E_M=29$keV is chosen to correspond with the maximum $q$ in the reconstruction space, namely $q_M=2$. That is
$$E_M=\frac{(hc)q_{M}}{\sin\paren{\frac{\omega_{M}}{2}}},\ \ \ \text{(see equation \eqref{qdef})}$$
where $\omega_{M}= 114.6^{\circ}$ is an upper bound on the scattering angles possible with the portal geometry of figure \ref{figmain}.

The input distribution $I_0(E)\propto\sum_{j=1}^{29}\delta(E-j)$ is chosen to be uniform in this study. This choice is well founded since, given the separability of $I_0(E,x_1)$ (equation \eqref{I0}), we can divide through by $I_0(E)>0$ in equation \eqref{equBG} to obtain the Bragg signal with uniform input spectra. We use 600 detector positions $d_1\in\{-300+(j-1) : 1\leq j\leq 600\}$mm spaced at 1mm intervals. A 1mm detector pitch is realistic given the pixel sizes of the energy sensitive detectors on the market \cite{egan2014dark} (the detectors in that paper have a 250$\mu$m pitch, and energy resolution $0.6$keV). The $\epsilon$ coordinate of the detectors varies with the chosen $x_2\in [-410,410]$mm, and is determined by $\epsilon=\Phi(x_2)$ (equation \ref{phi}). The detectors are assumed to be energy-resolved with energy resolution 1keV, so we are able to distinguish between all 29 energies in our range $E\in[1,29]$keV. The number of source positions (31) is low here as we anticipate longer scanning times per source projection so as to allow for sufficient photon counts. The increased number of detector positions (600) and energies bins (29) are at no detriment to the scan time. Note that due to the delta-comb nature of $I(E)$, there is no overlap in the energy bins. In total, the number of data points (or the number of rows of $A$) is $p=31\times 29\times 600 =539400$. We sample $f$ as a $750\times 600$ (high resolution) image, with $n=750$ $q$ samples in $\{q=\frac{2(j-1)}{750} : 1\leq j\leq 750\}\AA^{-1}$ and $m=600$ $x_1$ samples in $x_1\in\{-300+(j-1) : 1\leq j\leq 600\}$mm. We aim reconstruct $f$ on the box $[0,2]\AA^{-1}\times [-300,300]$mm. The length of $\vy$ (i.e. the number of reconstructed pixels) is therefore $n\times m=750\times 600=45000$ (also the number of columns of $A$).

We use the characteristic library consisting of $l=405$ elements, with widths $w_j\in [1,3]$cm at $5$mm intervals, and centers $x^c_j\in[-200,200]$mm at 5mm intervals. We are assuming $f(q,x_1)=0$ for $|x_1|>200$mm, i.e., $f$ is zero outside of the scanning tunnel depicted in figure \ref{figmain}.  We also assume that $f$ is composed of crystallite samples with width not exceeding 3cm, e.g., this could be a small sample of narcotics hidden amongst clutter in a mail parcel.

\subsection{Quantitative analysis and hyperparameter selection}
As a quantitative measure of the accuracy of our results we use the edge $F_1$ score \cite{webber2020joint,andrade2019shearlets,taha2015metrics}. Let $f$ and $\tilde{f}$ denote ground truth and reconstructed  images respectively. We employ the code of \cite{canny1986computational} here, which detects large jumps in the gradient images (i.e. in $\nabla f$ and $\nabla \tilde{f}$). The edge $F_1$ score is a measure of how well we have recovered $\nabla f$ and the locations of the Bragg peaks (i.e. the large spikes in the Bragg spectra of figure \ref{Fq1}). Let TP, FN and FP denote the proportion of true positives, false negatives and false positives in a classification result, respectively. Then the $F_1$ score is defined
\begin{equation}
\label{F1score}
F_1=\frac{2\text{TP}}{2\text{TP}+\text{FP}+\text{FN}}.
\end{equation}
Equivalently, the $F_1$ score is the geometric mean of the recall and precision of a classification result. The edge $F_1$ score is the $F_1$ score of \eqref{F1score} calculated between the edge maps $\nabla f$ and $\nabla\tilde{f}$. Specifically, we use the code of \cite{canny1986computational} (with high gradient threshold) to calculate $\nabla f$ and $\nabla\tilde{f}$ and convert to binary images. Then, we vectorize said binary images and use as input to equation \eqref{F1score}. 

Refer back to equation \ref{equ1} and the definitions of the hyperparameters $\alpha,\gamma,\lambda$. Through experimentation we found that the hyperparameter values $\alpha=1e+6$ (multibang parameter) and $\gamma=1e+10$ (no overlap constraint) worked well for the examples considered here. We choose the smoothing parameter $\lambda$ ($L^1$ smoothing) for the best results in terms of edge $F_1$ score.

\subsection{The materials considered}
Let $Z$ denote the atomic number of the material. Then we consider the $F(q,Z)$ curves for three materials here, namely Carbon with a graphite structure (denoted C-graphite), salt and Carbon with a diamond structure (denoted C-diamond). We choose these three materials as they fall into our $Z<20$ range, and exhibit a variety of crystalline structures (e.g. hexagonal for graphite, and face-center cubic for salt), and are thus suitable for testing our algorithm. Here $Z$ replaces $\vx\in\mathbb{R}^2$ from section \ref{braggtrans} to represent the material. Furthermore, the $F(q,Z)$ curves for C-graphite, salt and C-diamond are well known and readily available in the literature (e.g. see \cite{zhu2015theoretical} for the crystal structure of salt). Calculating $F(q,Z)$ for general compounds is a difficult task and there is work to be done in the spectroscopy literature to compile a wider database of $F$ curves. 

\begin{figure}[!h]
\centering
\begin{subfigure}{0.32\textwidth}
\includegraphics[width=0.9\linewidth, height=4cm, keepaspectratio]{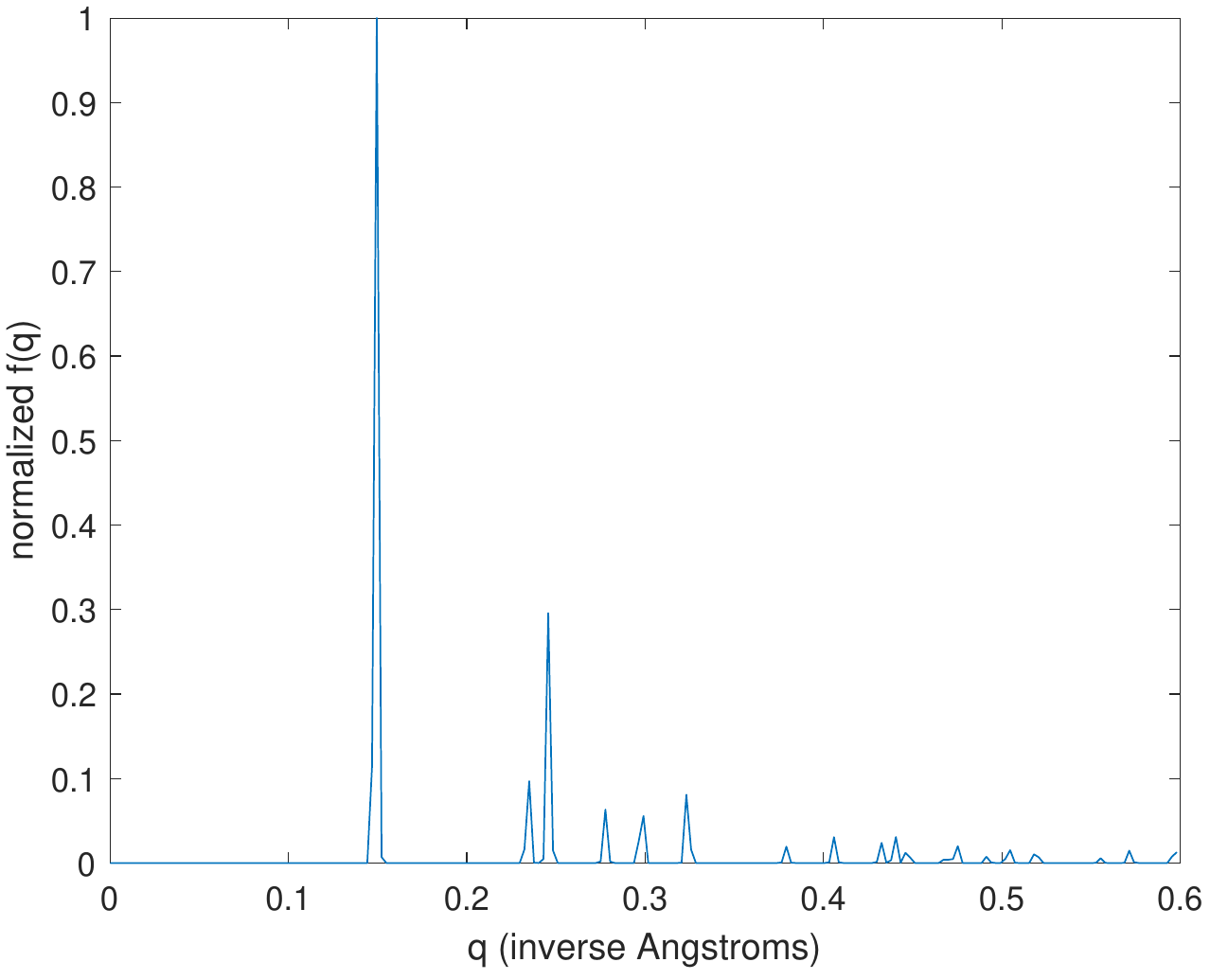}
\subcaption{C-graphite ($Z=6$)}
\end{subfigure}
\begin{subfigure}{0.32\textwidth}
\includegraphics[width=0.9\linewidth, height=4cm, keepaspectratio]{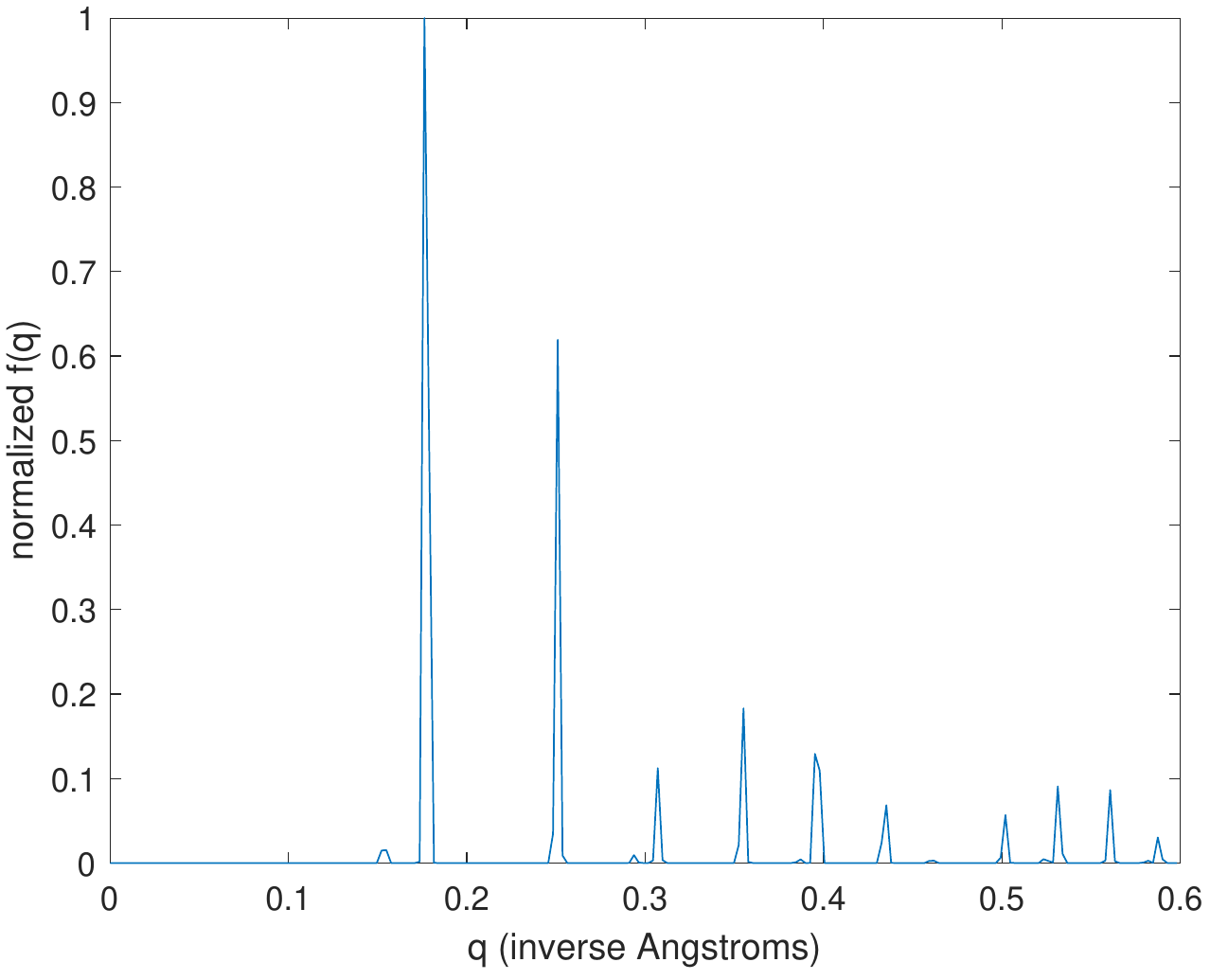} 
\subcaption{NaCl ($Z\approx 15$)}
\end{subfigure}
\begin{subfigure}{0.32\textwidth}
\includegraphics[width=0.9\linewidth, height=4cm, keepaspectratio]{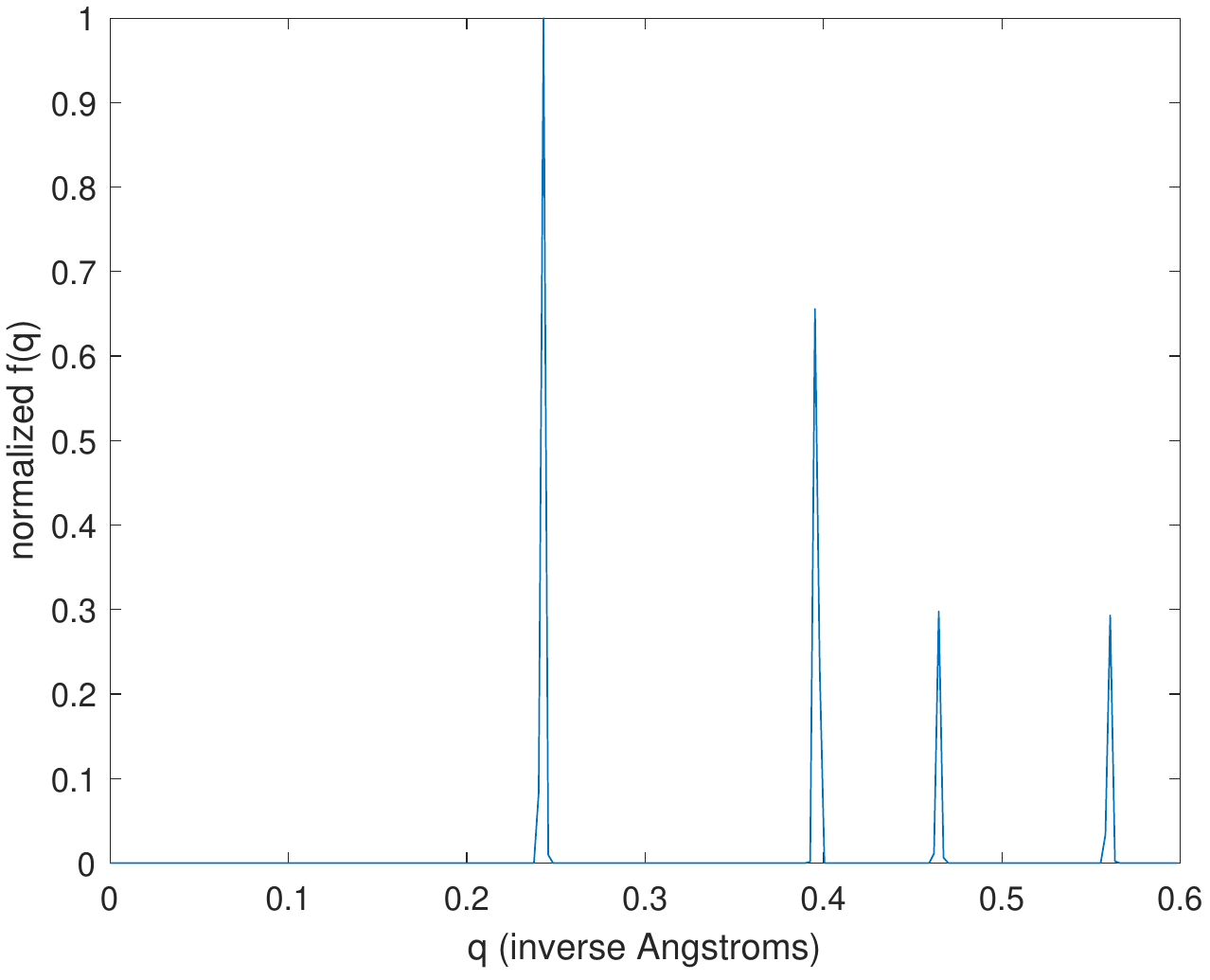}
\subcaption{C-diamond ($Z=6$)}
\end{subfigure}
\caption{$F(q)$ curve plots. The plots are normalized by $L^{\infty}$ norm (max value). We restrict the $q$ range to $[0,0.6]\AA^{-1}$ ($1\AA^{-1}\approx 12.4$keV) to better highlight the more significant (larger) Bragg peaks.}
\label{Fq1}
\end{figure}

The Bragg  differential cross section $F(q,Z)$ is defined \cite{webber2020bragg,ITC,paper}
\begin{equation}
\label{Bgm1}
\begin{split}
F(q,Z)\propto\frac{1}{q}\sum_{H}\delta\left(\frac{1}{2d_H}-q\right)d_H\left|F_{H}\left(q\right)\right|^2=g(q,Z),
\end{split}
\end{equation}
where $\delta$ is the Dirac-delta function and 
\begin{equation}
\label{FH}
F_H\left(q\right)=\sum_{i=1}^{n_a}F_i\left(q\right) e^{-2\pi i \vx_i\cdot H},
\end{equation}
is the scattering factor, where $n_a$ is the number of atoms in a crystal cell, the $\vx_i\in [0,1]^3$ are the coordinates of the atoms within the cell and $F_i$ is the atomic form factor \cite{hubbell,hubbell1} of atom $i$. For a given $Z$ let $Q=\text{supp}(g(\cdot,Z))\cap \{q<2\}=\{q_1,\ldots,q_{n_q}\}$ be set of $q$ values for which $g(\cdot,Z)$ is non-zero in the range $q\in[0,2]$, with $|Q|=n_q$. Then we model the $F$ curves in our simulations as the Gaussian mixture
\begin{equation}
\label{BgG}
F(q,Z)\propto\sum_{j=1}^{n_q}g(q_j,Z)e^{-\frac{(q-q_j)^2}{\sigma^2}},
\end{equation}
where $\sigma^2=10^{-6}$ is chosen to be small relative to $q_{\text{max}}$ so that the Gaussians of \eqref{BgG} accurately represent the delta functions of \eqref{Bgm1}. See figure \ref{Fq1} for plots of the $F$ curves for C-graphite, salt and C-diamond.

\subsection{Comparison to other methods}
For comparison, we consider a TV regularized solution, and use the microlocal filtering techniques explained in section \ref{microlocal} as data pre-processing to remove boundary artifacts. Specifically we aim to find
\begin{equation}
\label{equTV}
\argmin_{\vy}\sum_k\left[\paren{A\vy}_k -b_k\log\paren{A\vy}_k\right]+\lambda\text{TV}_{\beta}(\vy),\ \ \vy\in [0,\infty)^{n\times m}
\end{equation}
where
\begin{equation}
\text{TV}_{\beta}(\vy)=\sum_{j=1}^m\sum_{i=1}^n\paren{\paren{\nabla Y}_{ij}^2+\beta^2}^{\frac{1}{2}}
\end{equation}
is a smoothed TV norm, where $Y$ is $\vy$ reshaped into a $n\times m$ matrix (image), and $\nabla Y$ is the gradient image. The parameter $\beta>0$ is included so that the gradient of $\text{TV}_{\beta}$ is defined at zero. The smoothing parameter $\lambda$ controls the level of TV regularization. TV regularization is applied to good effect in the BST literature \cite{greenberg2013snapshot,webber2020bragg}, and thus why it is chosen as a point of comparison. We also combine TV ideas with those of microlocal analysis in section \ref{microlocal}, to assist in the reduction of boundary artifacts. To solve \ref{equTV} we implement the code ``JR\_PET\_TV" of \cite{ehrhardt2014joint} with non-negativity constraints, which is applied in that paper to PET imaging. That is, we input filtered data $\vb$ (using the filters of section \ref{microlocal} as pre-processing) to JR\_PET\_TV, constraining $\vy\in [0,\infty)^{n\times m}$, and choose $\lambda,\beta$ for the best results in terms of edge $F_1$ score. We denote this method as ``FTV", which stands for Filtered Total Variation.

\subsection{Analytic data testing}
\label{analytic}
\begin{figure}[!h]
\centering
\begin{subfigure}{0.4\textwidth}
\includegraphics[ width=1\linewidth, height=1\linewidth, keepaspectratio]{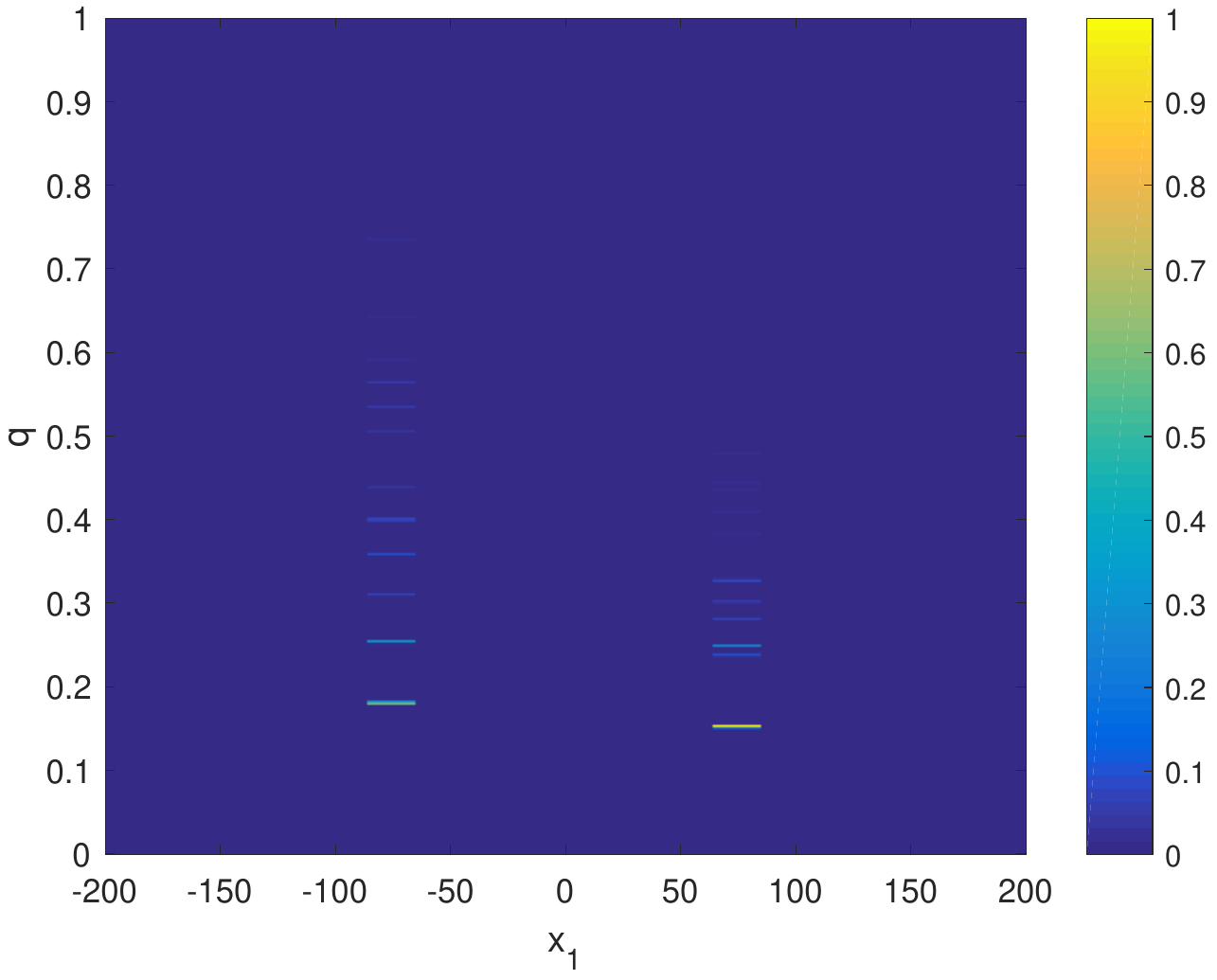} 
\subcaption{Phantom 1}
\end{subfigure}
\begin{subfigure}{0.4\textwidth}
\includegraphics[ width=1\linewidth, height=1\linewidth, keepaspectratio]{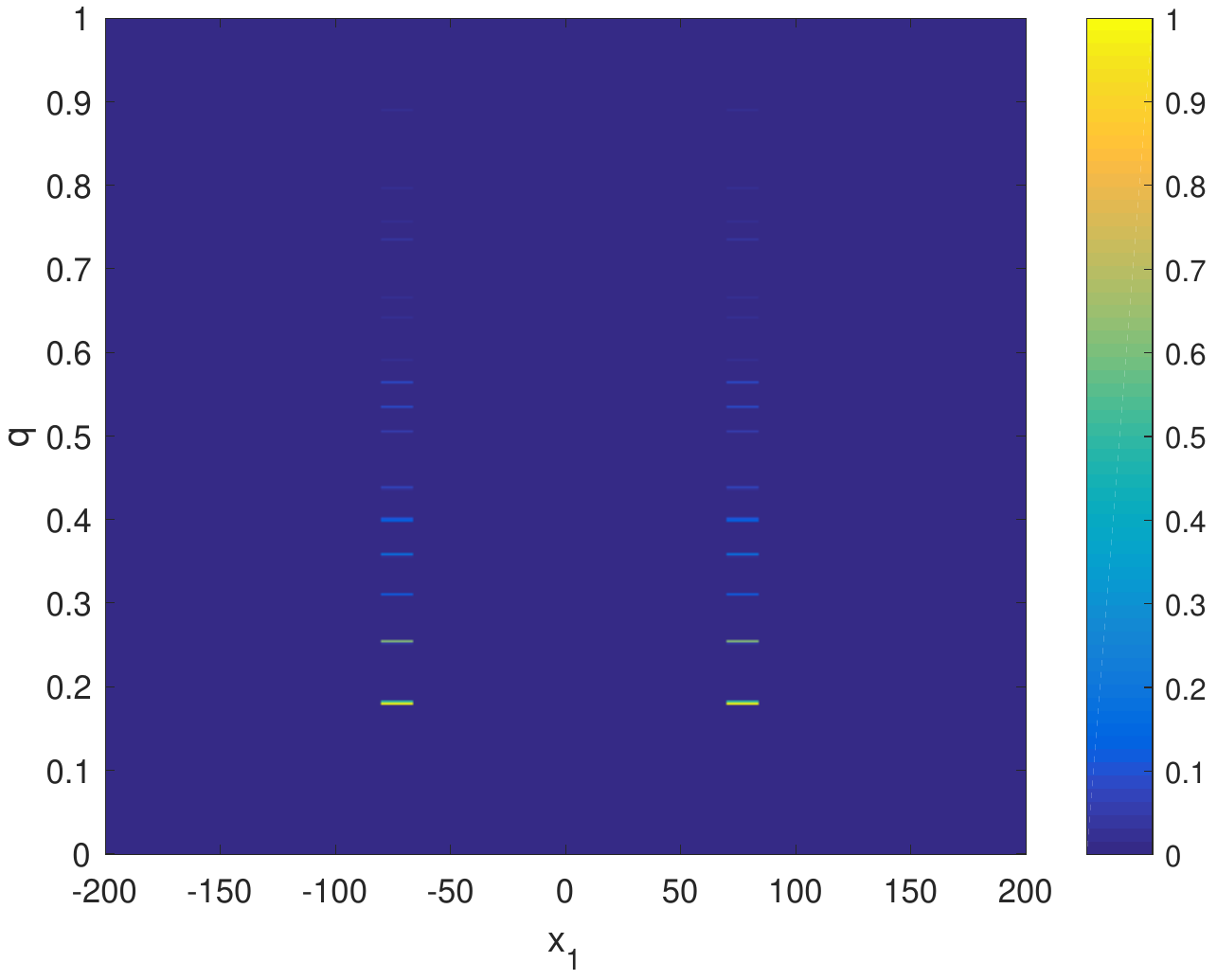} 
\subcaption{Phantom 2}
\end{subfigure}
\caption{Imaging phantoms.}
\label{figPhan}
\end{figure}

Here we consider the analytic data testing of Algorithm \ref{algorithm1} by means of a Poisson noise model. That is, we use a scaled, matched model as mean to a multivariate Poisson distribution to generate data, where the scaling  factor controls the level of noise. We consider two imaging phantoms for reconstruction. See figure \ref{figPhan}. The left-hand phantom is comprised of an NaCl and C-graphite object with centers $x_c=-75$mm and $x_c=75$mm respectively, both with width $w=20$mm. The widths and centers used are chosen to be part of the imaging basis used (as detailed in section \ref{IP}). The right-hand phantom is comprised of two NaCl objects with centers $x_c=-72.5$mm and $x_c=77.5$mm, both with width $w=1.25$cm. The widths and centers in this case are chosen to lie outside of the imaging basis. Let $\vy$ denote the vetorized form of $f(\cdot,\cdot,x_2)$ for a given $x_2\in[-410,410]$mm. Then the phantoms are normalized to have max value $\|\vy\|_{\infty}=1$. 

The analytic data is generated by
\begin{equation}
\vb_{\eta}\sim \text{Poisson}\paren{\eta_c l_{\vb}\frac{A\vy}{\|A\vy\|_1}},
\end{equation}
where $l_{\vb}$ is the length of $\vb= A\vy$. That is, the exact model $A\vy$ is normalized in $L^1$ (by total photon counts) and scaled by $\eta_cl_{\vb}$, then used as mean to a multivariate Poisson distribution.
\begin{table}[!h]
\begin{subtable}{.49\linewidth}\centering
{\begin{tabular}{| c | c | c | c | c | c | c | c |}
\hline
Edge $F_1$ score & 2DBSR  & FTV    \\ \hline
$x_2=0$     & $.86$ & $.74$  \\ 
$x_2=205$     & $.88$ & $.93$  \\  
$x_2=-205$     & $.88$ & $.93$  \\  \hline
\end{tabular}}
\caption{$\eta_c=10$ ($\eta_{\text{ls}}\approx 15\%$)}
\end{subtable}
\begin{subtable}{.49\linewidth}\centering
{\begin{tabular}{| c | c | c | c | c | c | c | c |}
\hline
Edge $F_1$ score & 2DBSR  & FTV    \\ \hline
$x_2=0$     & $.86$ & $.74$  \\ 
$x_2=205$     & $.88$ & $.73$  \\  
$x_2=-205$     & $.88$ & $.84$  \\  \hline
\end{tabular}}
\caption{$\eta_c=1$ ($\eta_{\text{ls}}\approx 50\%$)}
\end{subtable}
\caption{Phantom 1 results. The measurements in the left column are given in mm.}
\label{T1}
\end{table}
\begin{table}[!h]
\begin{subtable}{.49\linewidth}\centering
{\begin{tabular}{| c | c | c | c | c | c | c | c |}
\hline
Edge $F_1$ score & 2DBSR  & FTV    \\ \hline
$x_2=0$     & $.77$ & $.89$  \\ 
$x_2=205$     & $.79$ & $.89$  \\  
$x_2=-205$     & $.79$ & $.89$  \\  \hline
\end{tabular}}
\caption{$\eta_c=10$ ($\eta_{\text{ls}}\approx 15\%$)}
\end{subtable}
\begin{subtable}{.49\linewidth}\centering
{\begin{tabular}{| c | c | c | c | c | c | c | c |}
\hline
Edge $F_1$ score & 2DBSR  & FTV    \\ \hline
$x_2=0$     & $.78$ & $.89$  \\ 
$x_2=205$     & $.79$ & $.89$  \\  
$x_2=-205$     & $.79$ & $.89$  \\  \hline
\end{tabular}}
\caption{$\eta_c=1$ ($\eta_{\text{ls}}\approx 50\%$)}
\end{subtable}
\caption{Phantom 2 results. The measurements in the left column are given in mm.}
\label{T2}
\end{table}
\begin{figure}[!h]
\centering
\begin{subfigure}{0.4\textwidth}
\includegraphics[ width=1\linewidth, height=1\linewidth, keepaspectratio]{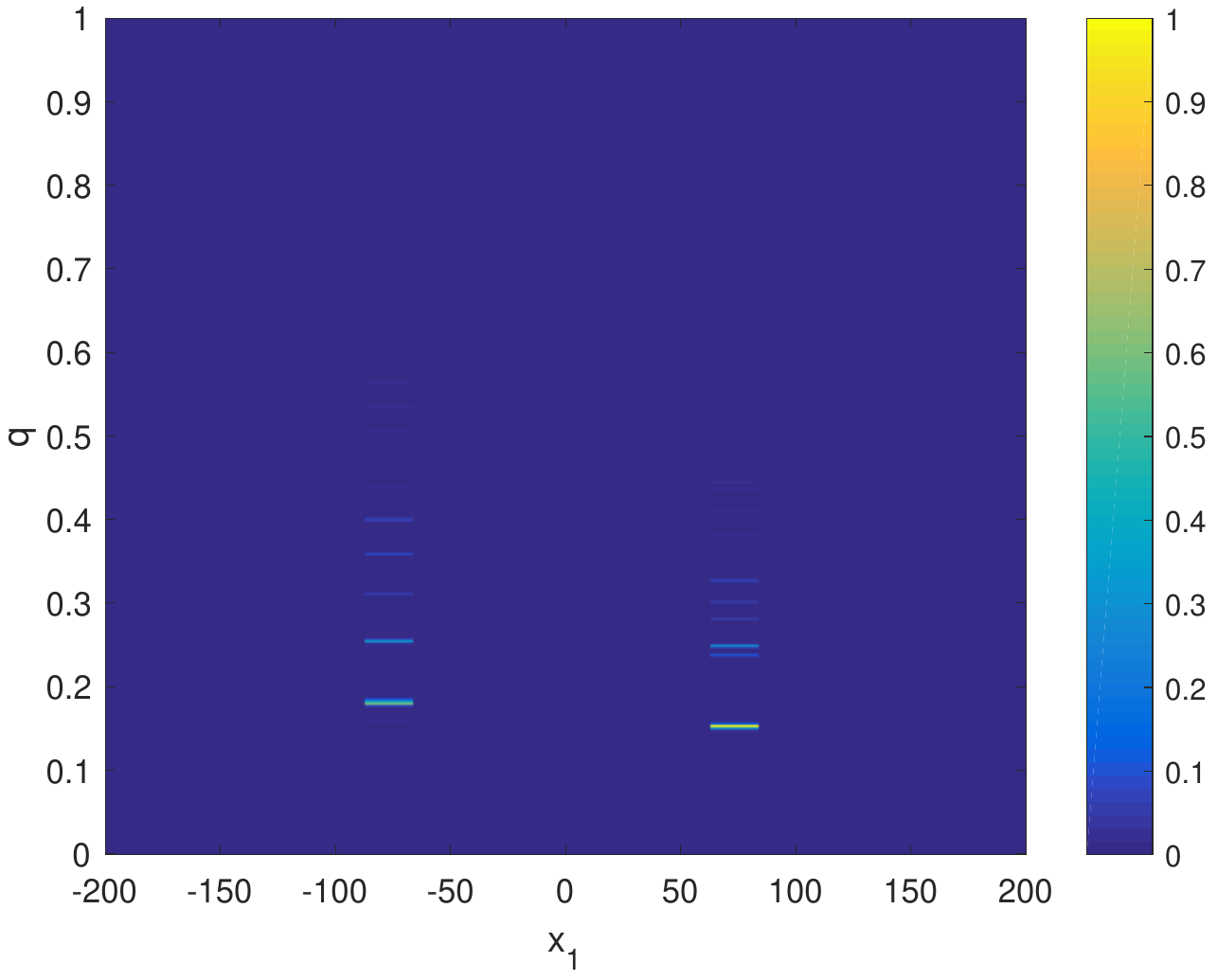} 
\subcaption{2DBSR}
\end{subfigure}
\begin{subfigure}{0.4\textwidth}
\includegraphics[ width=1\linewidth, height=1\linewidth, keepaspectratio]{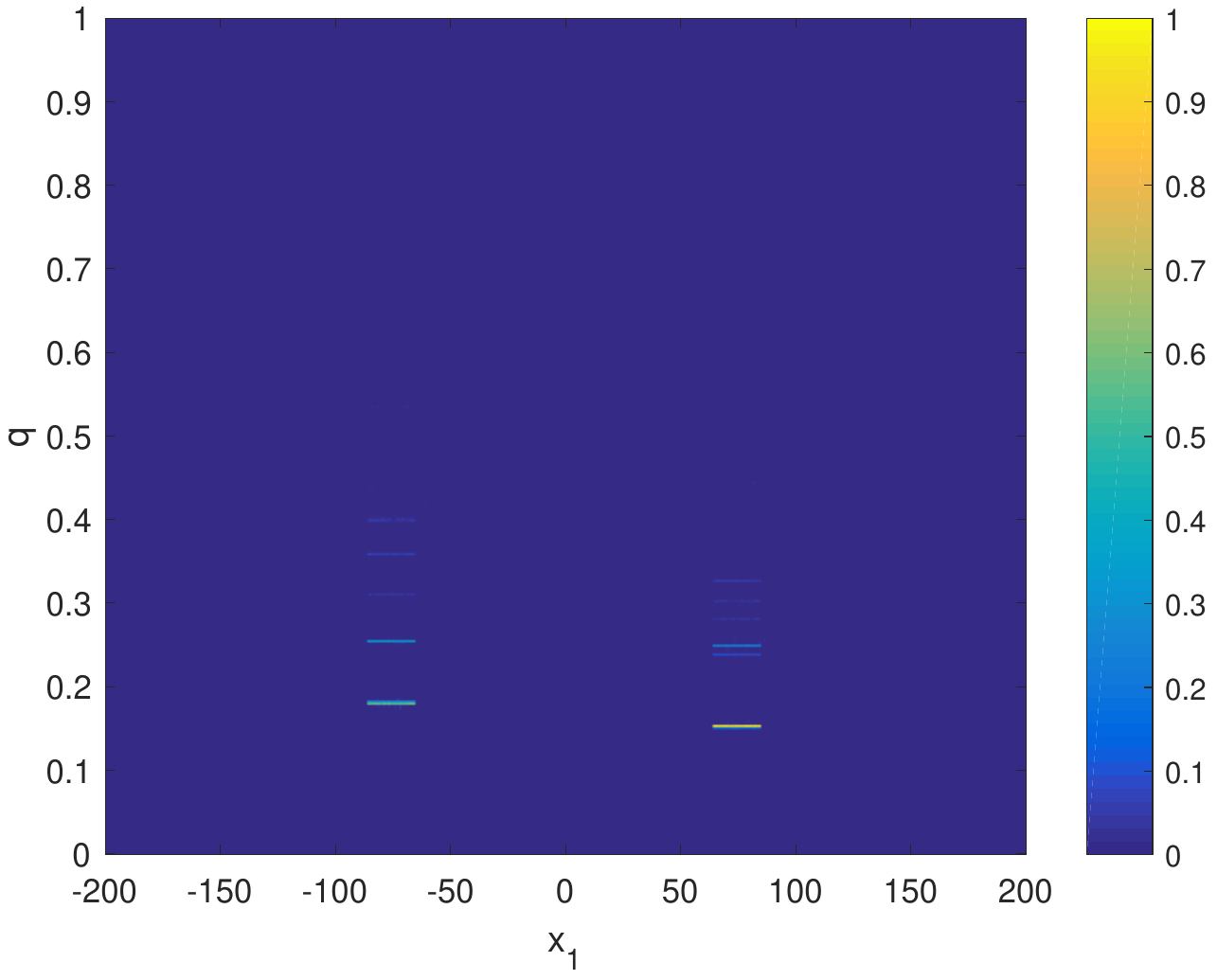} 
\subcaption{FTV}
\end{subfigure}
\caption{Phantom 1 reconstructions using 2DBSR (left) and FTV (right) along the central line profile ($x_2=0$), with count level $\eta_c=1$ ($50\%$ noise).}
\label{Phan1}
\end{figure}
\begin{figure}[!h]
\centering
\begin{subfigure}{0.4\textwidth}
\includegraphics[ width=1\linewidth, height=1\linewidth, keepaspectratio]{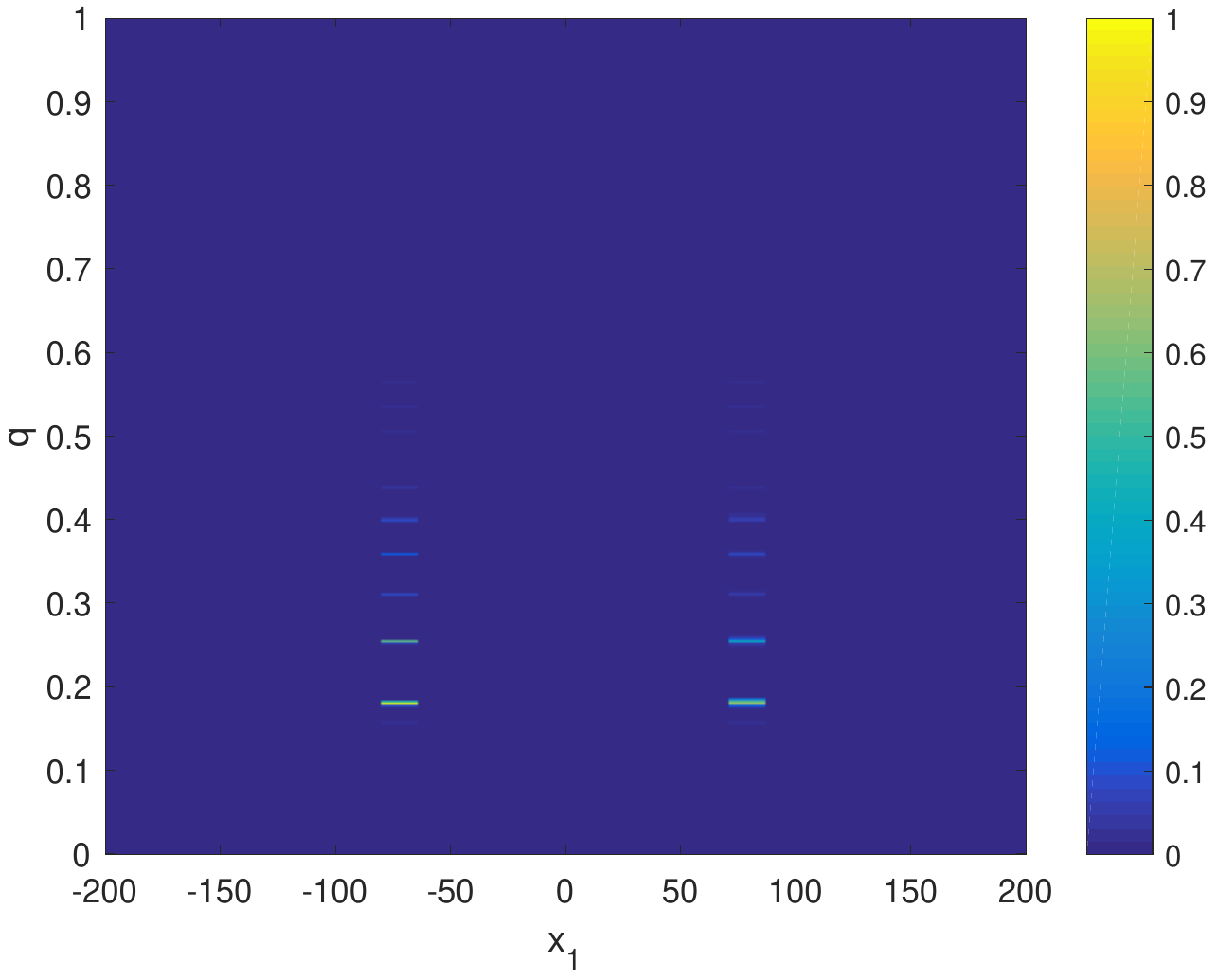} 
\subcaption{2DBSR}
\end{subfigure}
\begin{subfigure}{0.4\textwidth}
\includegraphics[ width=1\linewidth, height=1\linewidth, keepaspectratio]{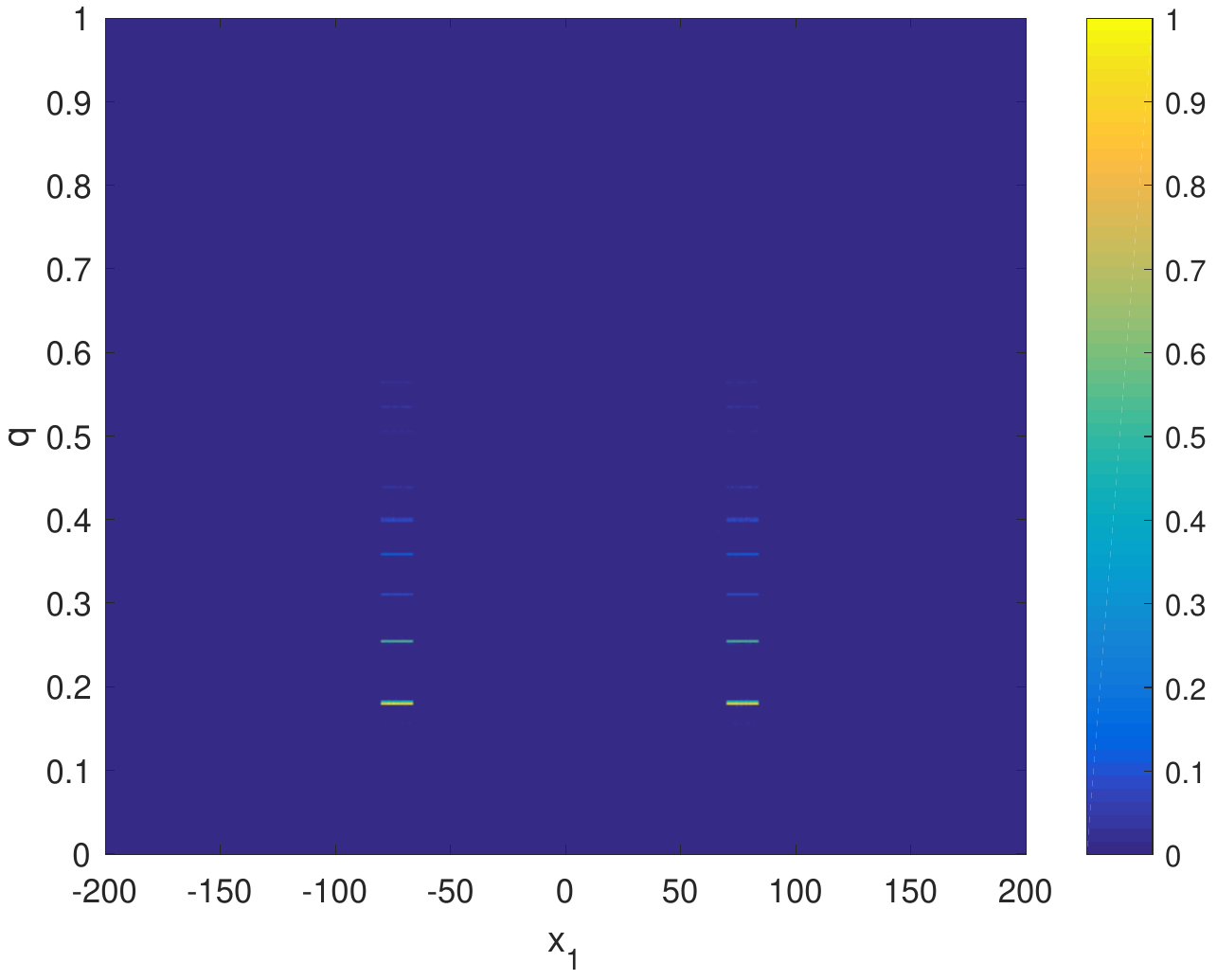} 
\subcaption{FTV}
\end{subfigure}
\caption{Phantom 2 reconstructions using 2DBSR (left) and FTV (right) along the central line profile ($x_2=0$), with count level $\eta_c=1$ ($50\%$ noise).}
\label{Phan2}
\end{figure}
The noisy data $\vb_{\eta}$ is generated as a Poisson draw. The scaling parameter $\eta_c$ is the average counts per detector and controls the level of noise, i.e., larger $\eta_c$ implies less noise and vise-versa. We consider two average count levels $\eta_c=10$ and $\eta_c =1$, which approximately equate to relative least square errors of $\eta_{\text{ls}}=15\%$ and $\eta_{\text{ls}}=50\%$ respectively, where $\eta_{\text{ls}}=\frac{\|\vb-\vb_{\eta}\|_2}{\|\vb\|_2}$. We consider three imaging line profiles $x_2\in [-410,410]$mm, namely $x_2=0$mm, $x_2=-205$mm, and $x_2=205$mm. Note that the operator $A$ varies with $x_2$, so we are considering three 2-D linear inverse problems with three different $A$ operators in total. 

See tables \ref{T1} and \ref{T2} for our results using 2DBSR (Algorithm \ref{algorithm1}) and FTV in terms of edge $F_1$ score. In figures \ref{Phan1} and \ref{Phan2} we present 2-D image reconstructions along the central line profile ($x_2=0$), for $\eta_c=1$. The image quality and $F_1$ scores are good and comparible using both methods. However we notice a significantly increased $F_1$ score using FTV on phantom 2. This is because the phantom 2 materials were chosen to lie outside imaging basis used for 2DBSR. In this case the 2DBSR algorithm chooses the centers and widths closest to the ground truth, namely $x_c=\pm75$mm and $w=1.5$cm. Thus in the 2DBSR reconstruction the centers and widths are slightly shifted from the ground truth causing a reduction in $F_1$ score. We could increase the size of the characteristic library to combat this, although this would be at a cost to the level of regularization and the implementation speed.  When using FTV, we have no such restrictions on imaging basis.

\subsection{Monte Carlo experiments}
\label{MC}
\begin{figure}[h!]
\centering
\includegraphics[ width=0.65\linewidth, height=0.65\linewidth, keepaspectratio]{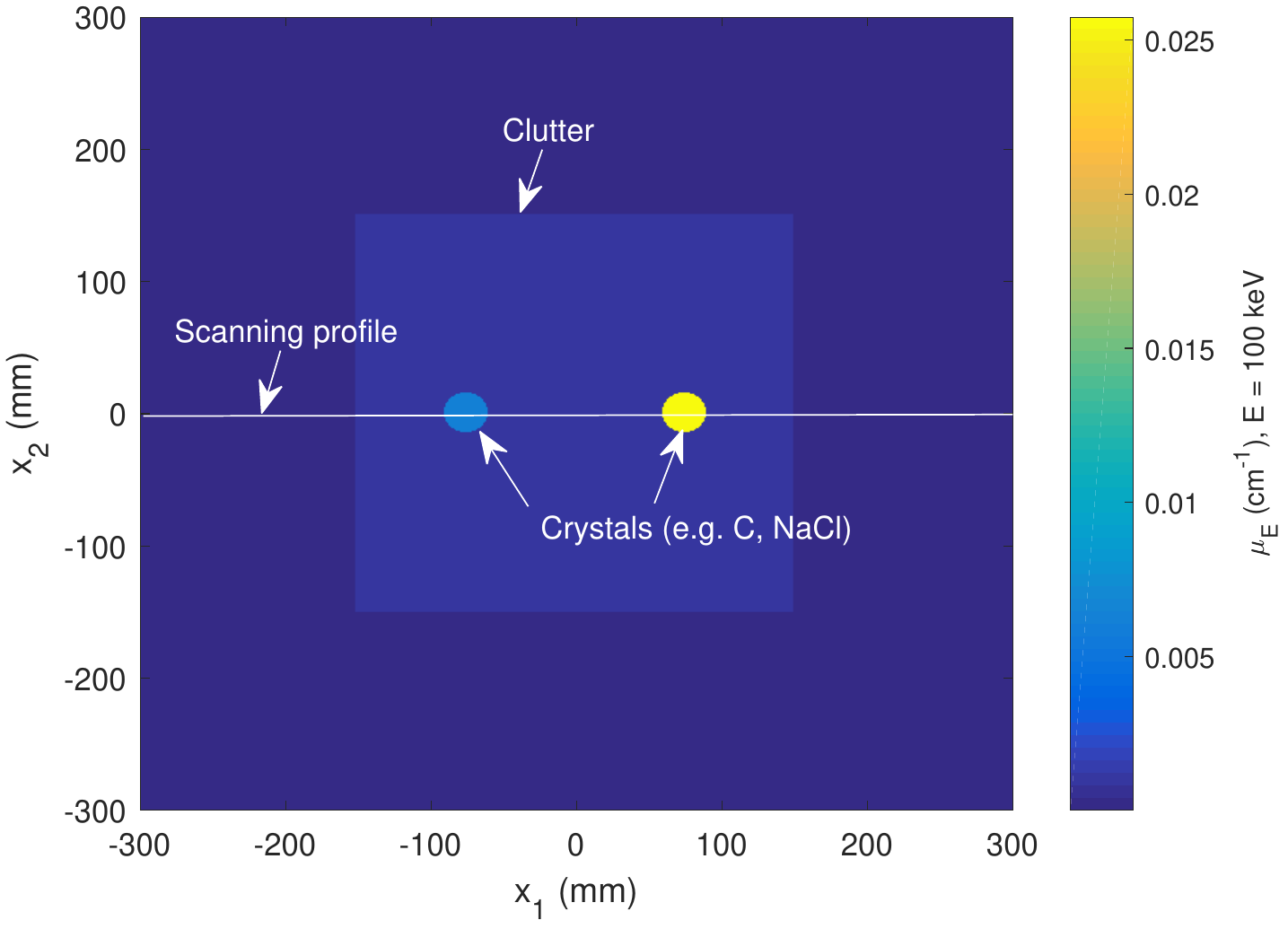} 
\caption{Test parcel.} \label{Parcel}
\end{figure}

In the Monte Carlo tests conducted in this section, we aim to recover spherical crystalline objects which are centered on the central ($x_2=0$) line profile.  Consider the schematic shown in figure \ref{Parcel}. Here we have displayed the attenuation image ($\mu_E$ for $E=100$keV) of two spherical crystallites centered on the $x_2=0$ line profile, which are embedded in clutter (e.g. paper, clothing). The experimental setup of figure \ref{Parcel} is chosen to simulate a parcel in commerce mail or a carry-on bag at an airport, whereby the sender (or passenger) has concealed a small sample of crystalline material in amongst a large volume of clutter. This could be intentional (e.g. attempted concealment and smuggling of narcotics) or accidental concealment (e.g. cluttered airport baggage with over-the-counter medicine inside).

We use the imaging parameters (e.g. source/detector positions, source energies) specified in section \ref{IP}. The $\epsilon$ coordinate of the detector array is determined by $\Phi$ (equation \ref{phi}), and is fixed at $\epsilon_0=\Phi(0)=37.5$mm throughout this section. $\epsilon_0$ is the $\epsilon$ coordinate of the detector array which corresponds to $x_2=0$. As discussed in section \ref{IP}, we use 29 source energies $E\in [1,29]$keV at 1keV intervals. The source intensity used here is $I_0(E)=6\times 10^{10}$ counts per projection, per energy, which amounts to $1.74\times 10^{12}$ counts per projection over all 29 energy bins. Equivalently the source intensity is $2.94\times 10^9$ counts, per cm$^2$, per energy bin, at 1m from the source. 
\begin{figure}[!h]
\centering
\begin{subfigure}{0.4\textwidth}
\includegraphics[ width=1\linewidth, height=1\linewidth, keepaspectratio]{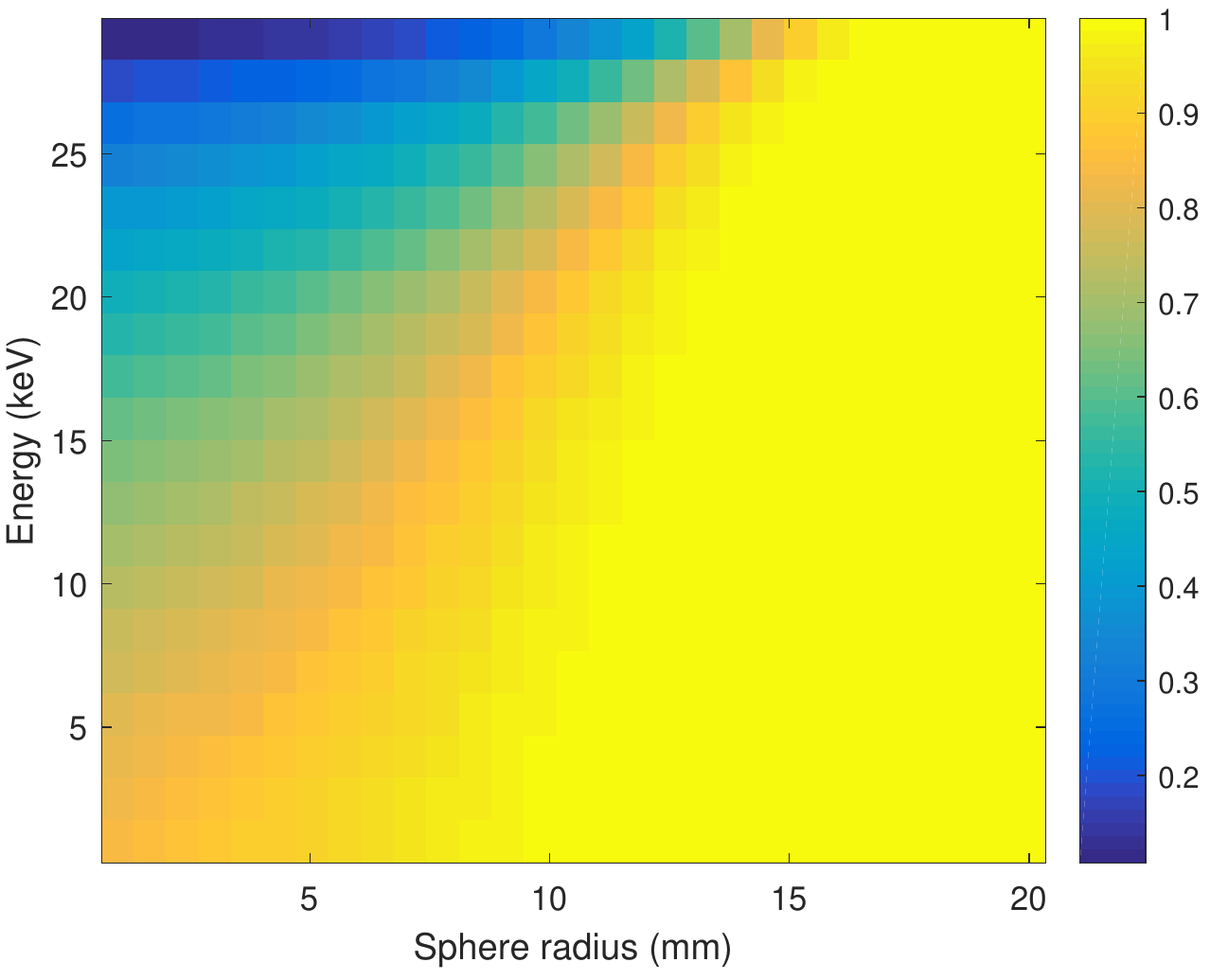} 
\subcaption{relative $\mu$ error}\label{muErr}
\end{subfigure}
\begin{subfigure}{0.4\textwidth}
\includegraphics[ width=1\linewidth, height=1\linewidth, keepaspectratio]{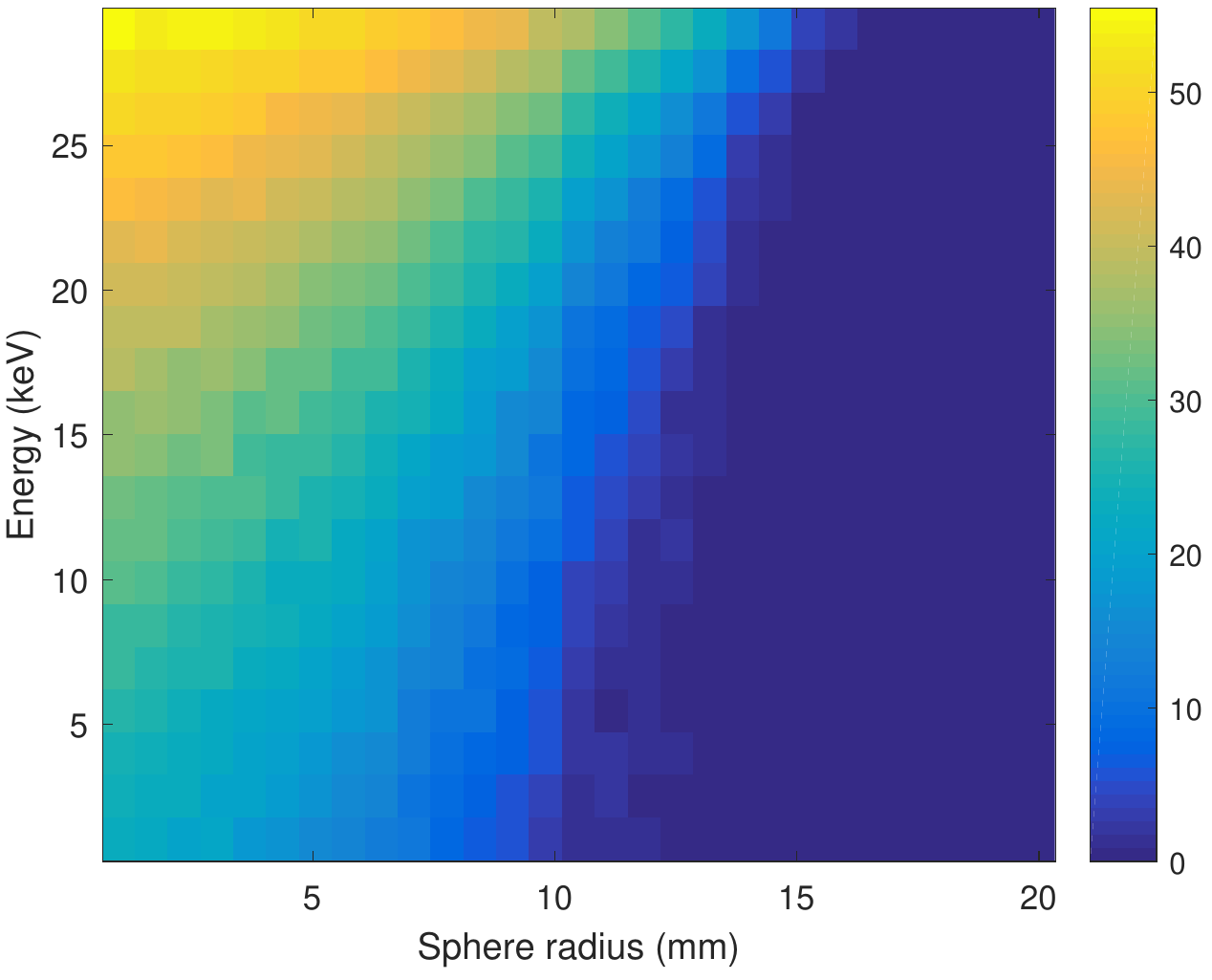} 
\subcaption{absolute Compton signal}\label{CompSig}
\end{subfigure}
\caption{Plots of relative, systematic attenuation error (left) and absolute Compton signal (right) with $E$ and sphere radius ($r$).}
\label{figERR}
\end{figure}

The clutter is simulated as a 30cm$^3$ box centered at the origin. The clutter attenuation coefficient is simulated as one-tenth that of cellulose. Cellulose is chosen as a material as it is a primary component of paper and clothing, which are commonly found in airport luggage and mail parcels. 
\begin{table}[!h]
\centering
\begin{tabular}{| c | c | c | c | c | c | c | c |}
\hline
Edge $F_1$ score & S1 & S2   & FTV    \\ \hline
St ($x_c=0, r=10$) & .82  & .83     & $.17$  \\ 
St ($x_c=0, r=15$) & .66  & .76  & $.08$  \\ 
\bf{C (\boldmath{$x_c=0, r=15$})} & \textbf{.64}  & \textbf{.73}      & \textbf{.05}\\ 
\bf{Di (\boldmath{$x_c=-6.25, r=8.75$})\boldmath{$^*$}} & \textbf{.60}  & \textbf{.71}      & \textbf{.42} \\ 
St ($x_c=-75, r=5$), Di ($x_c=75, r=5$) & .73  & .70    & $.54$  \\ 
\bf{St} (\boldmath{$x_c=-75, r=10$}), \bf{Di} (\boldmath{$x_c=75, r=10$}) & \textbf{.59}  & \textbf{.66}      & \bf{.10}  \\ 
\bf{St} (\boldmath{$x_c=-72.5, r=6.25$}), \bf{St} (\boldmath{$x_c=77.5, r=6.25$})\boldmath{$^*$} & \textbf{.57}  & \textbf{.52}      & \bf{.30}\\ 
\bf{St} (\boldmath{$x_c=-75, r=10$}), \bf{C} (\boldmath{$x_c=75, r=10$})& \textbf{.42}  & \textbf{.35}  & \bf{.05}  \\
C ($x_c=-75, r=10$), Di ($x_c=75, r=10$) & .60  & .67  & $.05$  \\ \hline
\end{tabular}
\caption{$F_1$ scores for each experiment conducted, comparing stage 1 (S1), and stage 2 (S2) of 2DBSR (ours) with FTV (literature). The measurements in the brackets in the left-hand column are in millimeters, and describe the centers ($x_c$) and radii ($r$) of the spheres scanned. Here ``St" denotes NaCl (salt), ``C" denotes C-graphite and ``Di" is C-diamond. The asterisked entries denote materials which are outside of the imaging basis.}
\label{T3}
\end{table}
The attenuation coefficient of the cellulose is scaled by $1/10$ following the assumption that the clutter will largely be comprised of air, thus reducing the attenuative effects. As cellulose is an amorphous material with significantly lower attenuation coefficient than that of the crystalline materials considered in this paper (i.e. C-graphite, NaCl, C-diamond), the scattering contribution from the clutter is modeled as zero.

\begin{figure}[!h]
\centering
\begin{subfigure}{0.32\textwidth}
\includegraphics[width=0.9\linewidth, height=4cm, keepaspectratio]{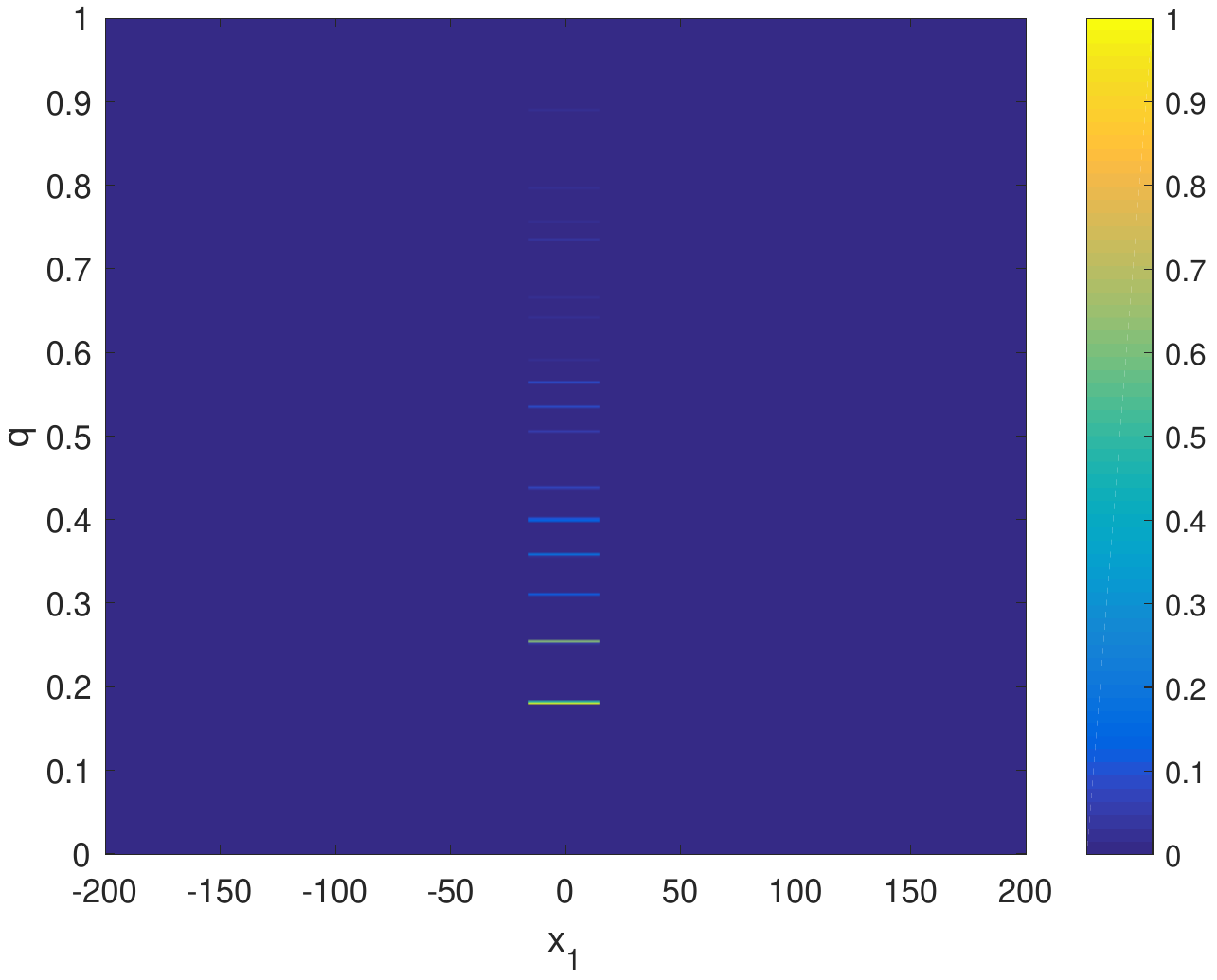}
\subcaption{ground thruth}
\end{subfigure}
\begin{subfigure}{0.32\textwidth}
\includegraphics[width=0.9\linewidth, height=4cm, keepaspectratio]{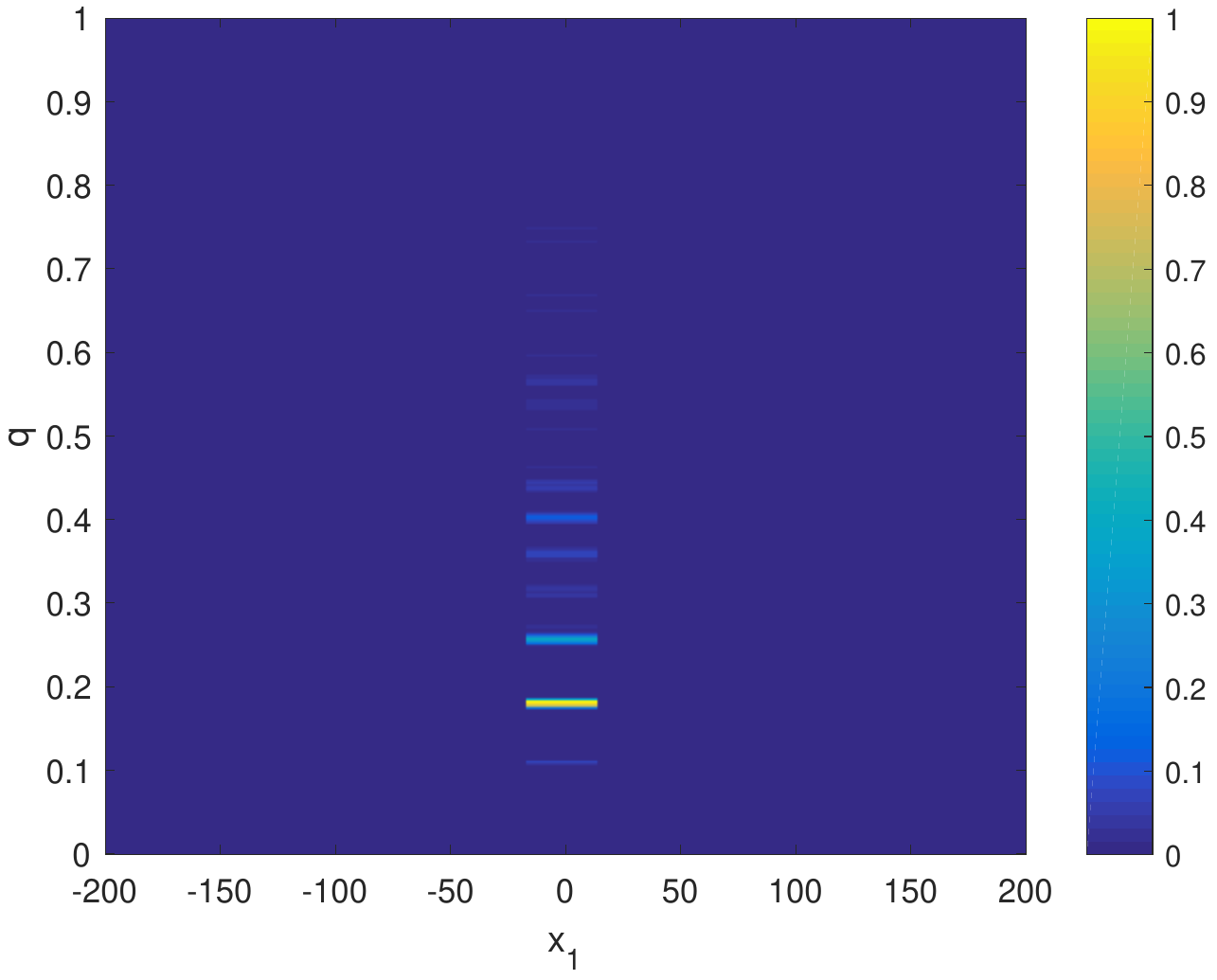}
\subcaption{stage 1}
\end{subfigure}
\begin{subfigure}{0.32\textwidth}
\includegraphics[width=0.9\linewidth, height=4cm, keepaspectratio]{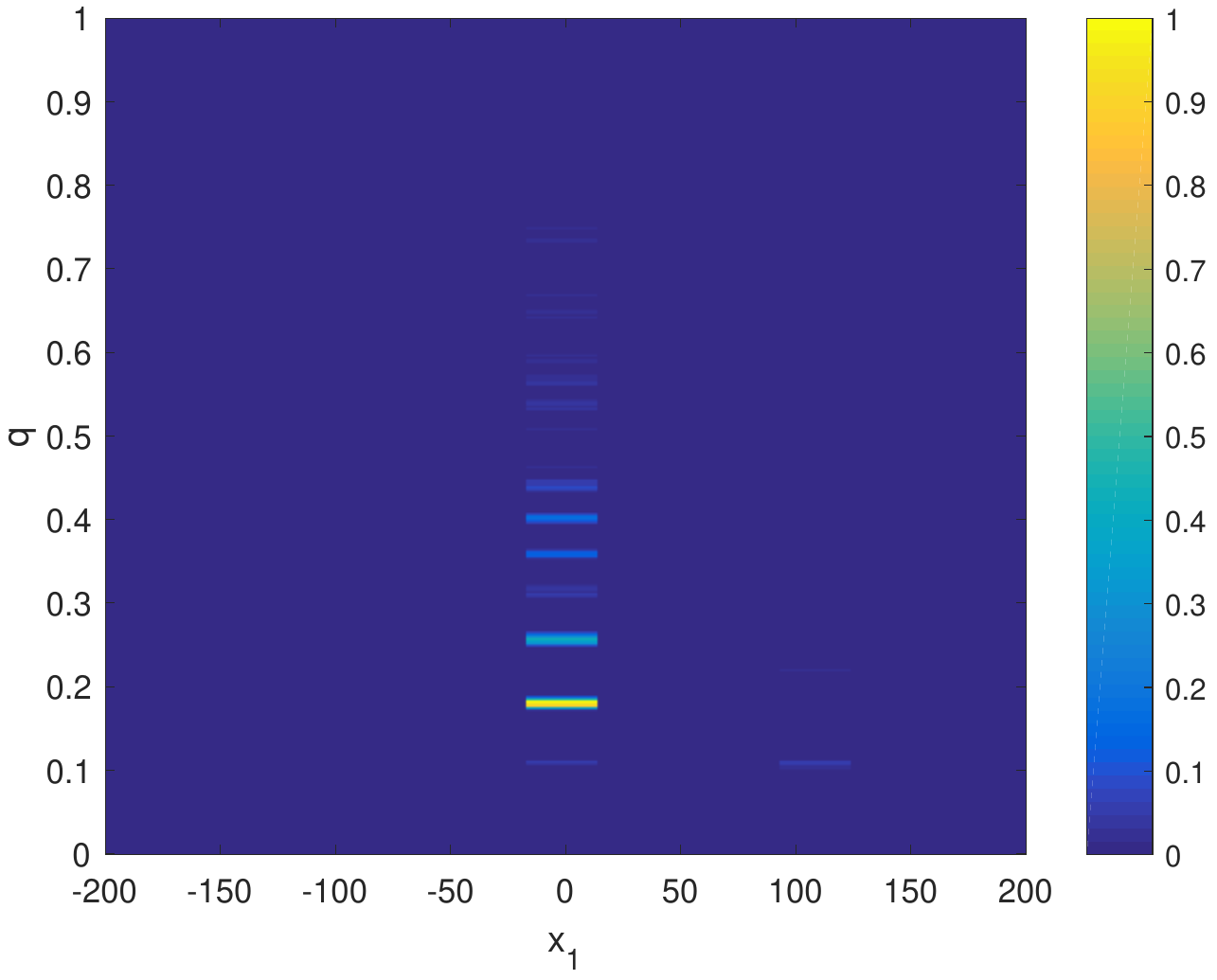} 
\subcaption{stage 2}
\end{subfigure}
\begin{subfigure}{0.35\textwidth}
\includegraphics[width=1\linewidth, height=4.8cm, keepaspectratio]{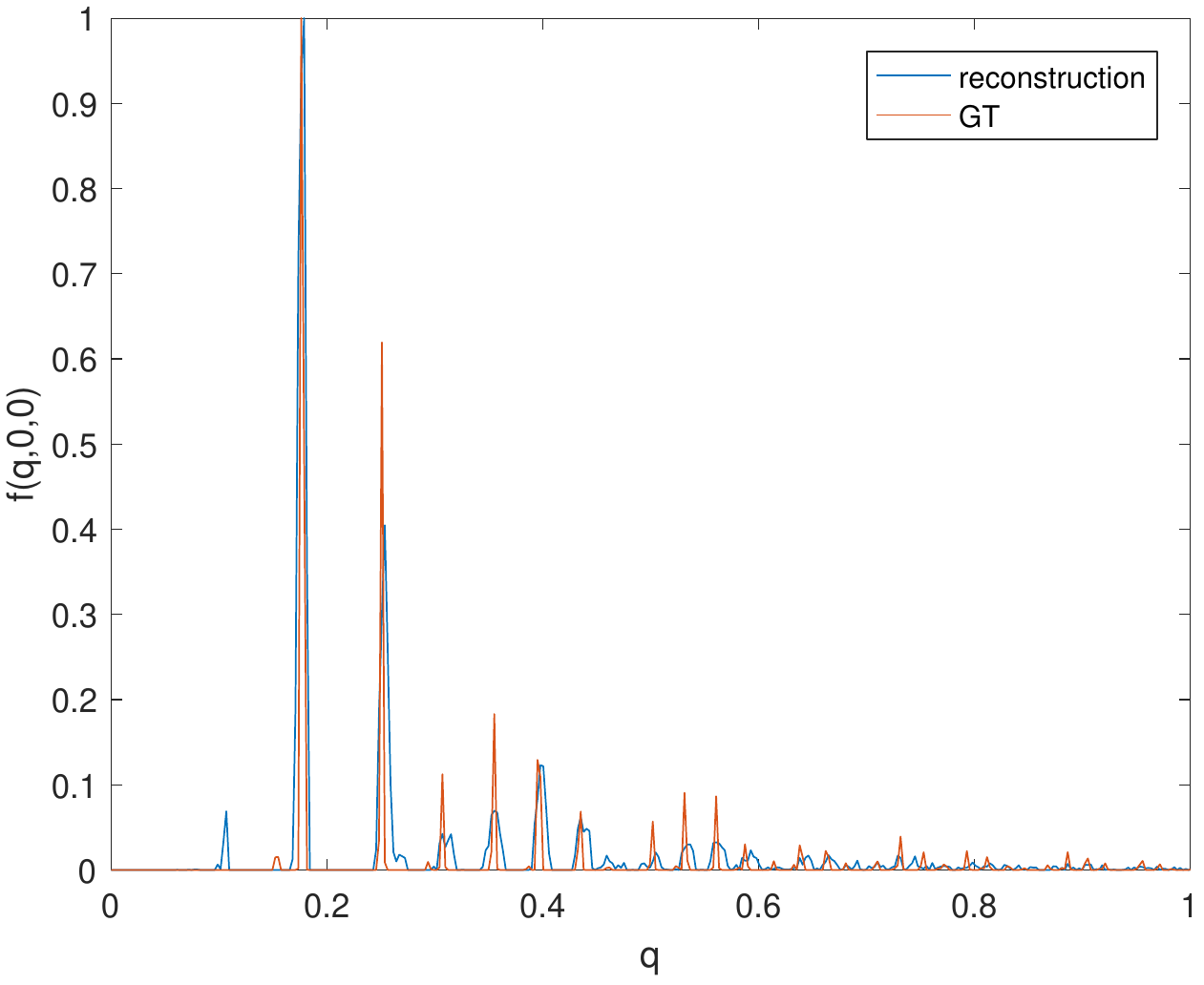}
\subcaption{stage 1}
\end{subfigure}
\begin{subfigure}{0.35\textwidth}
\includegraphics[width=1\linewidth, height=4.8cm, keepaspectratio]{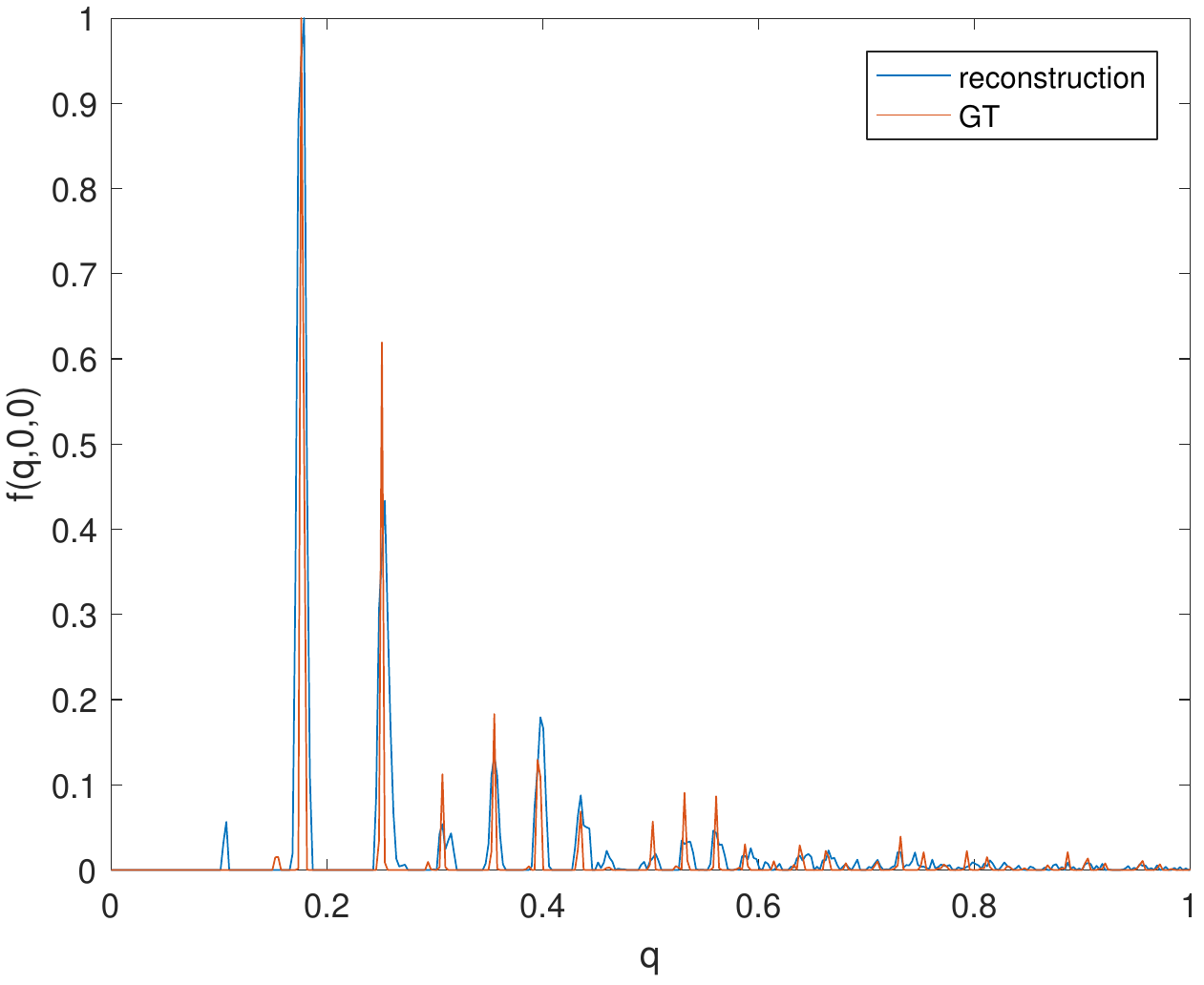} 
\subcaption{stage 2}
\end{subfigure}
\caption{Illustration of two-stage reconstruction process used for 2DBSR. We show reconstructions of an NaCl sphere, center $x_c=0$mm and $r=15$mm at each stage of the reconstruction process. Figures (A)-(C) show the ground truth and the 2-D reconstructions, and figures (D)-(F) show the central line profiles of the 2-D reconstructions corresponding to each stage. The $F_1$ scores corresponding to each stage are, stage 1 - $F_1=0.66$, stage 2 - $F_1=0.76$.}
\label{stages}
\end{figure}

\begin{figure*}
\centering
\setlength{\tabcolsep}{5pt}
\begin{tabular}{ c|cc }  
  &  2DBSR (S2) & FTV \\ \hline \\[-0.4cm] 
 \rotatebox{90}{\hspace{2.5cm} GT}  & 
   \includegraphics[ width=0.4\linewidth, height=0.4\linewidth, keepaspectratio]{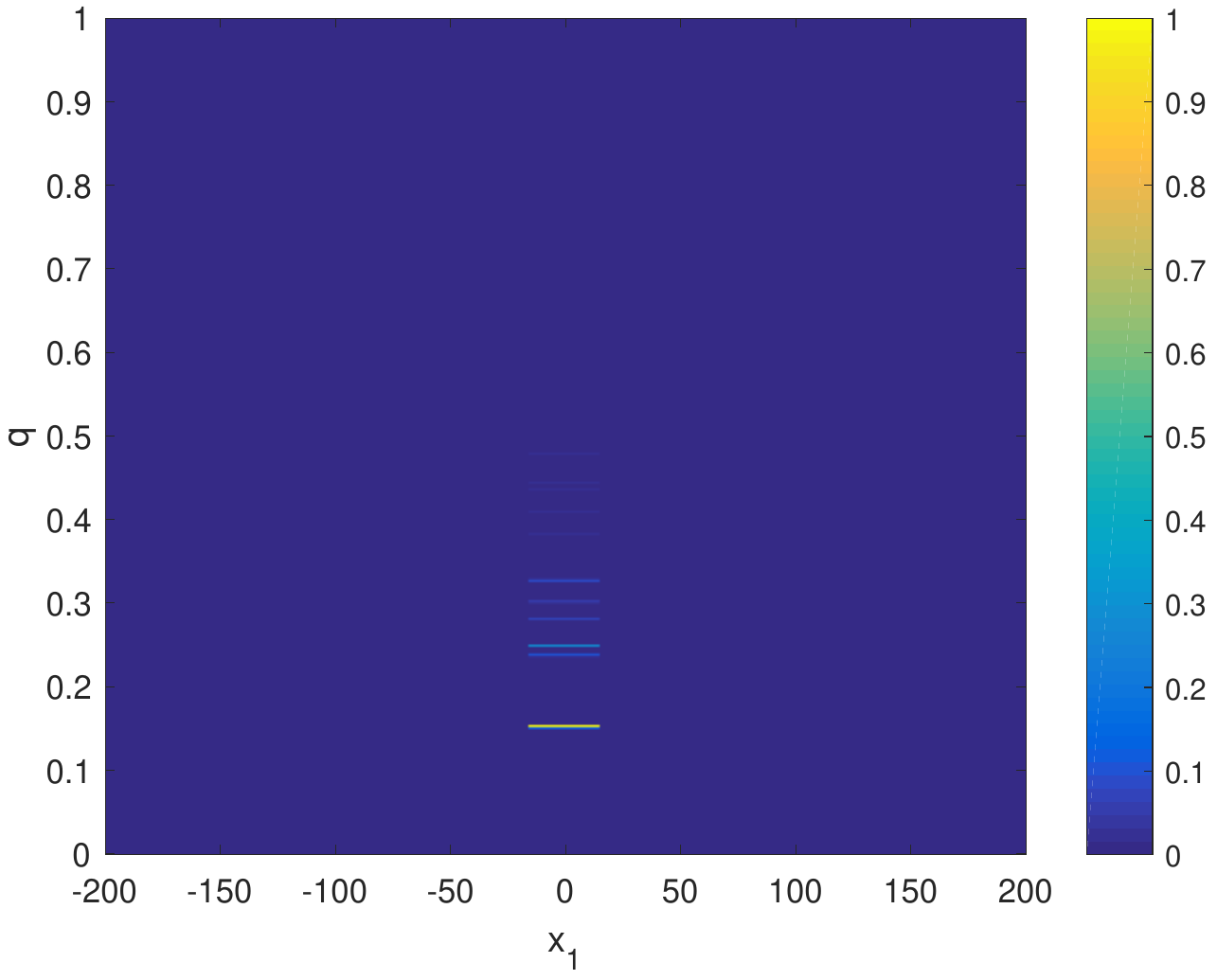} &
   \includegraphics[ width=0.4\linewidth, height=0.4\linewidth, keepaspectratio]{C_r15_GT} \\
 \rotatebox{90}{\hspace{1cm} 2-D reconstruction} &
  \includegraphics[ width=0.4\linewidth, height=0.4\linewidth, keepaspectratio]{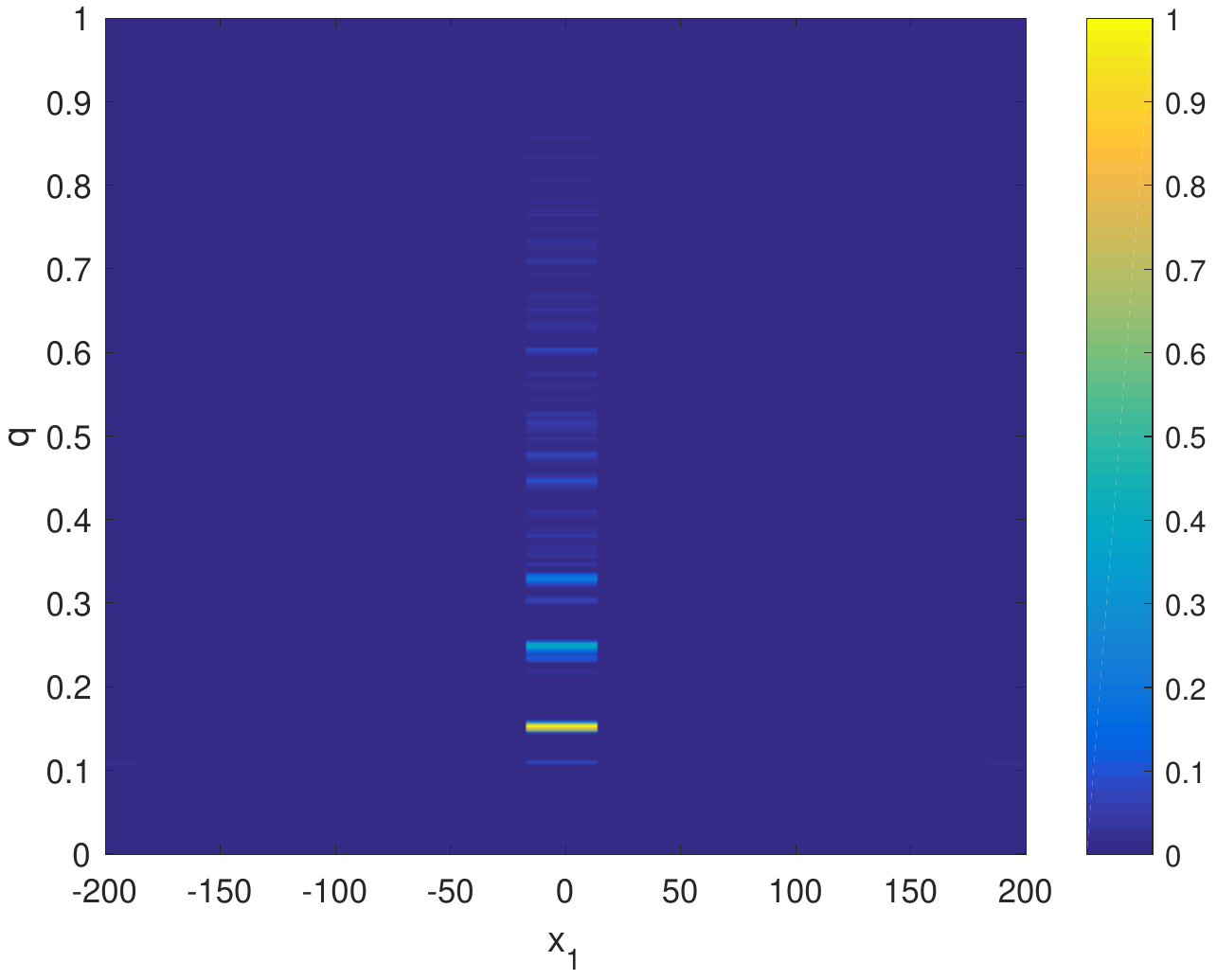} &
   \includegraphics[ width=0.4\linewidth, height=0.4\linewidth, keepaspectratio]{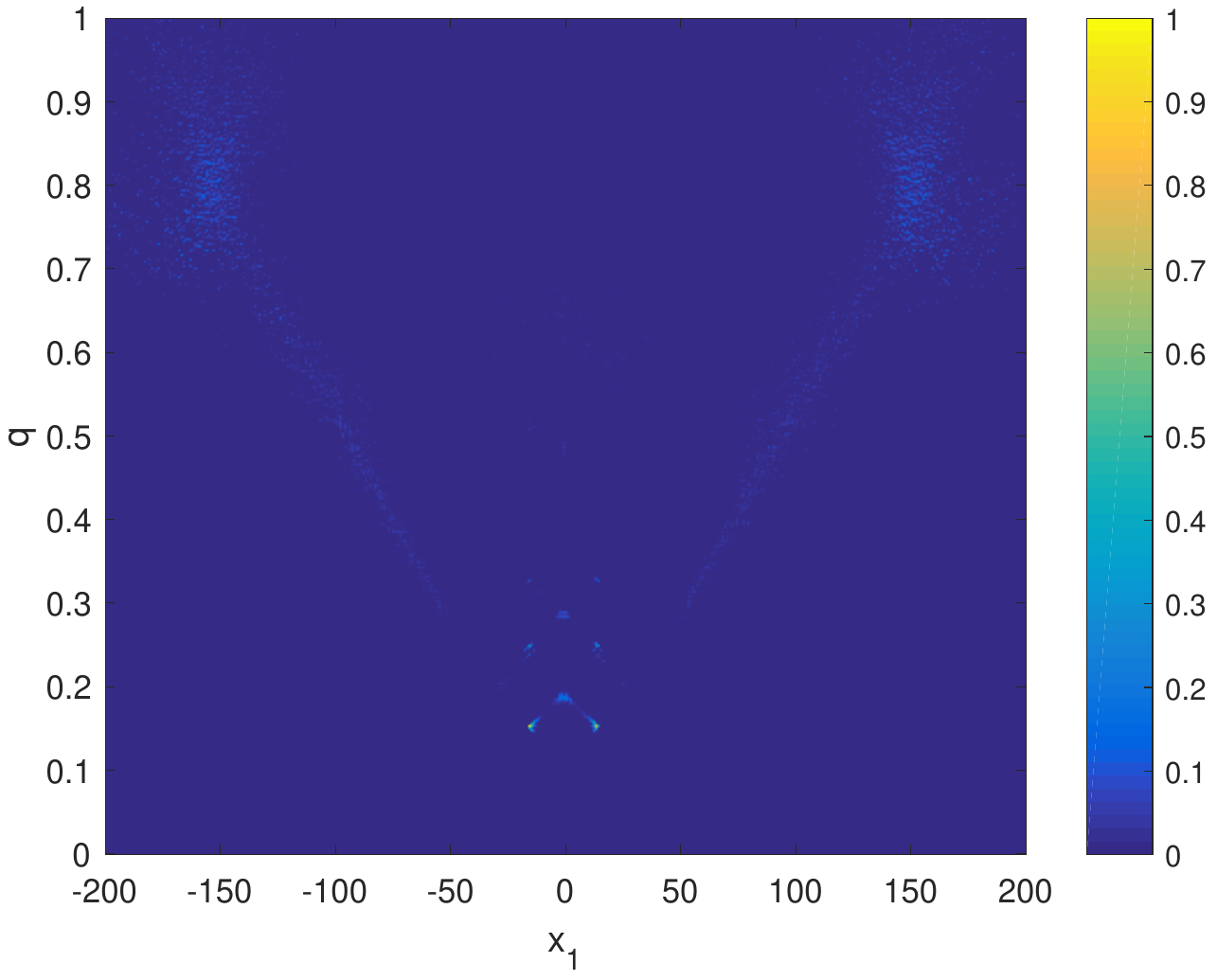} \\
 \rotatebox{90}{\hspace{1.3cm} LP ($x_1=0$mm)} &
   \includegraphics[ width=0.4\linewidth, height=0.4\linewidth, keepaspectratio]{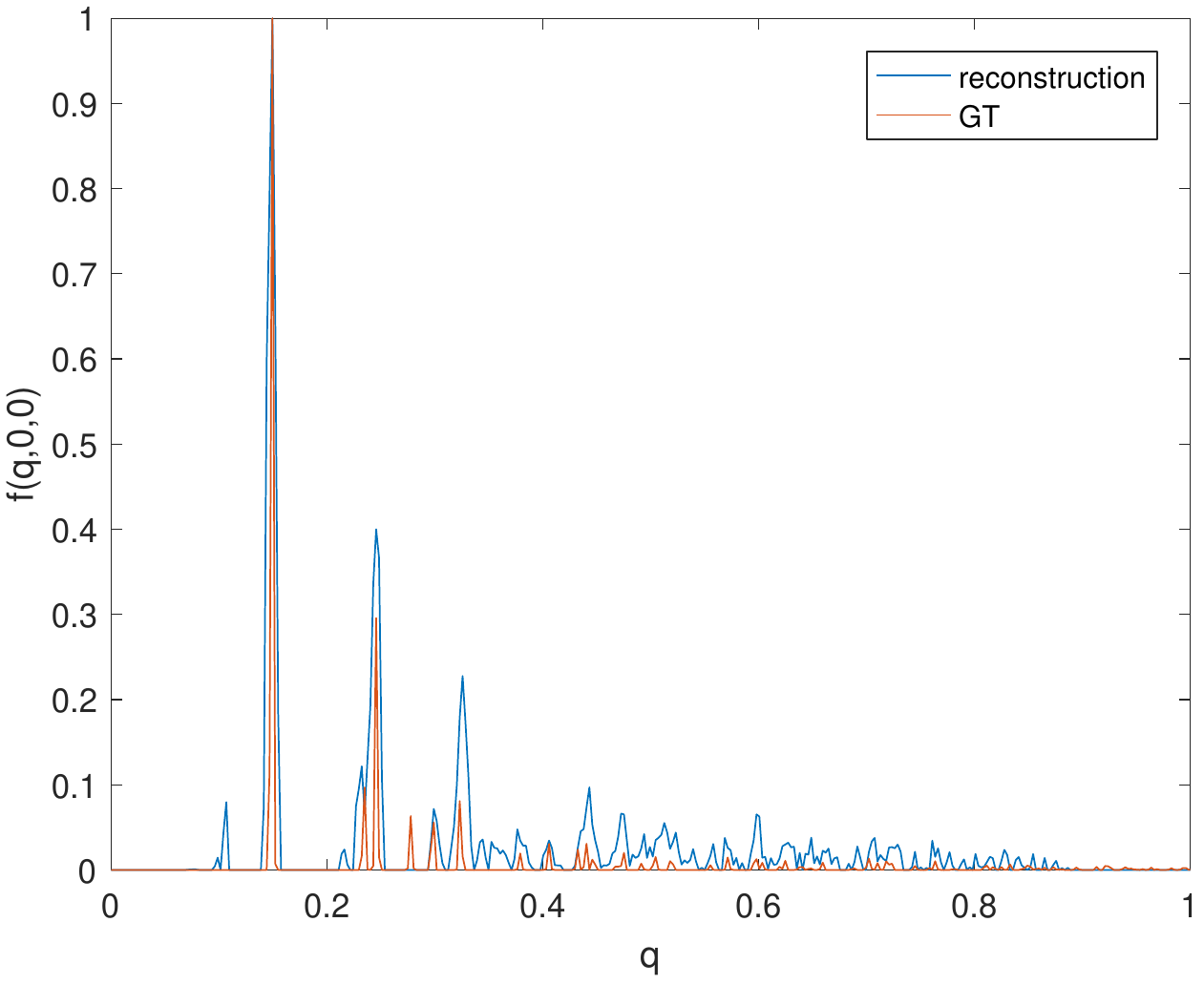} &
   \includegraphics[ width=0.4\linewidth, height=0.4\linewidth, keepaspectratio]{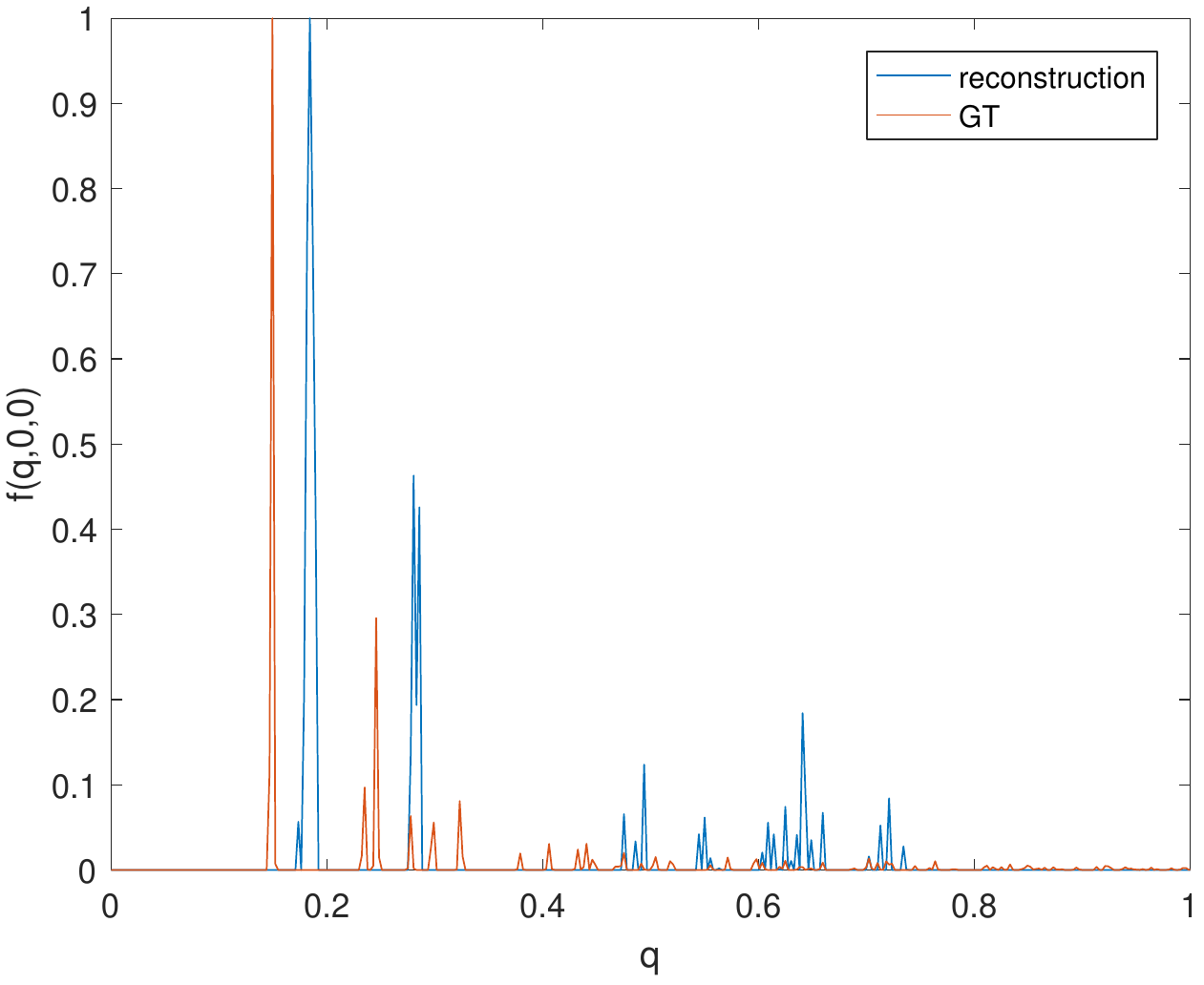}
\end{tabular}
\caption{C-graphite sphere reconstruction, $x_c=0$mm, $r=15$mm.}
\label{F0}
\end{figure*}

Let $x_c\in [-200,200]$mm denote the sphere center $x_1$ coordinate. The sphere center $x_2$ and $x_3$ coordinates are fixed at $x_2=x_3=0$mm. Let $r\in [5,15]$mm denote the sphere radius. For example, the spheres shown in figure \ref{Parcel} have center $x_1$ coordinates $x_c=\pm 75$mm and radius $r=15$mm. Then we consider nine Monte Carlo experiments in total which are designed to cover a range of $x_c$ and $r$, for materials contained inside and outside the imaging basis used (as detailed in section \ref{IP}).
\begin{figure*}
\centering
\setlength{\tabcolsep}{5pt}
\begin{tabular}{ c|cc }  
  &  2DBSR (S2) & FTV \\ \hline \\[-0.4cm] 
 \rotatebox{90}{\hspace{2.5cm} GT}  & 
   \includegraphics[ width=0.4\linewidth, height=0.4\linewidth, keepaspectratio]{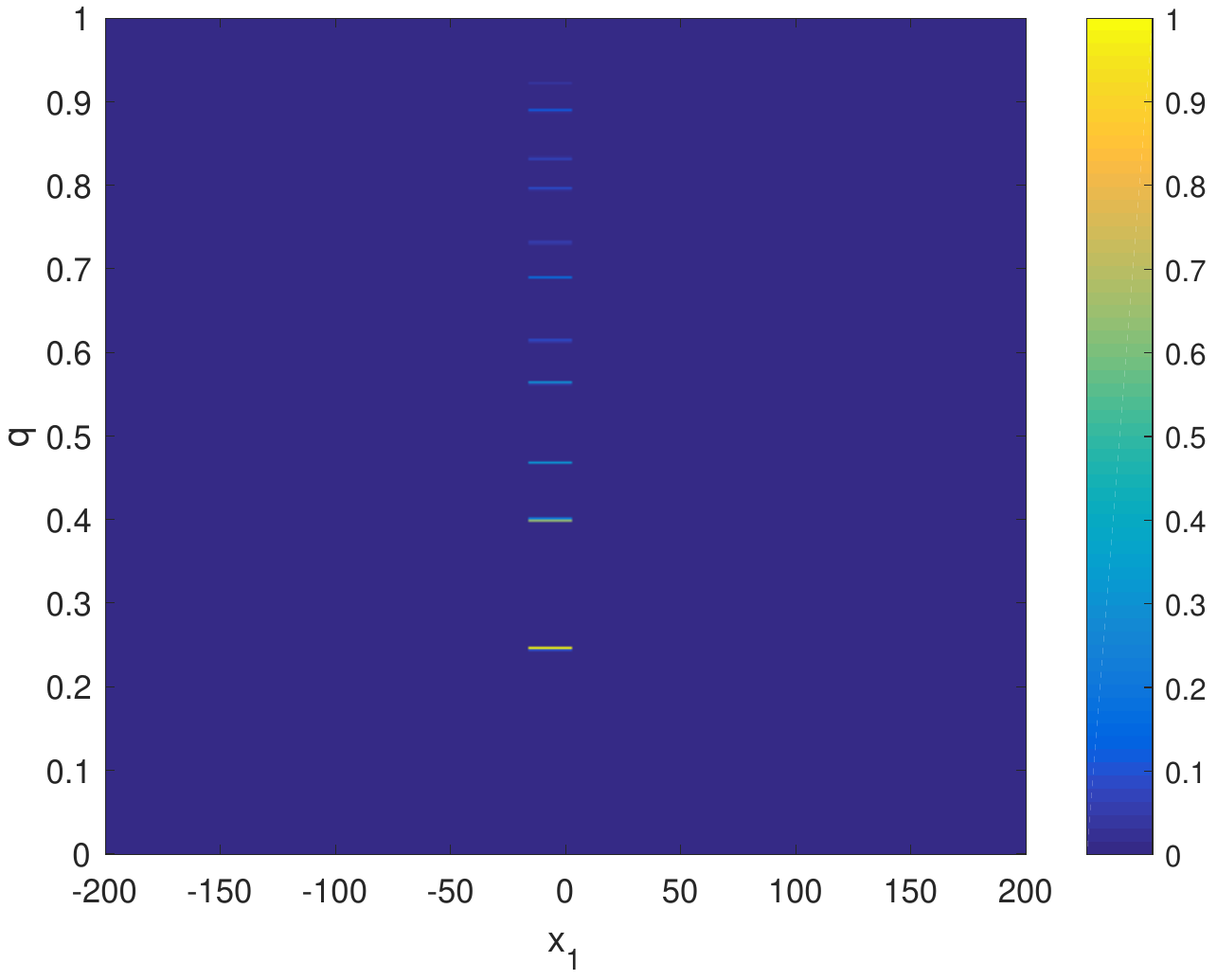} &
   \includegraphics[ width=0.4\linewidth, height=0.4\linewidth, keepaspectratio]{DiShift_GT} \\
 \rotatebox{90}{\hspace{1cm} 2-D reconstruction} &
  \includegraphics[ width=0.4\linewidth, height=0.4\linewidth, keepaspectratio]{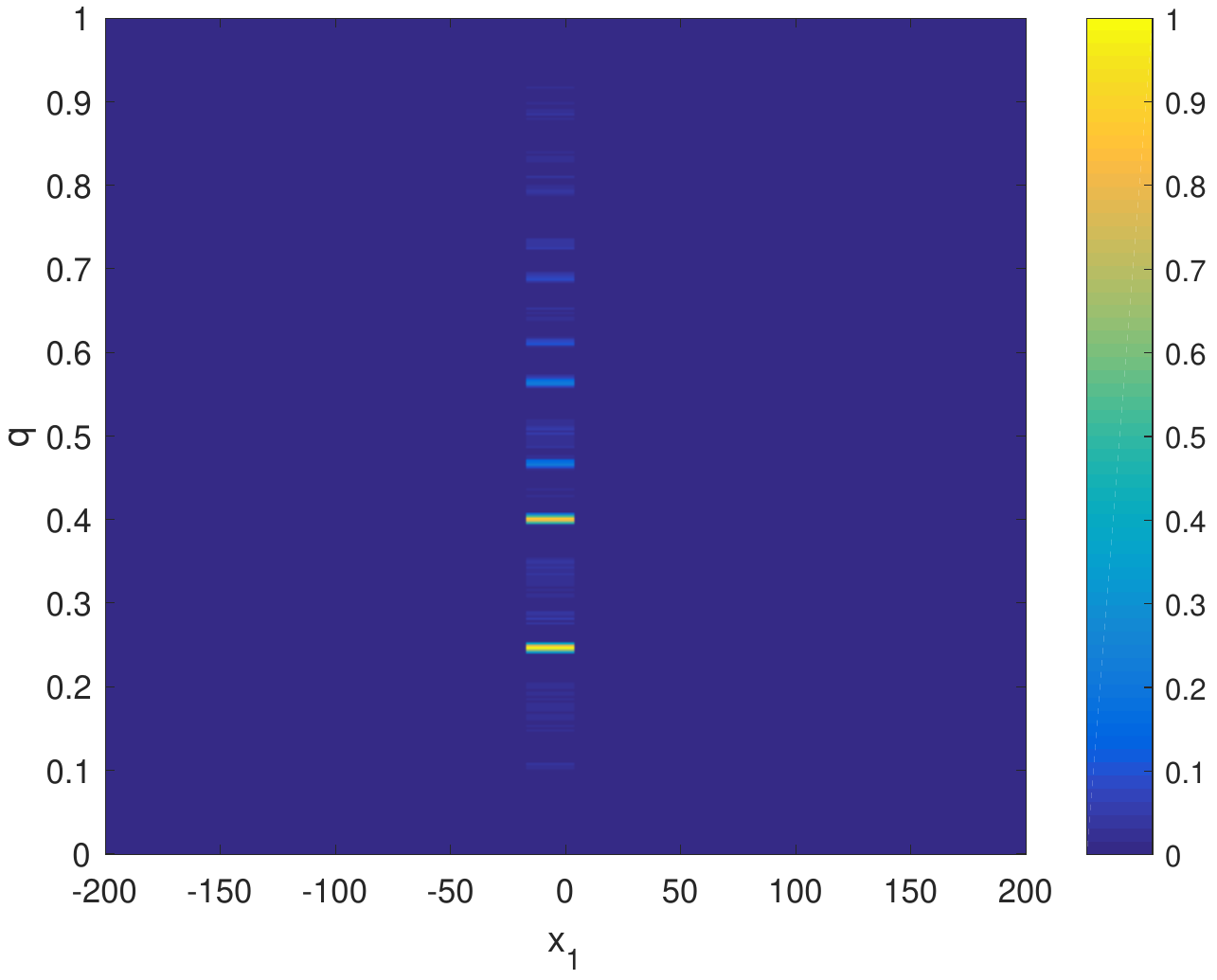} &
   \includegraphics[ width=0.4\linewidth, height=0.4\linewidth, keepaspectratio]{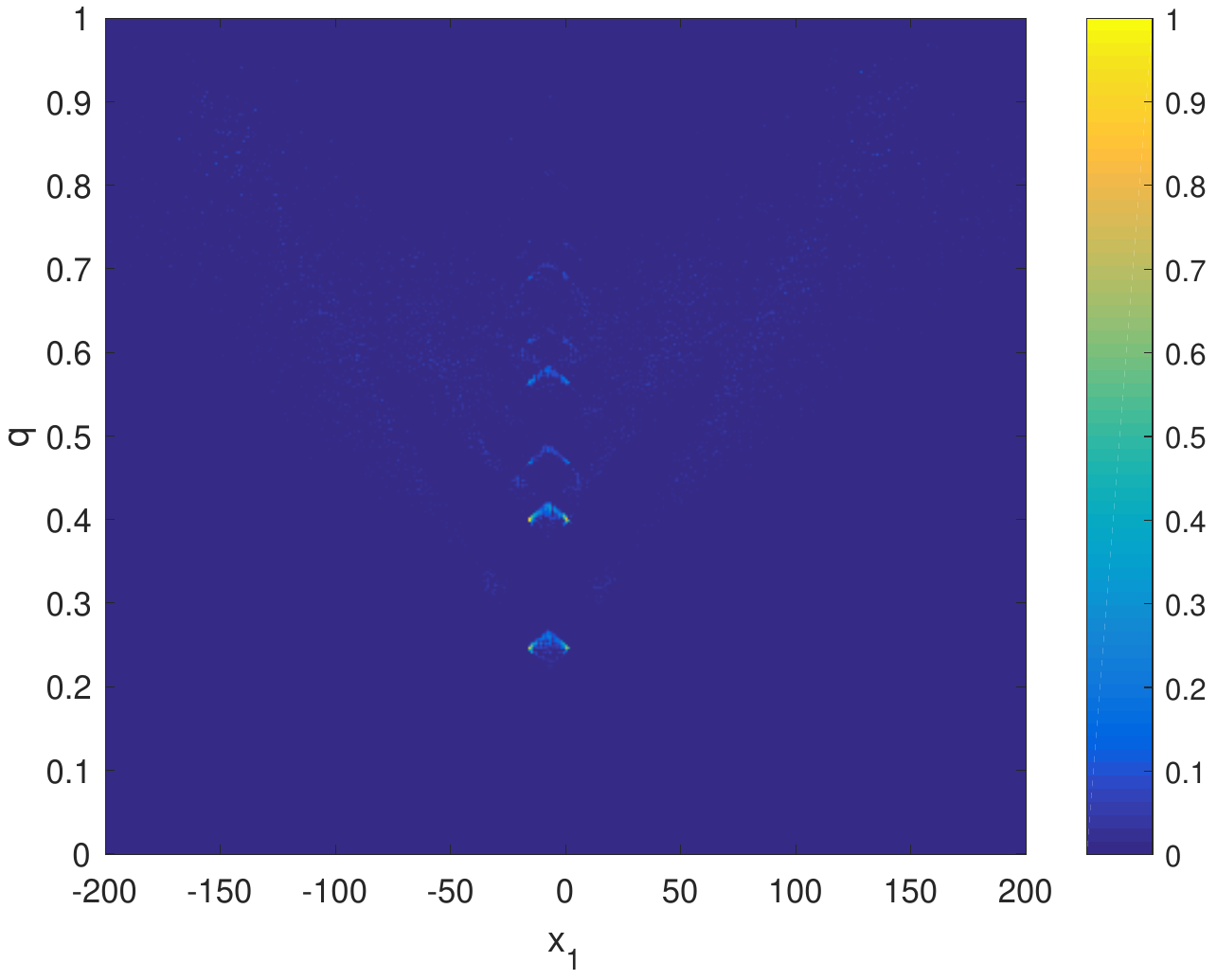} \\
 \rotatebox{90}{\hspace{1.3cm} LP ($x_1=0$mm)} &
   \includegraphics[ width=0.4\linewidth, height=0.4\linewidth, keepaspectratio]{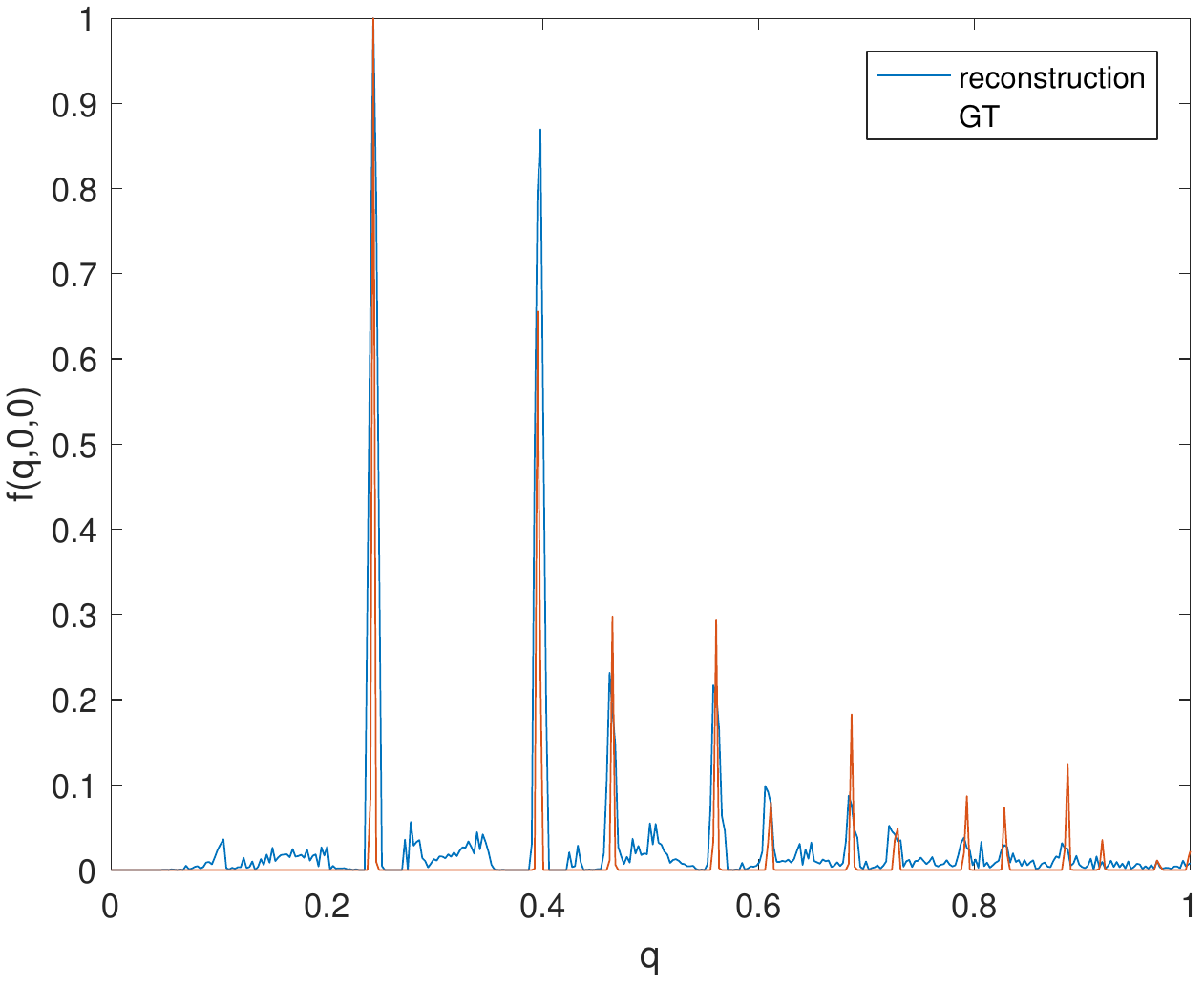} &
   \includegraphics[ width=0.4\linewidth, height=0.4\linewidth, keepaspectratio]{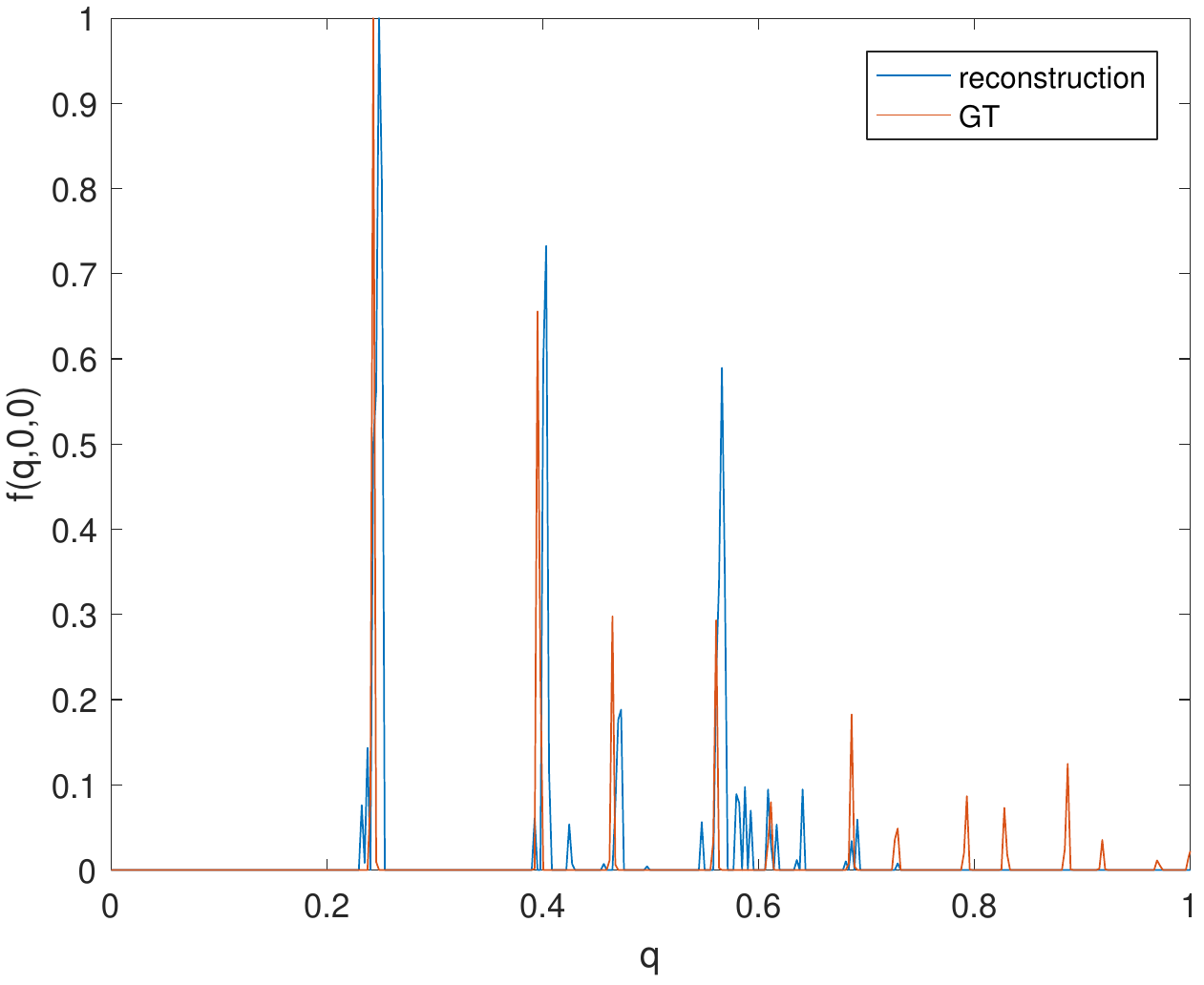}
\end{tabular}
\caption{C-diamond sphere reconstruction, $x_c=-6.25$mm, $r=8.75$mm.}
\label{F1}
\end{figure*}

\begin{figure*}
\centering
\setlength{\tabcolsep}{5pt}
\begin{tabular}{ c|cc }  
  &  2DBSR (S2) & FTV \\ \hline \\[-0.4cm] 
 \rotatebox{90}{\hspace{2.5cm} GT}  & 
   \includegraphics[ width=0.4\linewidth, height=0.4\linewidth, keepaspectratio]{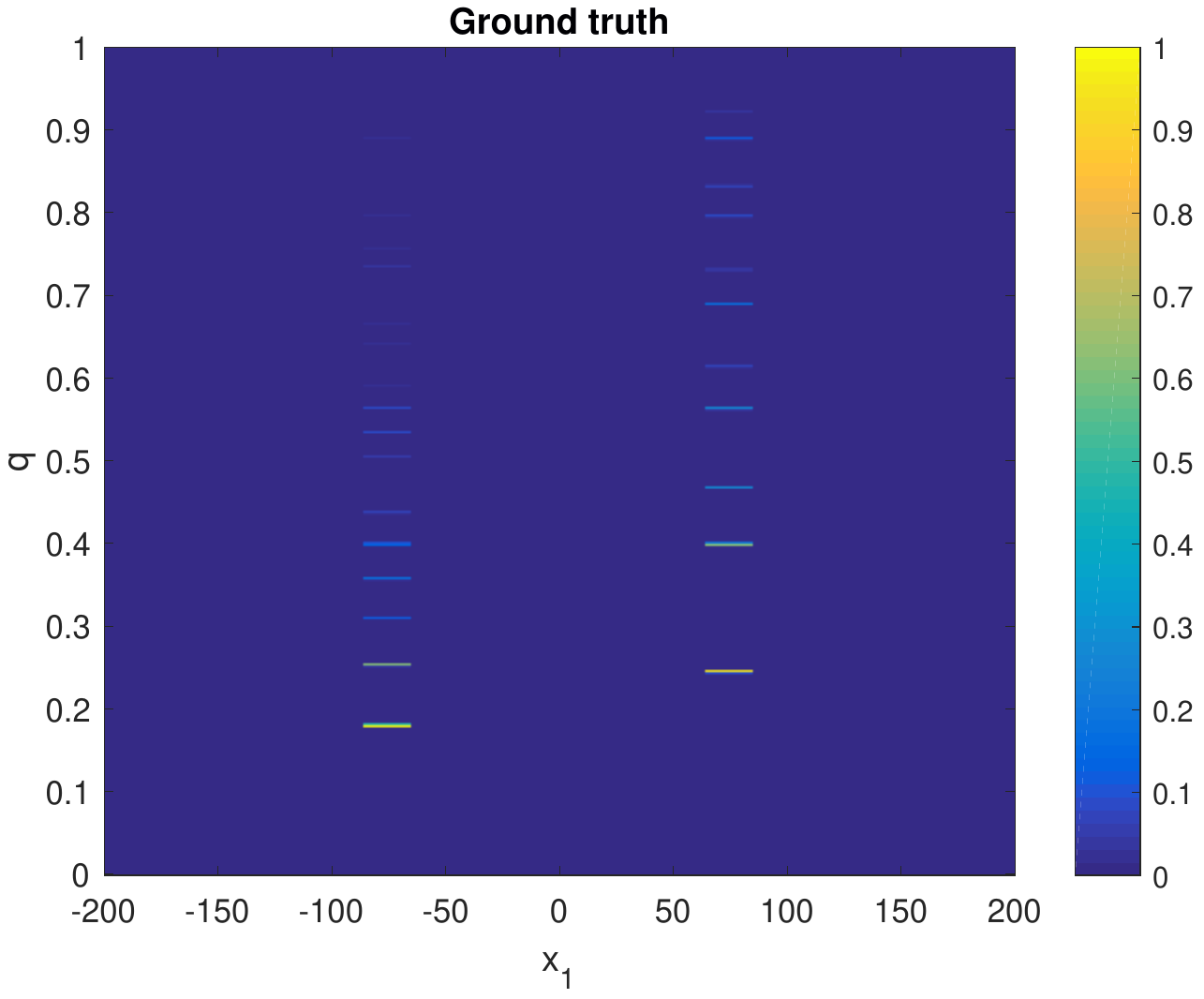} &
   \includegraphics[ width=0.4\linewidth, height=0.4\linewidth, keepaspectratio]{StDi_r10_GT} \\ 
 \rotatebox{90}{\hspace{1cm} 2-D reconstruction} &
  \includegraphics[ width=0.4\linewidth, height=0.4\linewidth, keepaspectratio]{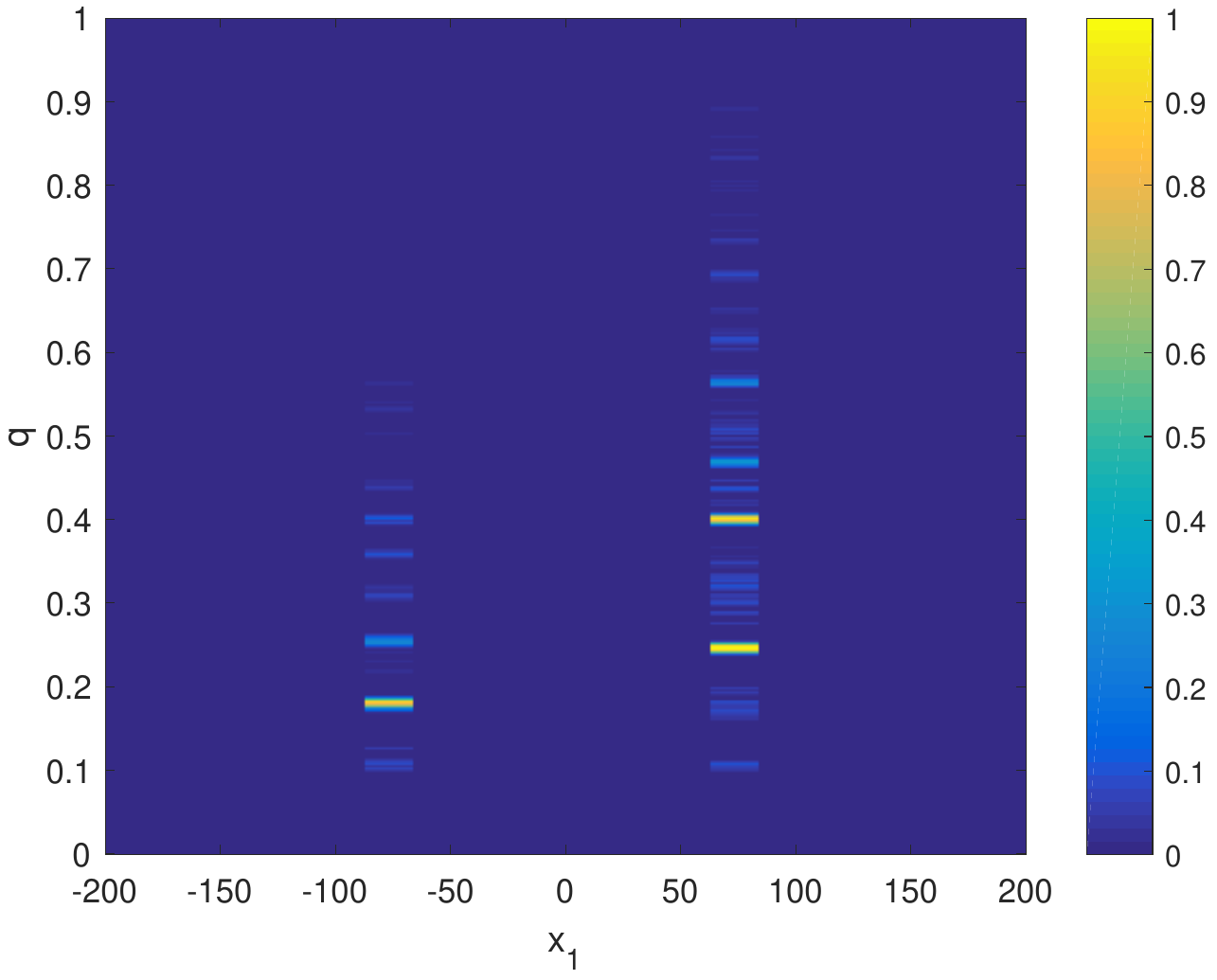} &
   \includegraphics[ width=0.4\linewidth, height=0.4\linewidth, keepaspectratio]{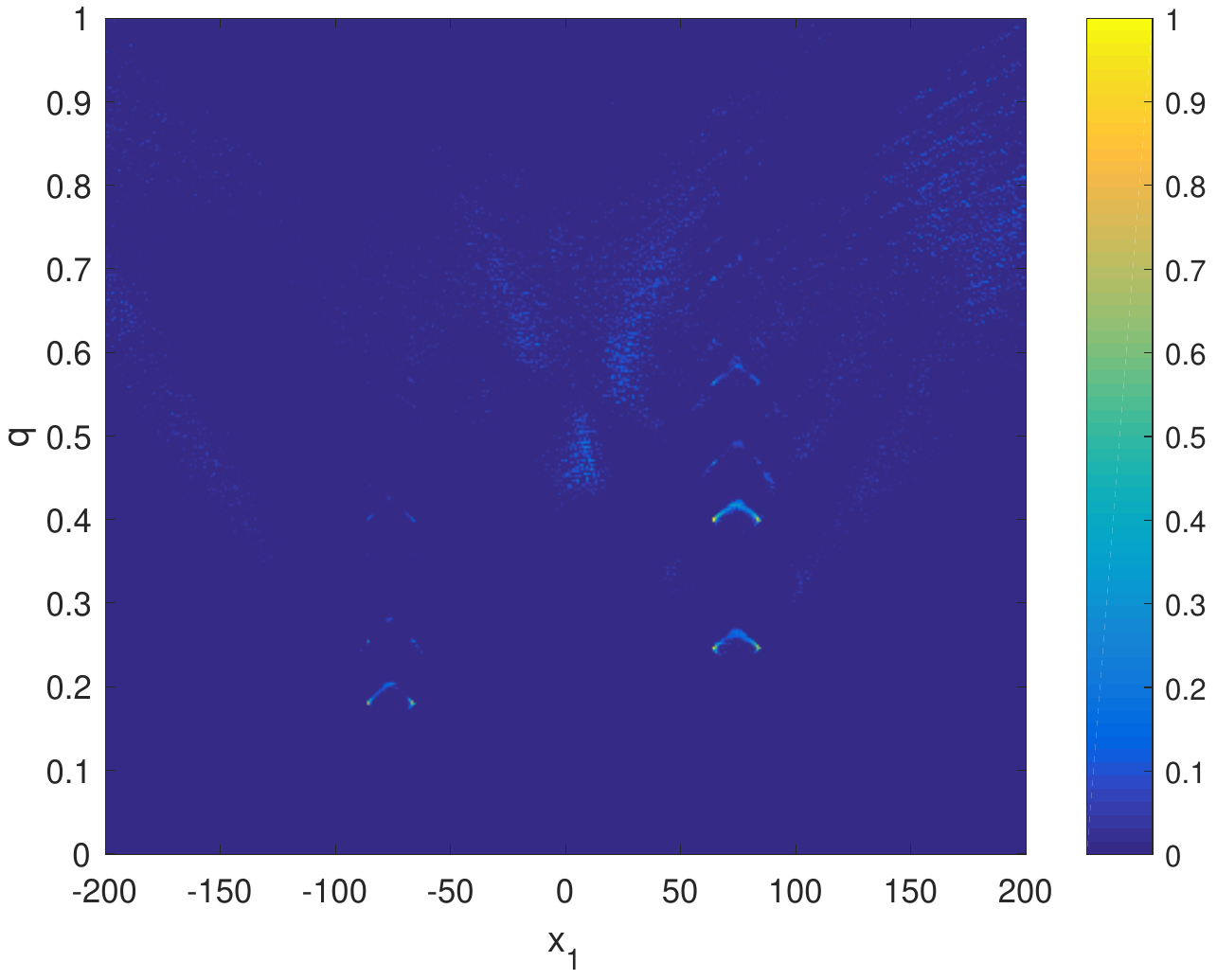} \\ 
 \rotatebox{90}{\hspace{1cm} LP1 ($x_1=75$mm)} &
   \includegraphics[ width=0.4\linewidth, height=0.4\linewidth, keepaspectratio]{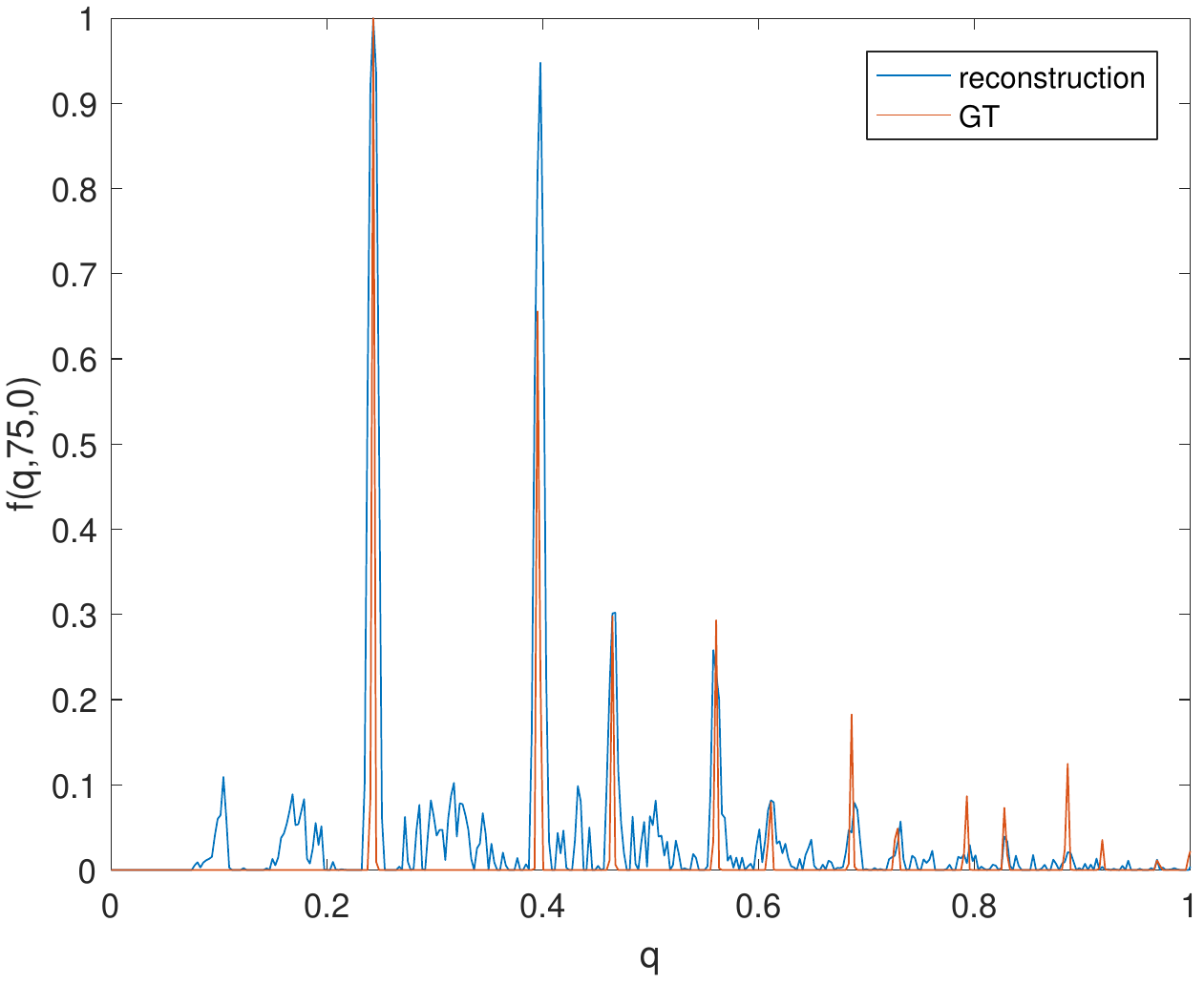} &
   \includegraphics[ width=0.4\linewidth, height=0.4\linewidth, keepaspectratio]{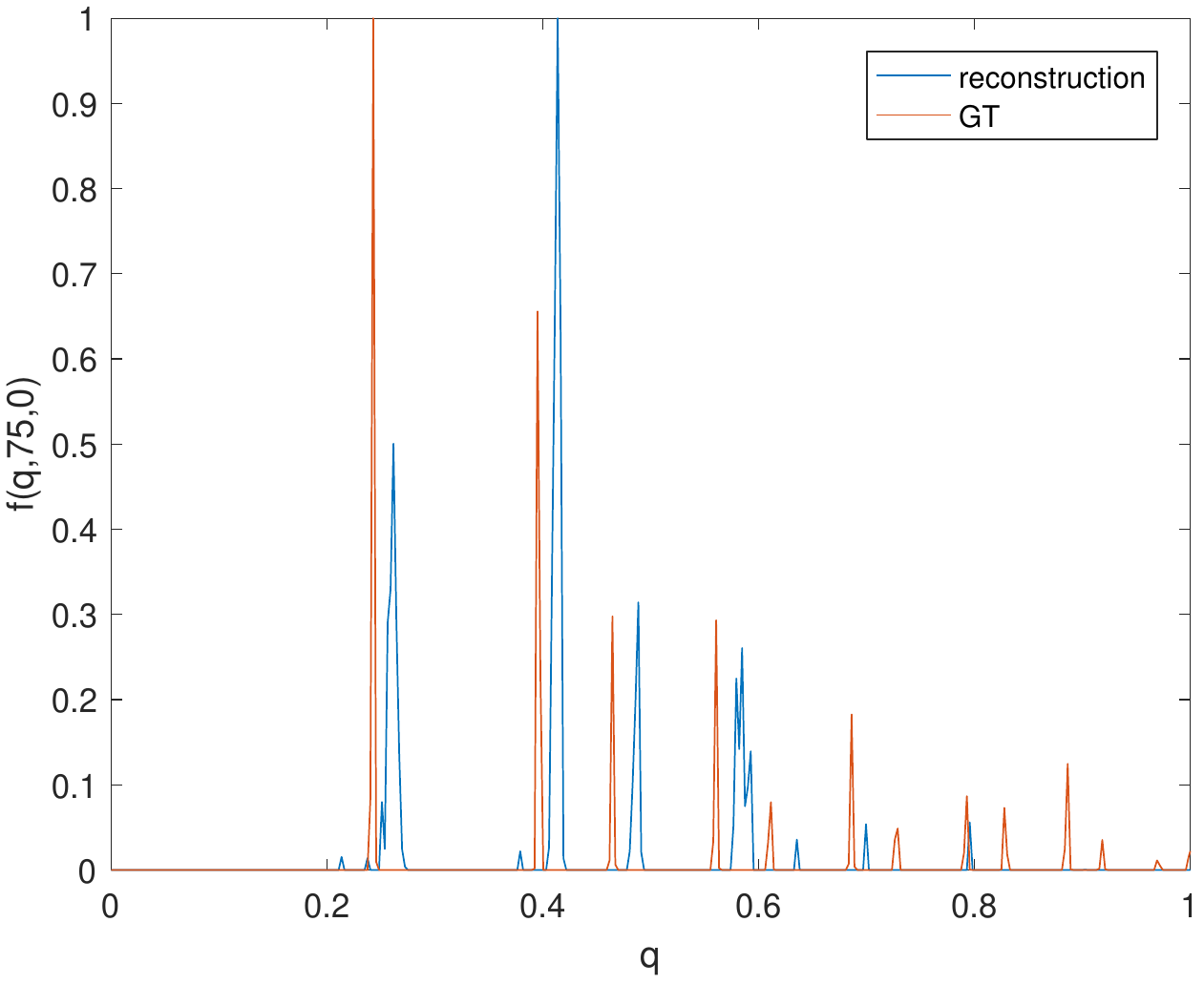} \\ 
 \rotatebox{90}{\hspace{1cm} LP2 ($x_1=-75$mm)} &
   \includegraphics[ width=0.4\linewidth, height=0.4\linewidth, keepaspectratio]{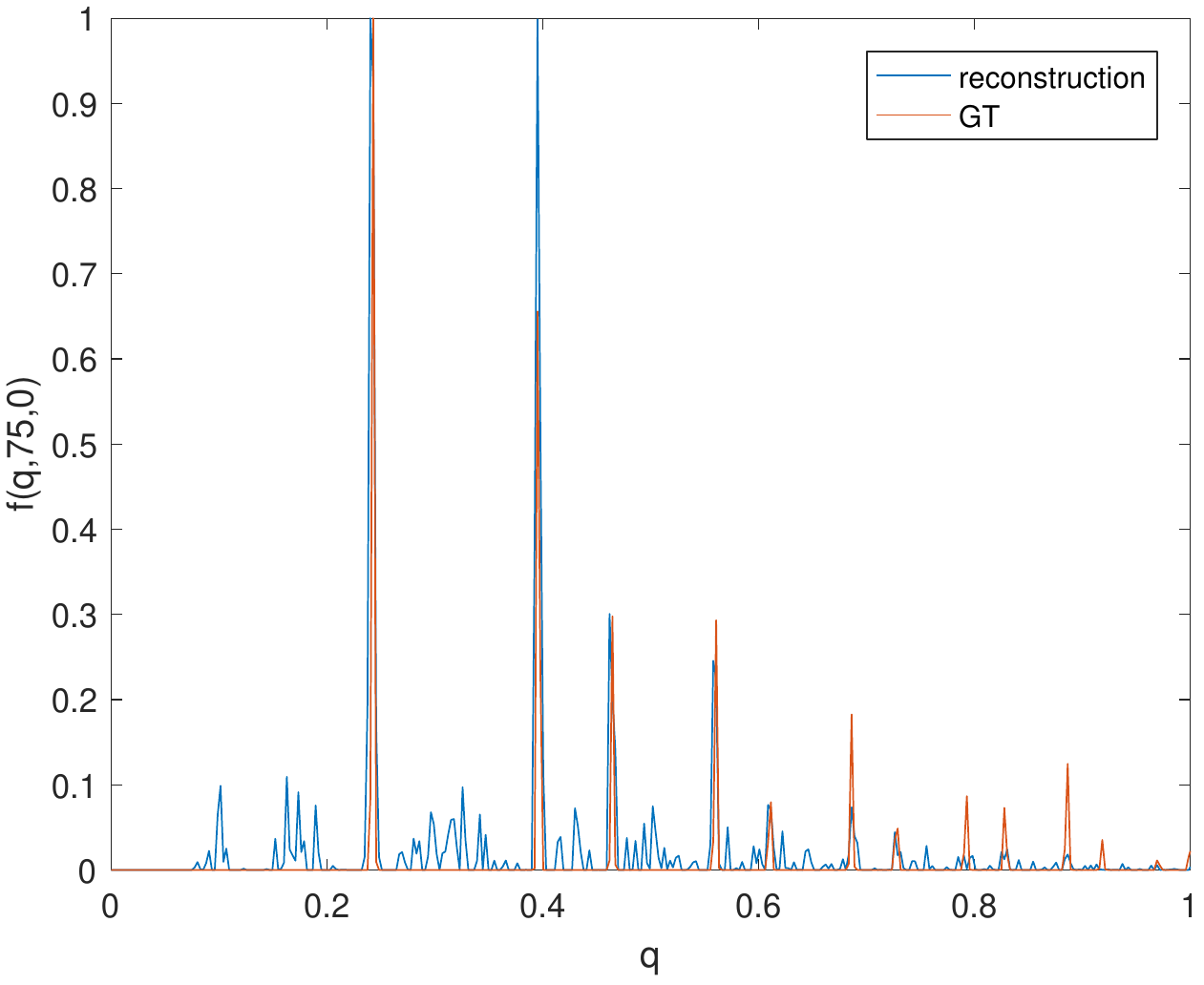} &
   \includegraphics[ width=0.4\linewidth, height=0.4\linewidth, keepaspectratio]{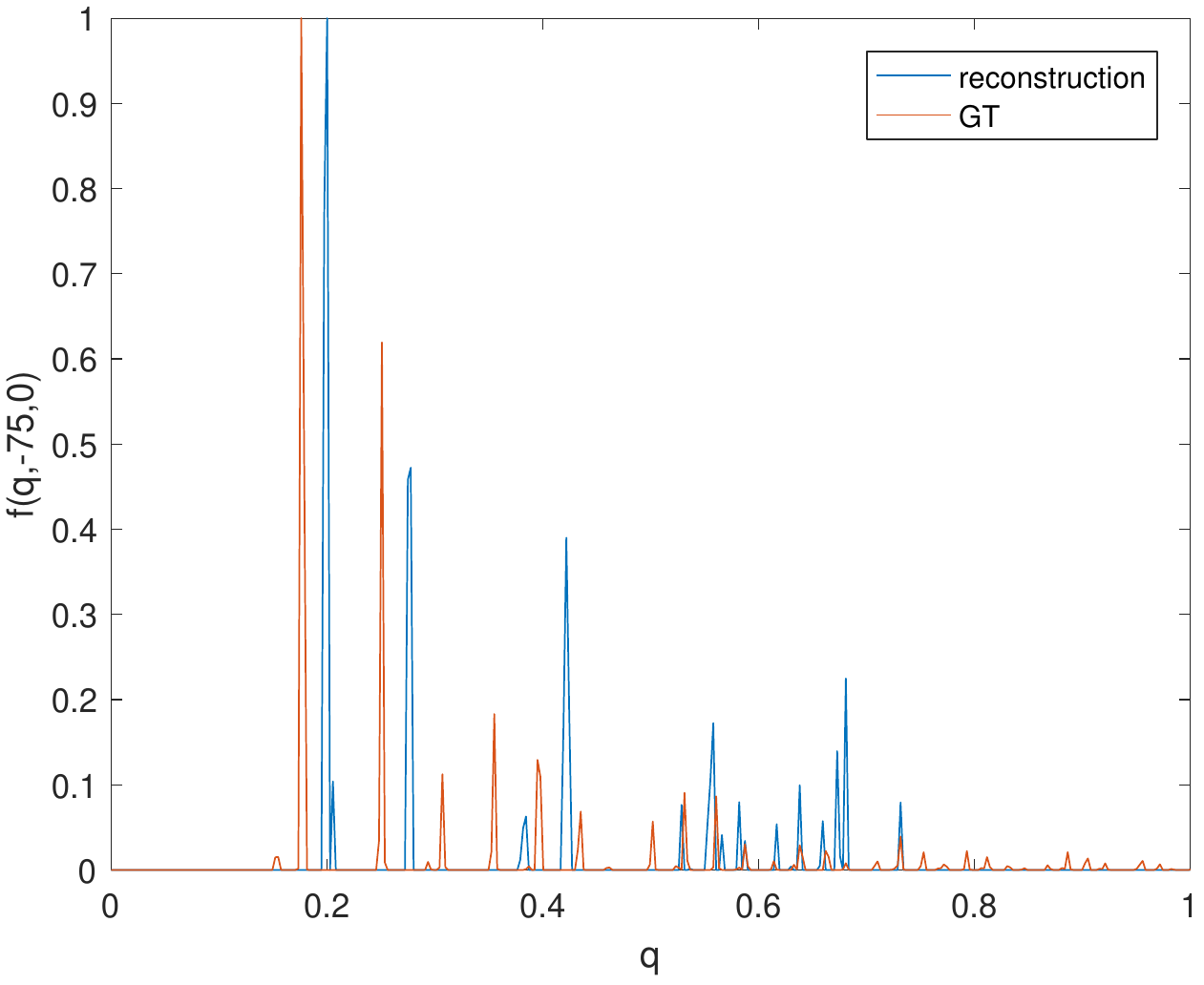} \\ 
\end{tabular}
\caption{NaCL and C-diamond sphere reconstruction, centers $x_c=-75$mm and $x_c=75$mm respectively, both with $r=10$mm.}
\label{F2}
\end{figure*}

\begin{figure*}
\centering
\setlength{\tabcolsep}{5pt}
\begin{tabular}{ c|cc }  
  &  2DBSR (S2) & FTV \\ \hline \\[-0.4cm] 
 \rotatebox{90}{\hspace{2.5cm} GT}  & 
   \includegraphics[ width=0.4\linewidth, height=0.4\linewidth, keepaspectratio]{StC_r10_GT} &
   \includegraphics[ width=0.4\linewidth, height=0.4\linewidth, keepaspectratio]{StC_r10_GT} \\ 
   
 \rotatebox{90}{\hspace{1cm} 2-D reconstruction} &
  \includegraphics[ width=0.4\linewidth, height=0.4\linewidth, keepaspectratio]{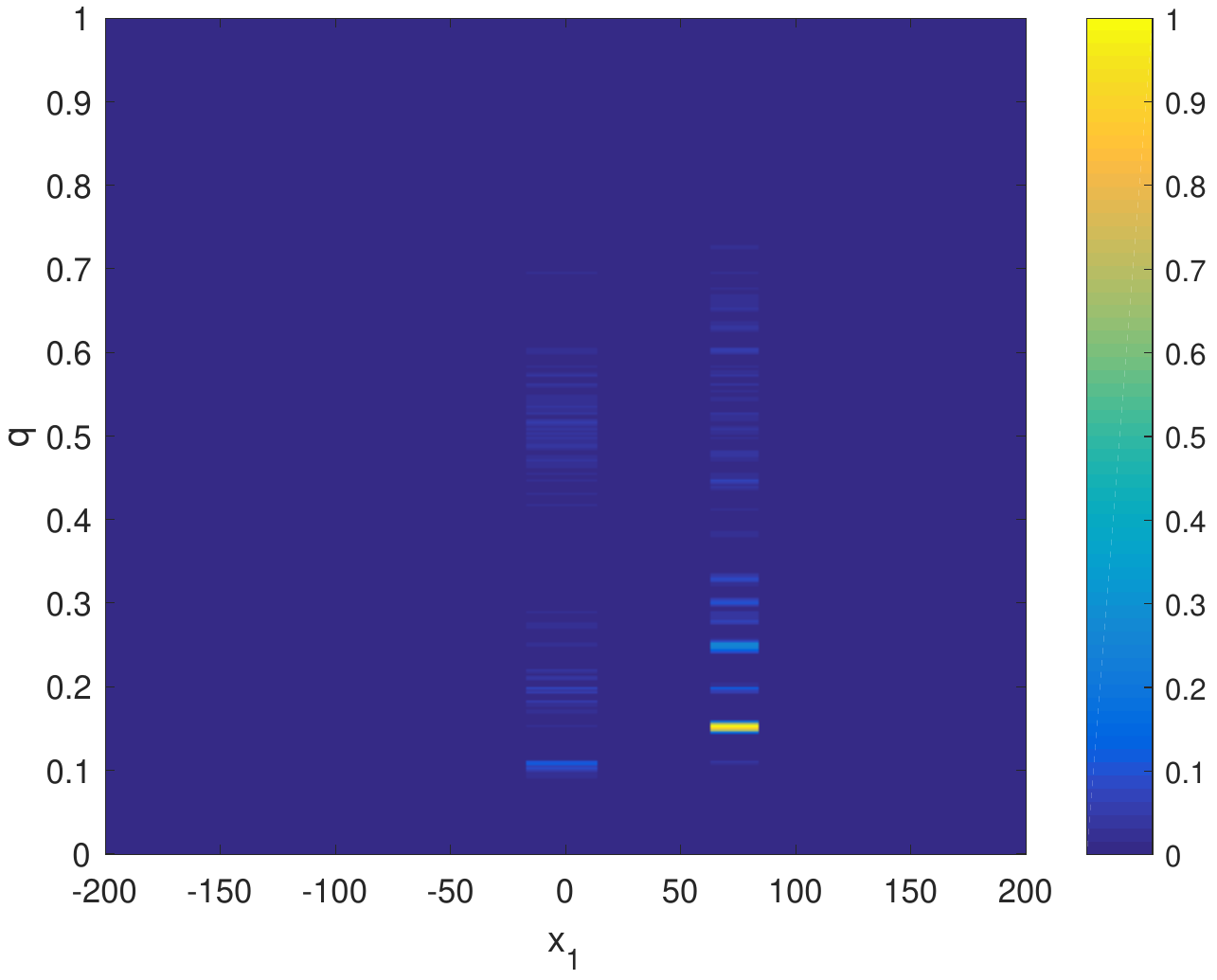} &
   \includegraphics[ width=0.4\linewidth, height=0.4\linewidth, keepaspectratio]{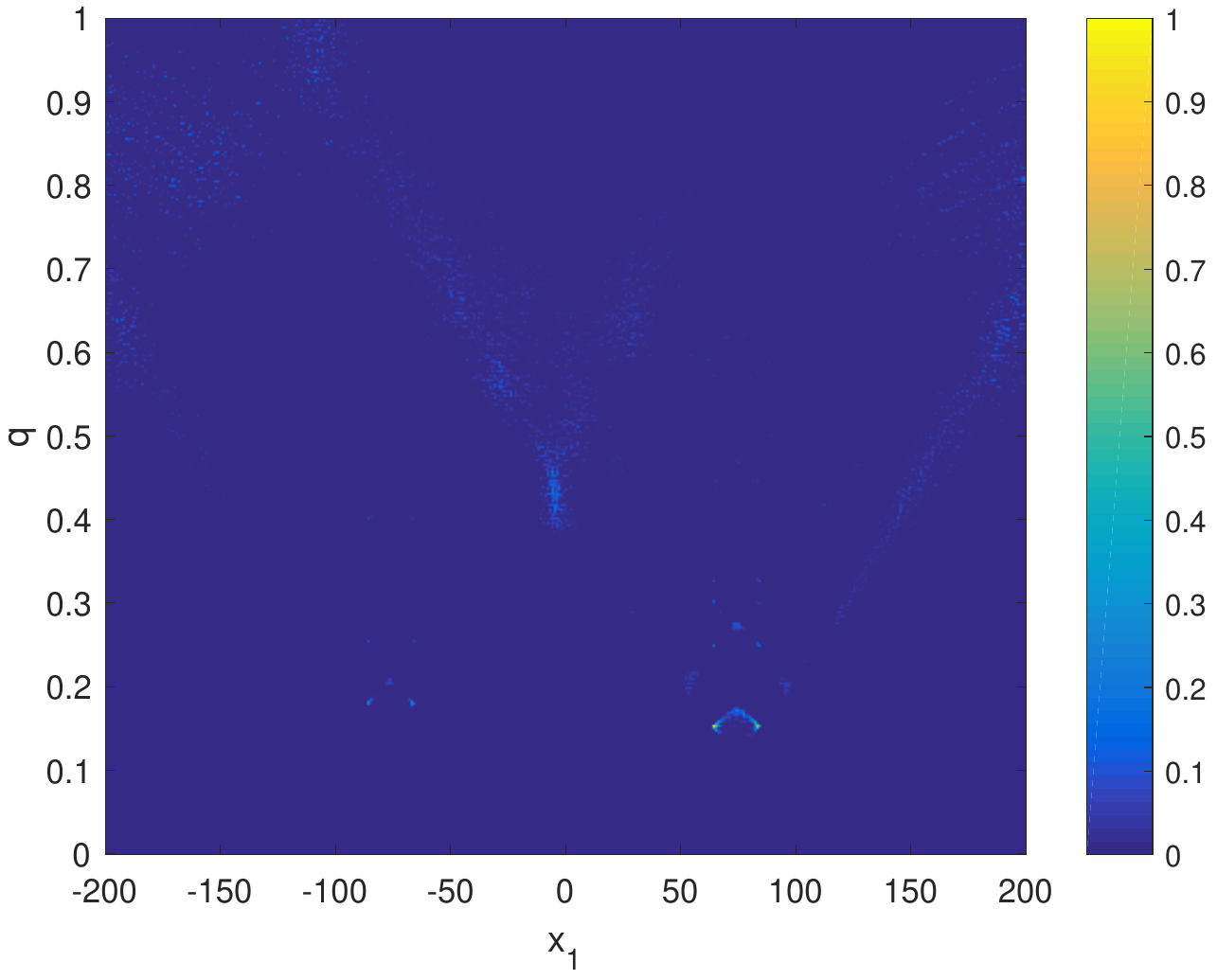} \\ 
   
 \rotatebox{90}{\hspace{1cm} LP1 ($x_1=75$mm)} &
   \includegraphics[ width=0.4\linewidth, height=0.4\linewidth, keepaspectratio]{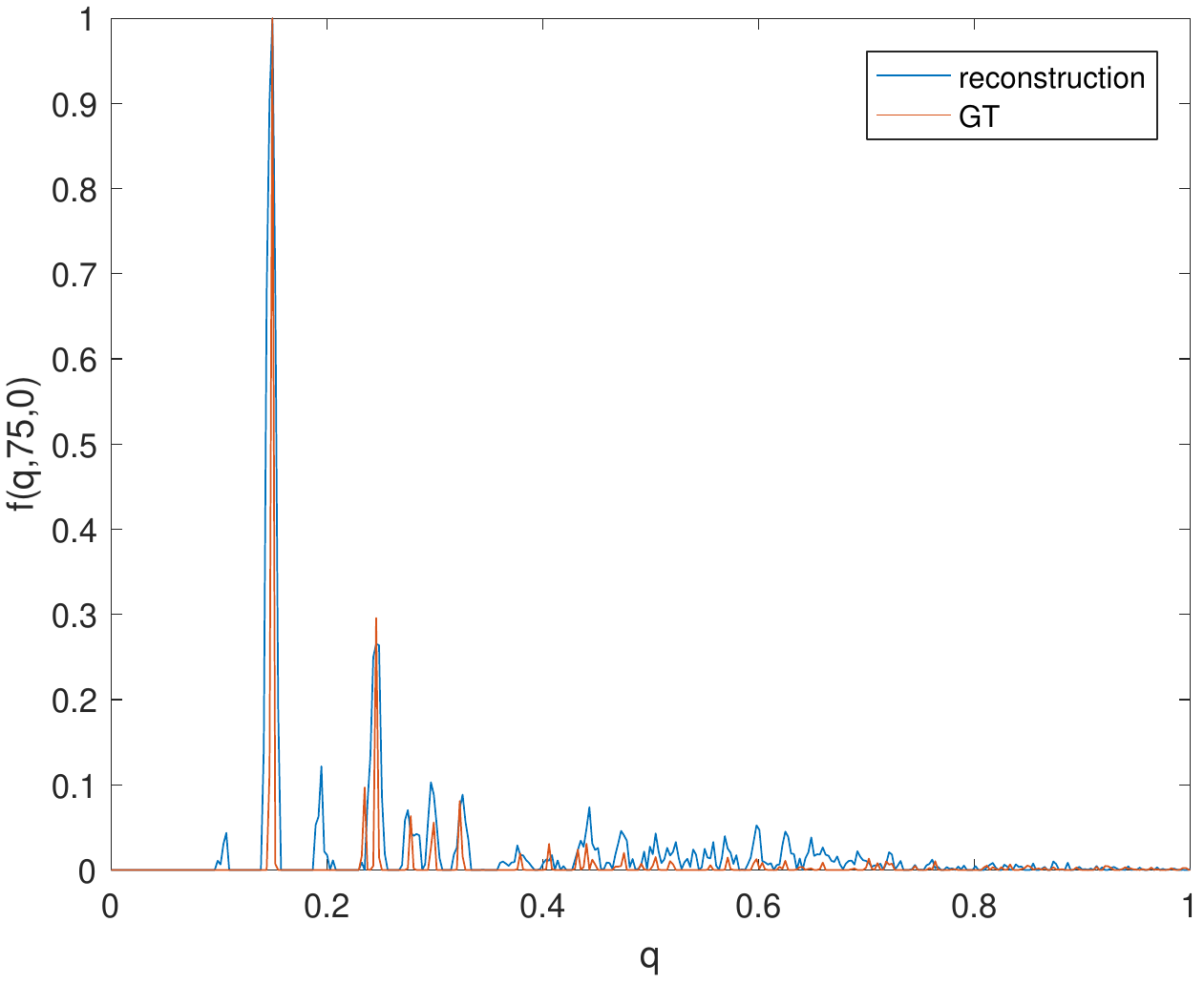} &
   \includegraphics[ width=0.4\linewidth, height=0.4\linewidth, keepaspectratio]{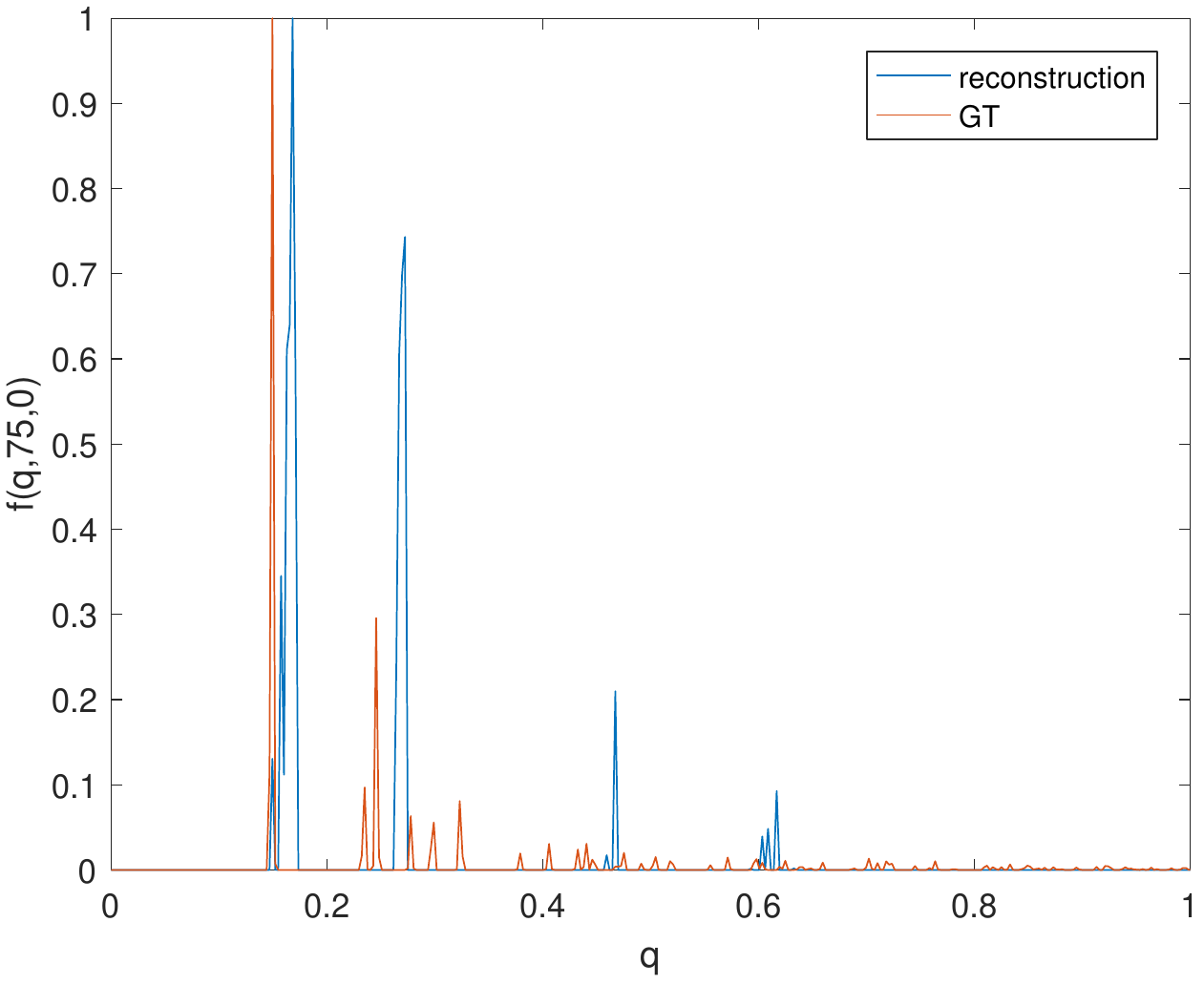} \\ 
   
 \rotatebox{90}{\hspace{1cm} LP2 ($x_1=-75$mm)} &
   \includegraphics[ width=0.4\linewidth, height=0.4\linewidth, keepaspectratio]{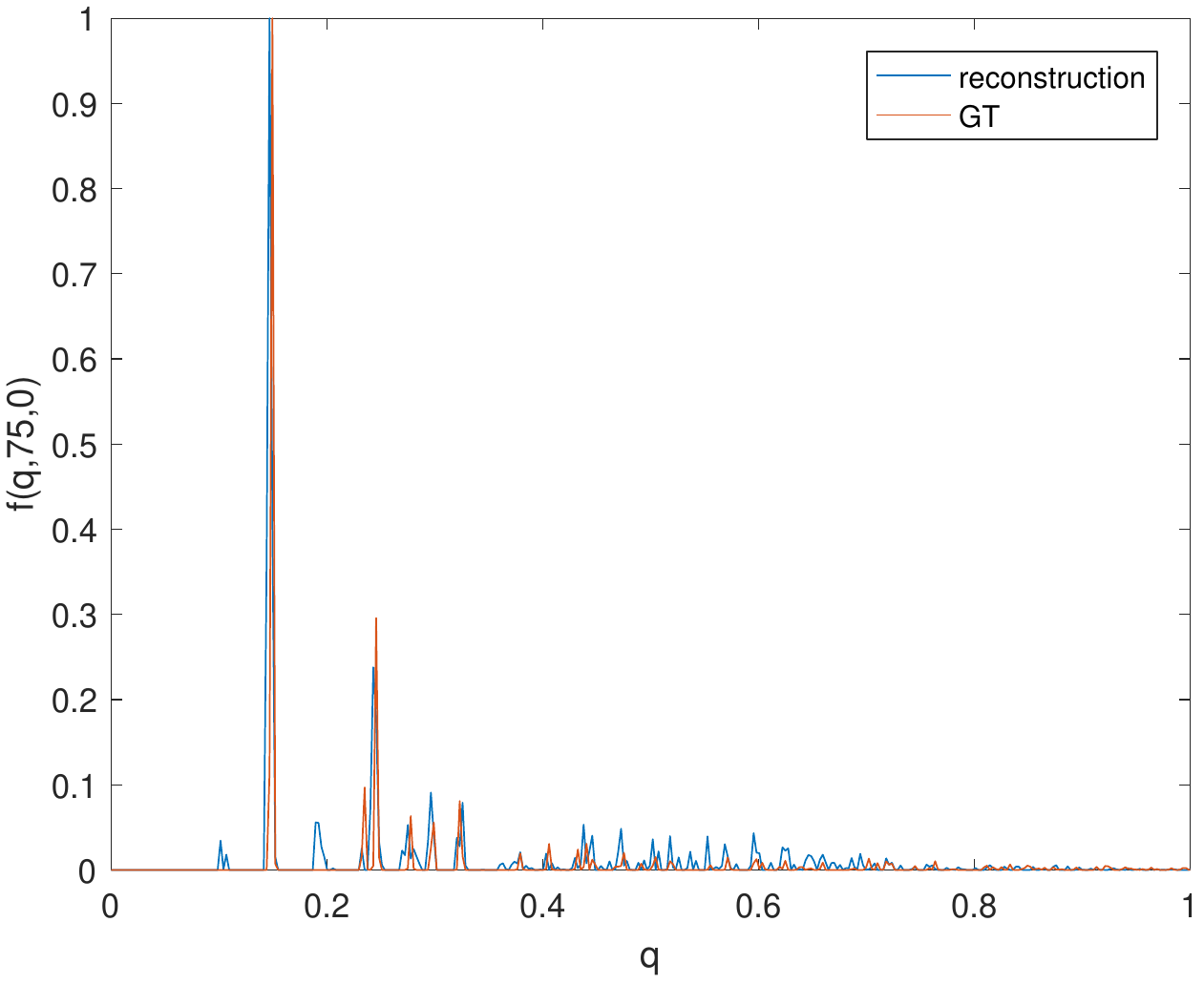} &
   \includegraphics[ width=0.4\linewidth, height=0.4\linewidth, keepaspectratio]{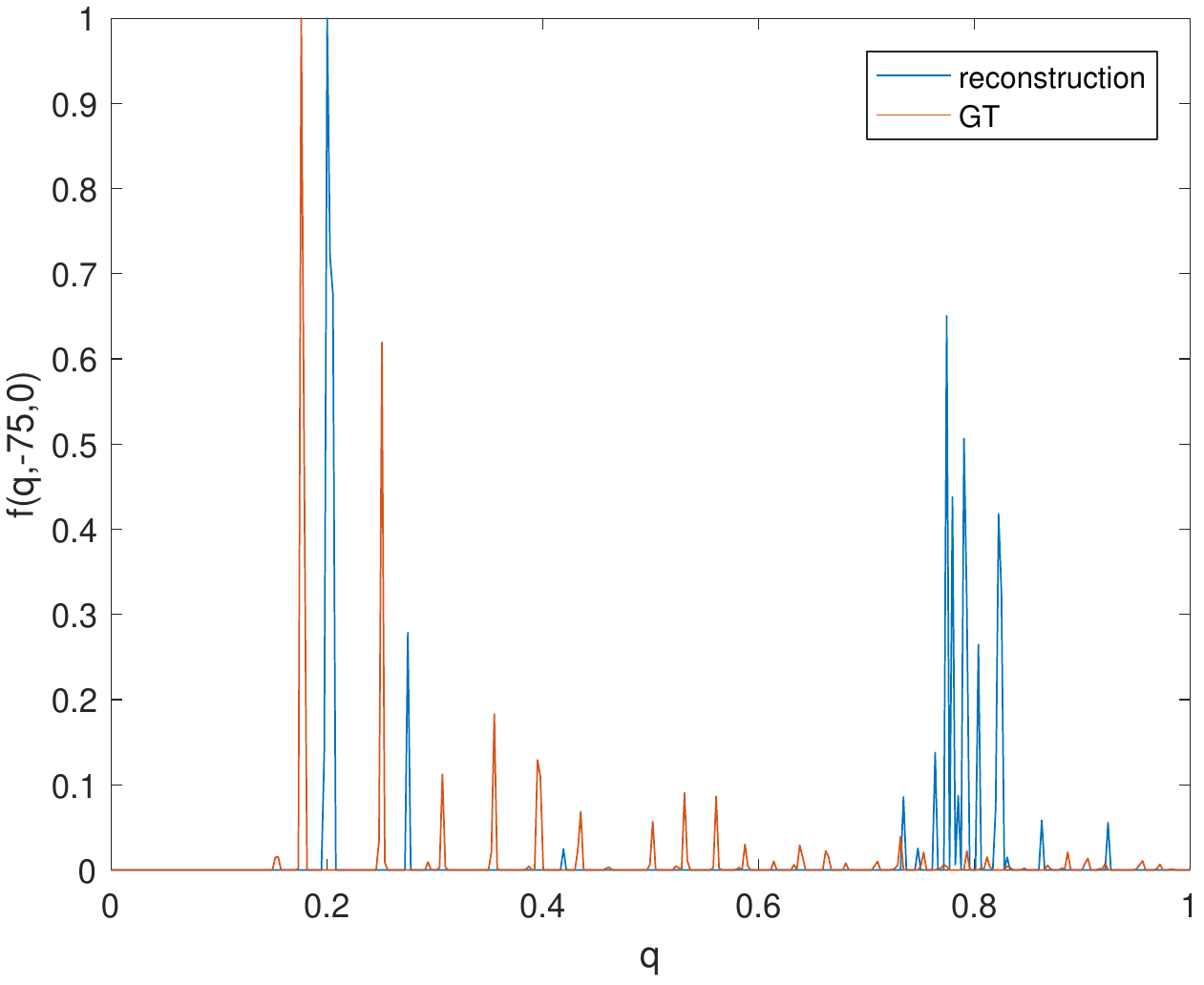} \\ 
 
\end{tabular}
\caption{NaCL and C-graphite sphere reconstruction, centers $x_c=-75$mm and $x_c=75$mm respectively, both with $r=10$mm.}
\label{F3}
\end{figure*}

\begin{figure*}
\centering
\setlength{\tabcolsep}{5pt}
\begin{tabular}{ c|cc }  
  &  2DBSR (S2) & FTV \\ \hline \\[-0.4cm] 
 \rotatebox{90}{\hspace{2.5cm} GT}  & 
   \includegraphics[ width=0.4\linewidth, height=0.4\linewidth, keepaspectratio]{StStShift_GT} &
   \includegraphics[ width=0.4\linewidth, height=0.4\linewidth, keepaspectratio]{StStShift_GT} \\ 
   
 \rotatebox{90}{\hspace{1cm} 2-D reconstruction} &
  \includegraphics[ width=0.4\linewidth, height=0.4\linewidth, keepaspectratio]{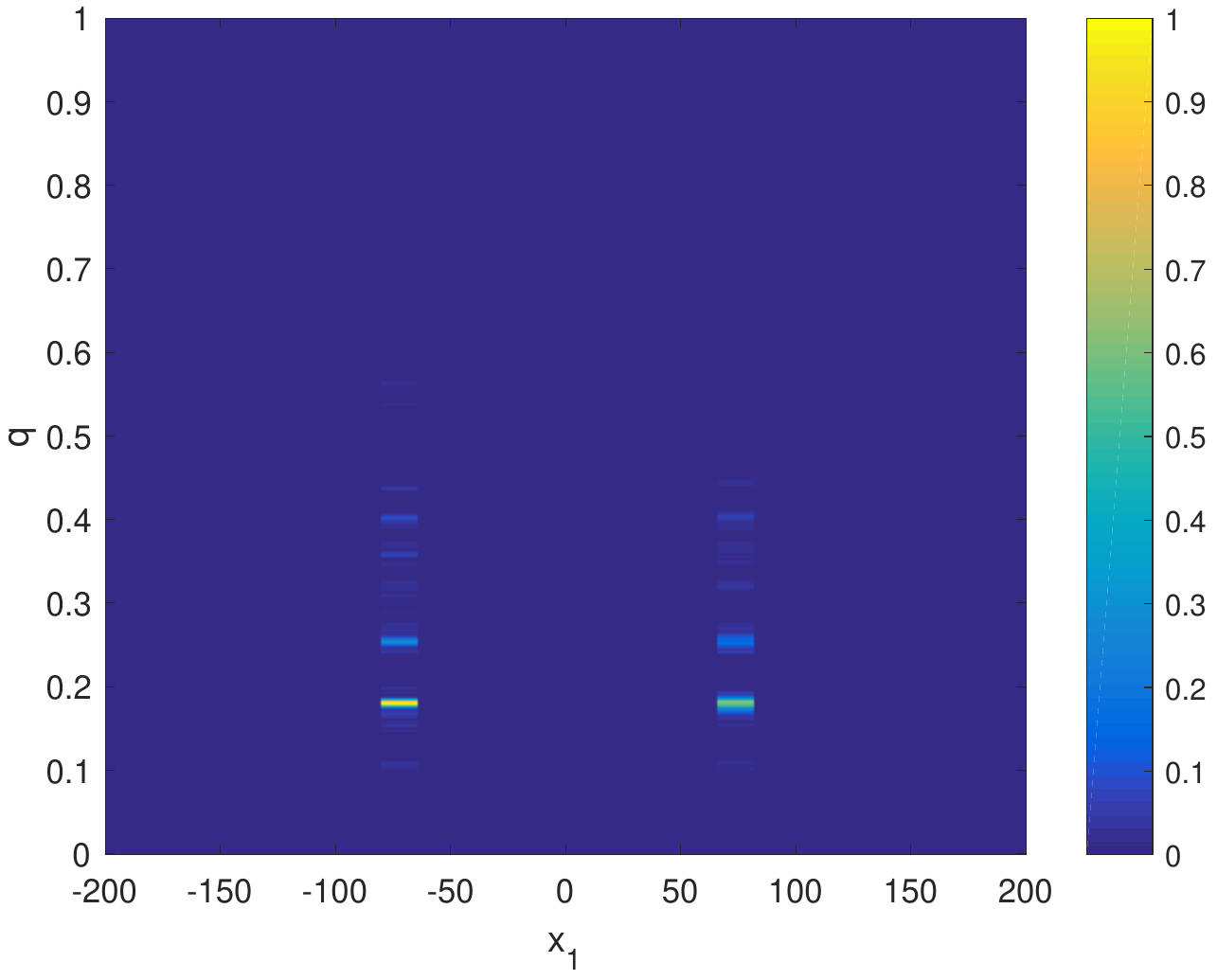} &
   \includegraphics[ width=0.4\linewidth, height=0.4\linewidth, keepaspectratio]{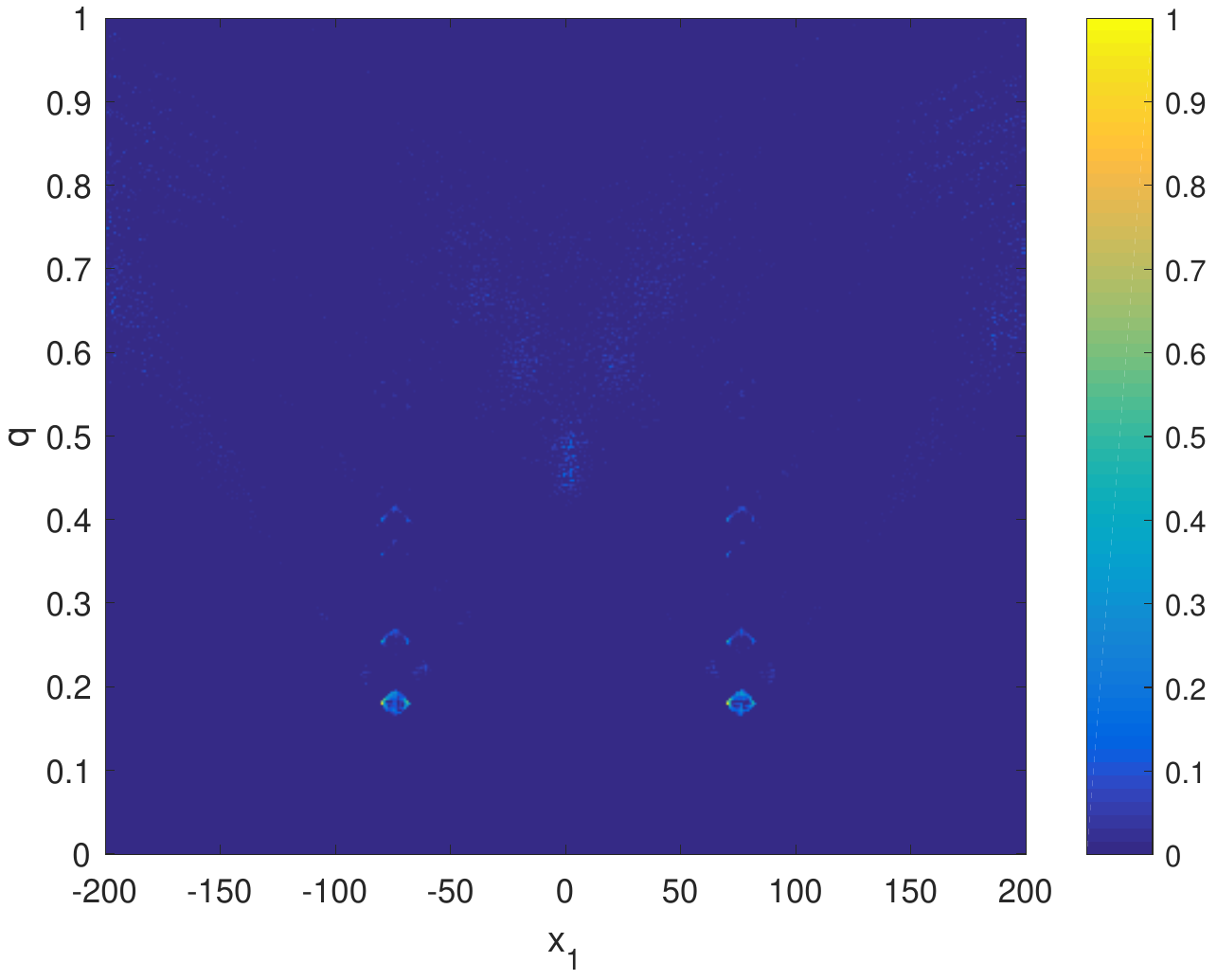} \\ 
   
 \rotatebox{90}{\hspace{1cm} LP1 ($x_1=75$mm)} &
   \includegraphics[ width=0.4\linewidth, height=0.4\linewidth, keepaspectratio]{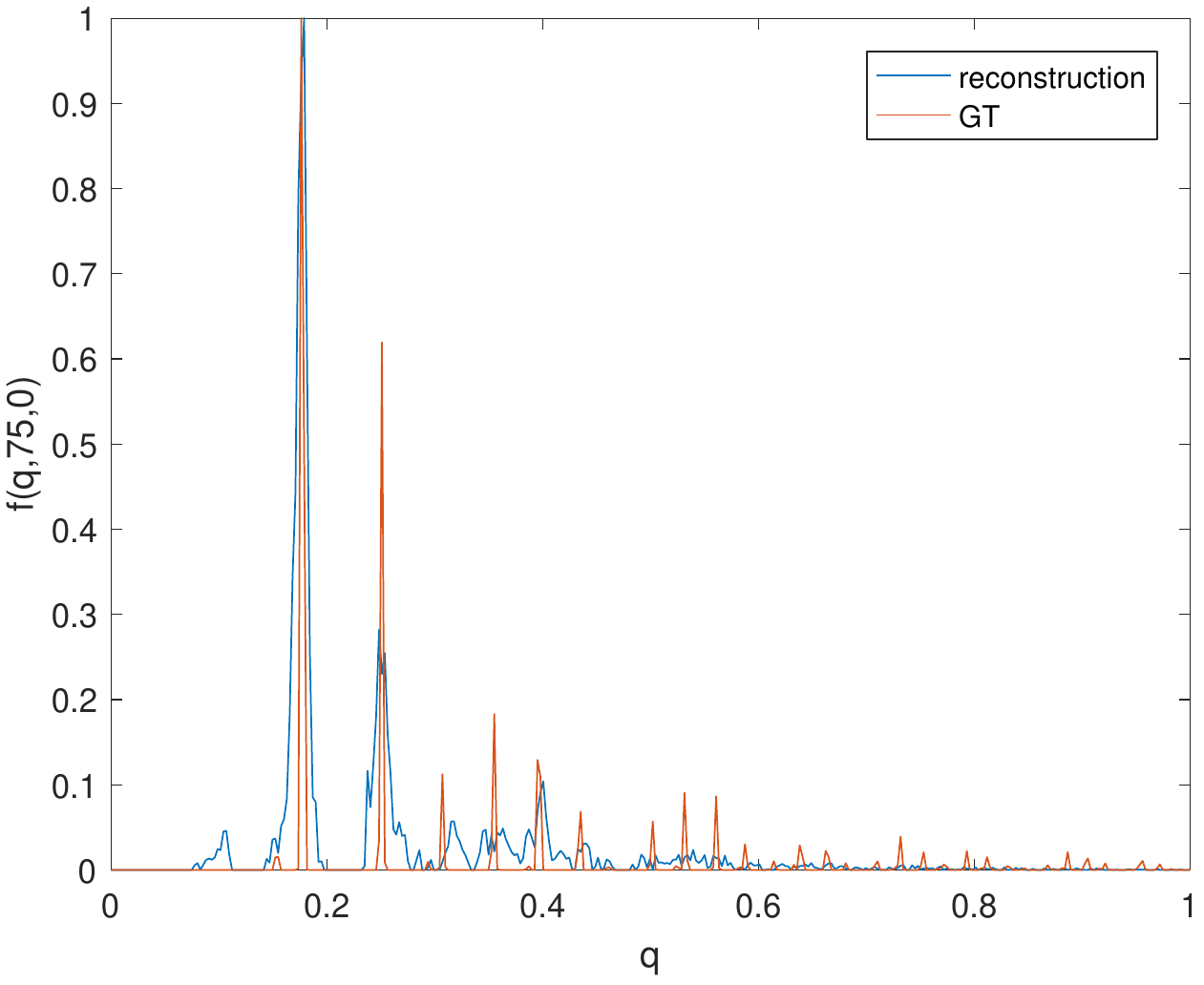} &
   \includegraphics[ width=0.4\linewidth, height=0.4\linewidth, keepaspectratio]{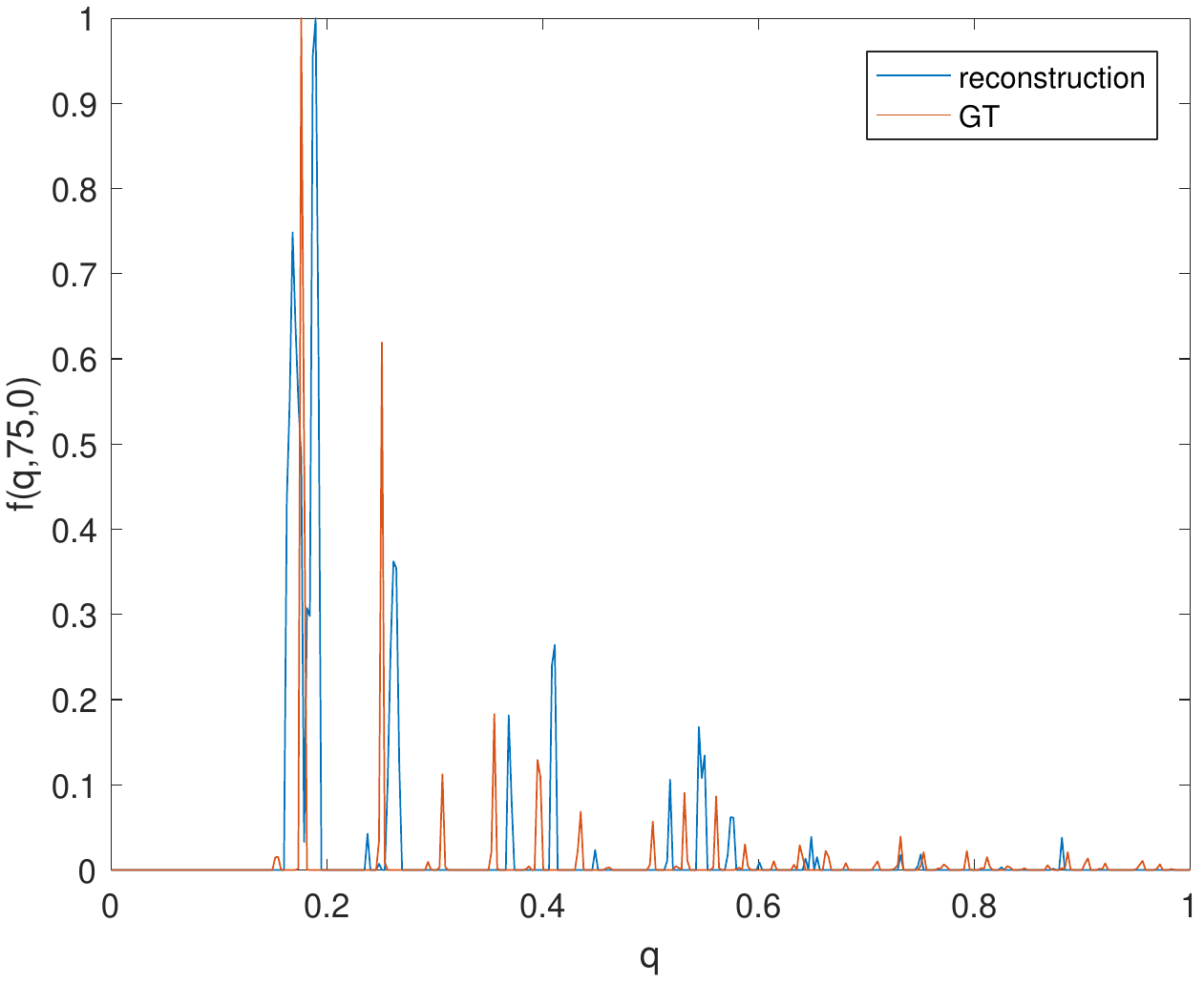} \\ 
  
 \rotatebox{90}{\hspace{1cm} LP2 ($x_1=-75$mm)} &
   \includegraphics[ width=0.4\linewidth, height=0.4\linewidth, keepaspectratio]{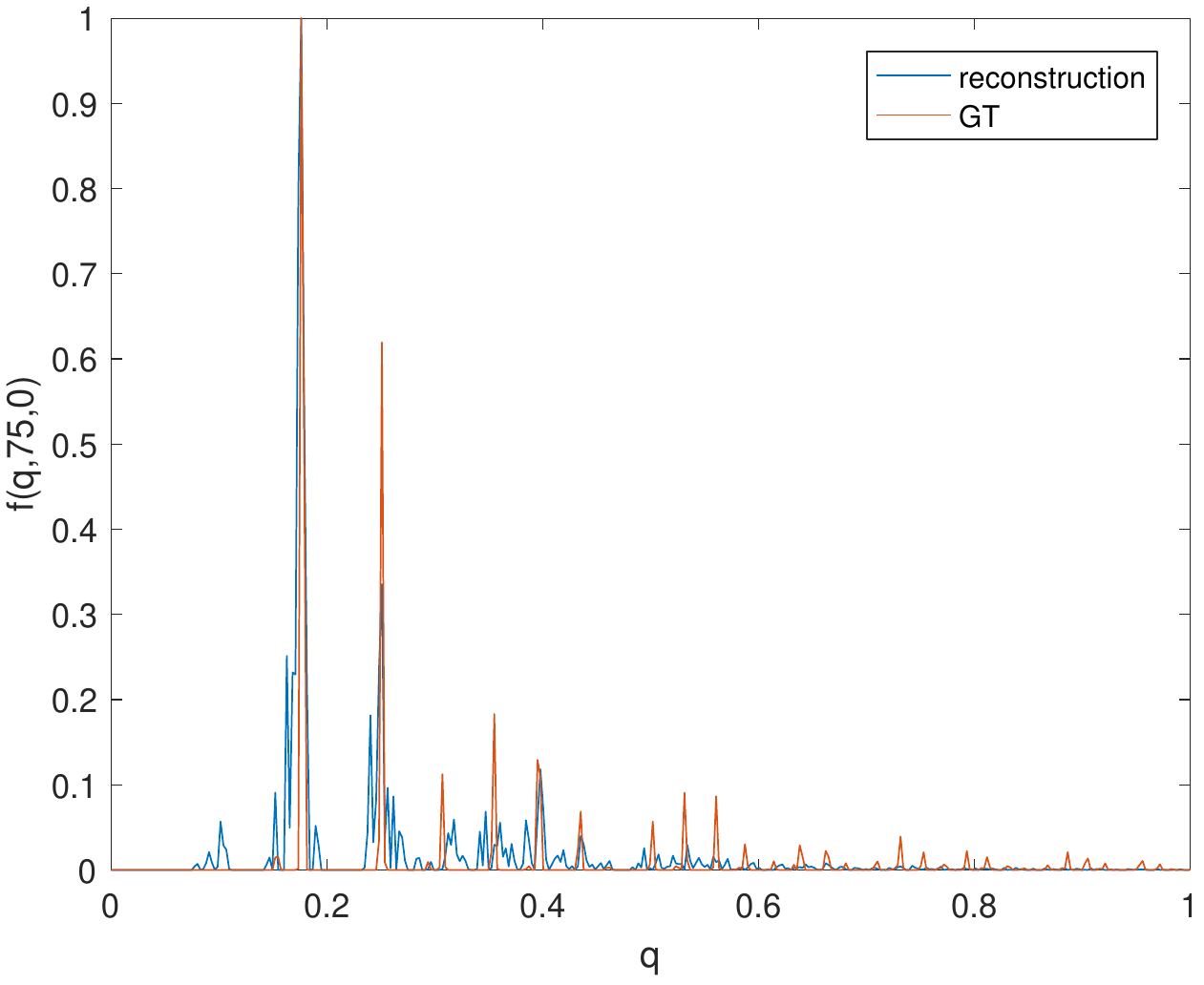} &
   \includegraphics[ width=0.4\linewidth, height=0.4\linewidth, keepaspectratio]{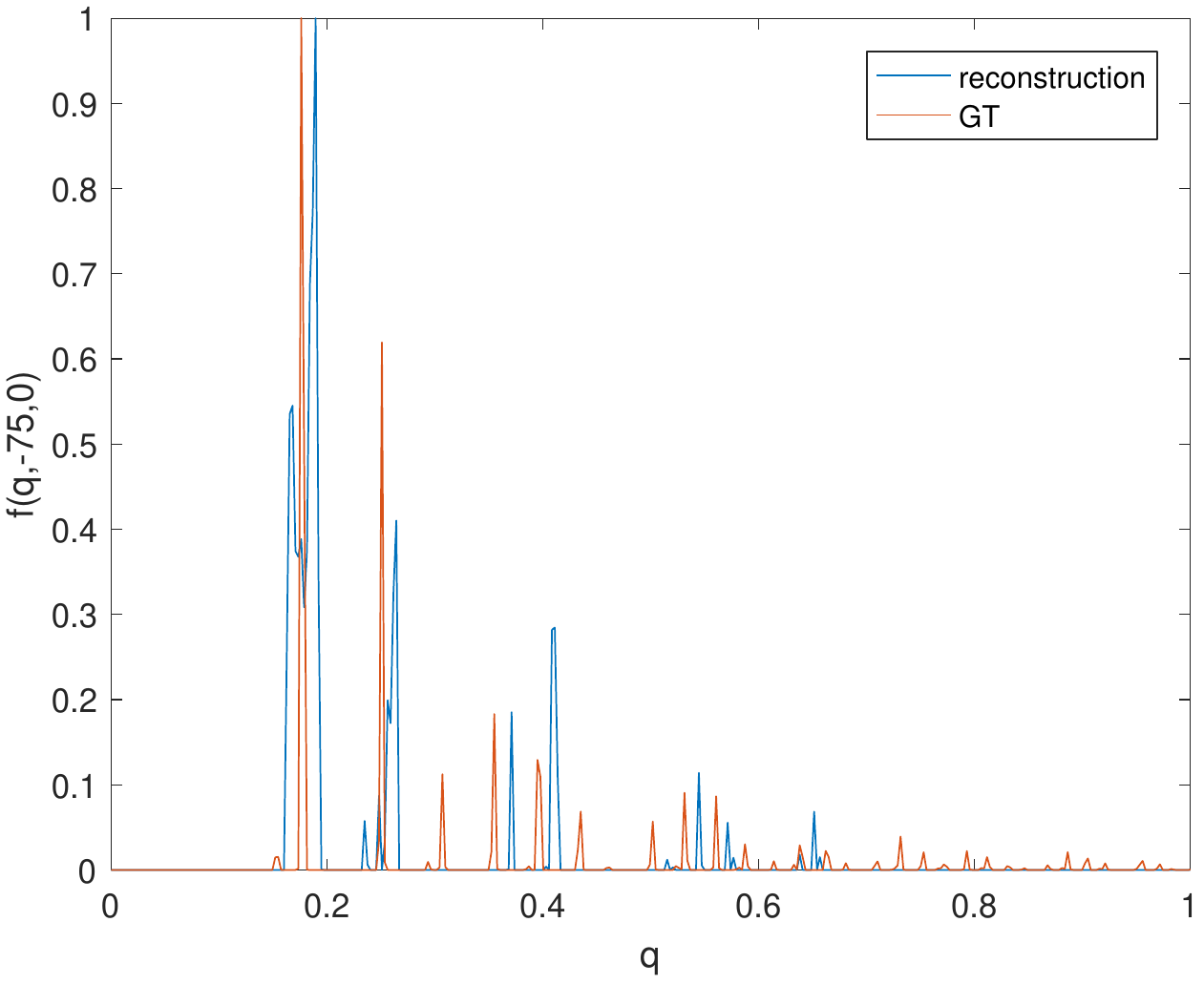} \\ 
  
\end{tabular}
\caption{Two NaCl spheres reconstruction, centers $x_c=-72.5$mm and $x_c=77.5$mm, both with $r=6.25$mm.}
\label{F4}
\end{figure*}

To generate the data, we use a novel, single scatter Monte Carlo (MC) code developed in Matlab, which we introduce in this paper. A step-by-step description of the MC code and a pseudo code is provided in the appendix. The full code is available from the authors upon request. The code fully incorporates noise effects due to attenuation and self-absorption, clutter effects and Compton scatter. We do not assume knowledge of the attenuation coefficient in our simulations, and thus there is systematic error in the data due to the neglection of attenuation modeling. In the energy and scattering angle range considered ($E<30$keV and $\omega<\frac{\pi}{2}$) the Compton scatter intensity is significantly lower than that of the Bragg signal. Hence the Compton scatter will act as an additional background noise in our simulations, along with the Poisson noise due to photon arrivals. 

In figure \ref{figERR} we have plotted the systematic attenuation error and Compton signal for a single sphere of Carbon-graphite with center $x_c=0$, for varying $r$ and $E$. We see in figure \ref{muErr} that as $r$ increases and $E$ decreases (towards the bottom-right corner of the figures) the systematic error due to attenuation increases due to higher self-absorption effects at low energy (fewer photons penetrate) and higher r (those that suffer greater attenuation due to Beer's Law). We see the converse effect for the Compton noise in figure \ref{CompSig} as much of the Compton scatter is self-absorbed for low $E$ and high $r$. This is as we'd expect and figure \ref{muErr} appears as the negative image of figure \ref{CompSig} (after scaling and translation). So there is a trade-off to consider here, and no size of object or photon energy is necessarily preferred over the other. The MC experiments considered here are designed to reflect a variety of levels of Compton noise and attenuation error in the data (i.e. we consider a variety of $r$ and $E$).  

To obtain the 2DBSR reconstructions presented in this section we implement the following two-stage reconstruction procedure. First, we apply Algorithm \ref{algorithm1}, with the filtered Bragg data $\vb$ (using the filters specified in section \ref{microlocal}), characteristic library and hyperparameters (as specified in section \ref{IP}) as input, choosing $\lambda$ for the best results in terms of edge $F_1$ score. Then, we perform a second full pass of Algorithm \ref{algorithm1}, replacing $\lambda\to \lambda\times 10$. We found that increasing $\lambda$ by an order of magnitude on the second pass of Algorithm \ref{algorithm1} helped to reduce background noise in the reconstruction, once the more significant artifacts (e.g. boundary artifacts) had been removed in the first pass of Algorithm \ref{algorithm1}. To show how this works, see figure \ref{stages} for image reconstructions of an NaCl sphere, center $x_c=0$mm, radius $r=15$mm, at each stage of the two-stage process. We notice an improvement in the image quality and $F_1$ score as we advance from stage 1 to stage 2.

Upon applying the same two-stage process to the FTV reconstructions, we did not see an improvement in the results, and hence why we did not apply the two-stage process to FTV. This is as expected since the objective of \eqref{equTV} is convex, and hence, as we chose $\lambda$ for FTV for the best results, increasing $\lambda\to \lambda\times 10$ gave the same result (i.e. the best result for FTV over all $\lambda$). 

In table \ref{T3} we have presented our results in terms of edge $F_1$ score for the nine MC experiments conducted. In table \ref{T3} we use the shorthand notation ``St" for NaCl (salt), ``C" for C-graphite and ``Di" for C-diamond. In brackets we specify the centers and radii of the spheres in millimeters. For example, ``St ($x_c=0$, $r=10$)" means a single NaCl sphere with center $x_c=0$mm and radius $r=10$mm, and ``St ($x_c=-75, r=5$), Di ($x_c=75, r=5$)" means an NaCl and C-diamond sphere with center $\pm 75$mm and radius $r=5$mm. Stages 1 and 2 of 2DBSR detailed above are denoted by S1 and S2 in table \ref{T3}. We see a significantly improved performance at all stages of 2DBSR over FTV, in all experiments conducted. We found that S2 offered better image quality and $F_1$ score overall, when compared to S1.

In table \ref{T3} we have highlighted five rows in bold. These correspond to the image reconstructions in figures \ref{F0}-\ref{F4}, which are chosen to highlight the best and worst results using our method (at stage 2), for materials within and outside of the imaging basis. In figures \ref{F0}-\ref{F4}, ``GT" denotes the ``Ground Truth", and ``LP" denotes ``Line Profile". We specify the line profile equation in brackets in each case. In figure \ref{F0} we have presented image reconstructions of a C-graphite sphere, center $x_c=0$mm and radius $r=15$mm. In this case the attenuation error is high and the Compton noise is low (see figure \ref{figERR}).  We see structured distortions (systematic error) of the GT Bragg peaks in the FTV reconstruction. That is, the horizontal line segments of the GT are stretched in the positive $q$ direction, and we see a shifting of the Bragg peaks in the $x_1=0$mm LP.  Given the assumed structure of $f$ when using 2DBSR (as in equation \eqref{fform1}), the 2DBSR algorithm does not allow for such structured distortions in the reconstruction. Hence we see a complete removal of the artifacts in the 2DBSR reconstruction, with an increase in $F_1$ score as a result. In figure \ref{F1} we see a similar effect although the $F_1$ using 2DBSR is decreased as, in this case, the C-diamond sphere is out of basis. The $F_1$ score using FTV is improved however, and the shifting of the Bragg peaks is reduced in the reconstruction. This is because the attenuation error is decreased due to the smaller sphere radius ($r=8.75$mm).

In figure \ref{F2}, we have presented reconstructions of NaCL and C-diamond spheres within the imaging basis, with centers $x_c=-75$mm and $x_c=75$mm respectively, both with $r=10$mm. Here we see some significant noise effects in-between the Bragg peaks in the 2DBSR reconstruction, at a detriment to the $F_1$ score ($F_1=0.65$) as some of the noise is mistaken for additional (unwanted) Bragg peaks. We may be able to combat this by continuing to increase $\lambda$ (the $L_1$ regularization parameter) and running further iterations of 2DBSR as in the second stage of the two-stage reconstruction process outlined above, with the idea that the sparse regularizers would remove the noise and focus on the more significant (true) Bragg peaks in the reconstruction. The centers and widths of the spheres in the 2DBSR reconstruction are recovered correctly, and the $F_1$ score is acceptable. Furthermore, the results using 2DBSR are much improved over FTV, where we see severe artifacts in the reconstruction, with low $F_1$ score ($F_1=0.10$).

In figure \ref{F3} we have presented $f(q,x_1)$ reconstructions of a NaCl and C-graphite sphere with centers $x_c=-75$mm and $x_c=75$mm, repsectively, both with radius $r=10$mm. We see that the NaCl sphere is effectively zeroed out in the 2DBSR and FTV reconstruction, leading to a low $F_1$ score using both methods. . The reason for the poor image quality is that the scattered signal contribution from the Carbon-graphite sphere is more than twice that of the salt sphere. The sum of counts (taken over all $p=539400$ data points) scattering from the salt sphere is $7.1\times 10^5$, and the sum of the counts from the Carbon-graphite sphere is $1.7\times 10^6$. The NaCl scatter therefore accounts for $29\%$ of the total signal. 
As the signal from the NaCl sphere is largely hidden by that of the Carbon, this makes the NaCl $f(q,x_1)$ function more difficult to recover, and it is set close to zero by 2DBSR and FTV. Similar effects are seen in conventional X-ray CT, for example, in head scans, where the bone and skull shows up strongly in the image (due to higher absorption) and the skull interior and soft tissue are harder to see \cite[figure 1]{ehrhardt2014joint}.

In figure \ref{F4} we show $f(q,x_1)$ reconstructions of two NaCl spheres with centers $x_c=-72.5$mm and $x_c=77.5$mm, respectively, both with radius $r=6.25$mm. Thus, in this example, the NaCl spheres are outside the 2DBSR imaging basis. In this case we notice that the Bragg peaks of the right-hand NaCl sphere are blurred in the 2DBSR reconstruction. The blurring effect is mild however and the image quality is improved when compared to FTV. In the FTV reconstruction the shifting of the Bragg peaks, although still prominent, is less significant than in figures \ref{F0}-\ref{F3}. We notice, in general, a reduction in the Bragg peak shift in the FTV reconstructions with decreasing $r$ (i.e. as the attenuation error decreases). This observation is supported by the results of table \ref{T3}, as the $F_1$ scores using FTV increase (and become more competitive with 2DBSR) with decreasing $r$. Thus it seems that the main cause of error in the FTV reconstructions is due to attenuation modeling, based on the error analysis of figure \ref{figERR}. 

\section{Conclusions and further work}
Here we have presented the 2DBSR algorithm - a novel reconstruction technique for two-dimensional Bragg scatter imaging. The regularization penalties we applied are based on ideas in multibang control and compressive sensing. We also incorporated filtering ideas from microlocal analysis for data pre-processing, with the aim to suppress boundary type artifacts (e.g. such as those observed in \cite{webber2020microlocal}) in the reconstruction. In section \ref{reconmethod} we formalized our approach and detailed the 2DBSR algorithm in section \ref{algorithm1}. The microlocal filters used for data pre-processing were defined as polynomials in $E$ in section \ref{microlocal}. In section \ref{results} we designed and conducted a variety of BST experiments, using analytic and Monte Carlo data, which are of interest in airport baggage screening and threat detection. Here we compared the performance of 2DBSR to a Filtered Total Variation (FTV) approach. 
In the analytic data experiments both methods performed well, with FTV slightly outperforming 2DBSR in terms of $F_1$ score, in most cases. In the Monte Carlo experiments, 2DBSR was shown to offer significantly higher performance than that of FTV in all cases considered (nine MC experiments in total). In the FTV reconstructions we saw severe image artifacts due to errors in the attenuation modeling. The artifacts appeared in the image as a shifting/distorting of the ground truth Bragg spectra in the positive $q$ direction, which caused a reduction in $F_1$ score.  In the majority of cases, the 2DBSR algorithm was successful in combatting the image artifacts due to attenuation modeling (among other errors, such as those due to Compton scatter), and the $F_1$ score was much improved over FTV. 

In future work we aim to consider joint electron density/Bragg spectra reconstruction using a combination of Bragg and Compton scatter data. Here we aim to use the models from the Compton scatter tomography literature \cite{2D1,2D2,2D3,2D4,3D1,3D3,3D4} to improve the forward model \ref{equBG}, and the recovery of $f(q,x_1)$ in $x_1$ space. That is, the Compton data will be primarily used for density recovery and locating the crystallites (i.e. recovering the centers and widths of the materials), with the recovery in $q$ space coming mainly from the Bragg data. We expect that such ideas will be effective in combatting the more significant errors observed in figures \ref{F3} and \ref{F4} in the 2DBSR reconstructions, where we saw the zeroing out of the NaCl sphere (in figure \ref{F3}) and blurring of the ground truth Bragg spectra (in figure \ref{F4}), which yielded poor $F_1$ scores. 

\section*{Funding}
This material is based upon work supported by the U.S. Department of Homeland Security Science and Technology Directorate, Office of University Programs, under Grant Awards 18STEXP00001-03-02, formerly 2013-ST-061-ED0001 and 70RSAT19FR0000155. The views and conclusions contained in this document are those of the authors and should not be interpreted as necessarily representing the official policies, either expressed or implied, of the U.S. Department of Homeland Security.

\section*{Acknowledgments}
We gratefully acknowledge the financial support from the U.S. Department of Homeland Security Science and Technology Directorate.

\section*{Disclosures}
The authors declare no conflicts of interest.

\section*{Data availability statement}
Data underlying the results presented in this paper are not publicly available at this time but may be obtained from the authors upon reasonable request.

\appendix
\section*{Appendix. Description of Monte Carlo code}
\label{A1}
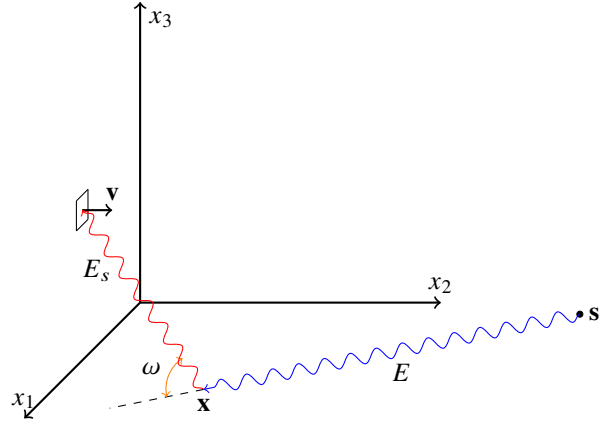
\begin{figure}[!htb]
\centering
\begin{tikzpicture}[scale=4]
\path (1.5,0,0.1) coordinate (S);
\path (0,0.5,0.5) coordinate (D);
\path (0.5,0,0.75) coordinate (w);
\path (0.25,0,0.9125) coordinate (a);
\draw [thin,dashed] (w)--(a);
\node at (0.05,0.3,0.5) {$E_s$};
\node at (1.1,0,0.6) {$E$};
\draw [fill] (S) circle (0.01);
\draw pic[draw=orange, <->,"$\omega$", angle eccentricity=1.5] {angle = D--w--a};
\draw[thick,->] (0,0,0) -- (1,0,0) node[above]{$x_2$};
\draw[thick,->] (0,0,0) -- (0,1,0) node[anchor=north west]{$x_3$};
\draw[thick,->] (0,0,0) -- (0,0,1) node[anchor=south]{$x_1$};
\draw (0,0.45,0.45)--(0,0.55,0.45);
\draw (0,0.45,0.45)--(0,0.45,0.55);
\draw (0,0.55,0.55)--(0,0.45,0.55);
\draw (0,0.55,0.55)--(0,0.55,0.45);
\draw[thick,->] (0,0.5,0.5)--(0.1,0.5,0.5)node[above]{$\mathbf{v}$};
\draw [snake it, ->,blue](S)node[right]{\textcolor{black}{$\vs$}}--(w)node[below]{\textcolor{black}{$\vx$}};
\draw [snake it, ->,red](w)--(D);
\end{tikzpicture}
\caption{A scattering event occurs at a scattering site $\vx$, for photons emitted from a source $\vs$ and recorded at a detector $\vd$ (displayed as a small square in the $x_1x_3$ plane). The initial photon energy is $E$ and the scattered energy is $E_s$. Here $\mathbf{v}$ is the direction normal to the detector surface. The scattering angle is $\omega$.}
\label{fig3}
\end{figure}

Consider the scattering event pictured in figure \ref{fig3}.
Let $L_{\vx_1\vx_2}=\{\vx_1+t(\vx_2-\vx_1) : t\in\mathbb{R}\}$ be the line through $\vx_1$ and $\vx_2$, and let us consider the scattering interaction space $\Omega=\{\text{NI},\text{PE},\text{INC},\text{COH}\}$, where NI denotes ``no interaction" (transmitted photons), PE denotes ``photoelectric absorption", INC denotes ``incoherent scatter", and COH denotes ``coherent scatter". Let $\mu(E,Z)=\sum_{\mathcal{E}\in \Omega/\{\text{NI}\}}\mu_{\mathcal{E}}(E,Z)$ denote the total attenuation coefficient, where $\mu_{\mathcal{E}}$ is the attenuation coefficient for the interaction type $\mathcal{E}\in\Omega/\{\text{NI}\}$, and $Z$ is the atomic number. Then the Monte Carlo algorithm, for a single photon trial, scattering site ($\vx$) and source position ($\vs$), reads as follows:
\begin{enumerate}
\item Sample the initial photon energy $E$ from the source spectrum. The photon is emitted from $\vs$ and travels in the direction $\vx-\vs$ towards $\vx$.
\item Sample from a Bernoulli distribution, where $e^{-\int_{L_{\vs\vx}}\mu(E,Z)}$ is the probablity of success, to decide if the photon reaches $\vx$.
\item Determine the photon interaction by sampling from a finite distribution on $\Omega$, with $P(\text{NI})=e^{-\mu(E,Z)}$ and
\begin{equation}
P(\mathcal{E})=(1-P(\text{NI}))\frac{\mu_{\mathcal{E}}(E,Z)}{\mu(E,Z)}
\end{equation}
for $\mathcal{E}\in \Omega/\{\text{NI}\}$.
\item Sample $\omega\in[0,\pi]$ (polar angle of the scatter direction) from the differential cross section distribution for the given interaction type, keeping a uniform spread over each scattering cone, namely the cone with central axis direction $\vx-\vs$ (the dashed line in figure \ref{fig3}), vertex $\vx$ and opening angle $\omega$. That is we sample the azimuth angle of the scatter direction (in $[0,2\pi]$) from a uniform distribution.
\item If the scattered photon travels in the direction $\vd-\vx$ towards $\vd$, sample from a second Bernoulli distribution, with $e^{-\int_{L_{\vx\vd}}\mu(E_s,Z)}$ as the probablity of success, to decide whether the photon reaches $\vd$.
\item If the scattered photon reaches $\vd$, count one photon with energy $E_s$, else continue. Repeat for all detectors in the array.
\end{enumerate}
The psuedo code detailed above is applied to the geometry of figure \ref{figmain}, with the differential cross section distribution set to one of the Bragg spectra of figure \ref{Fq1}, depending on the material. For Compton scatter, the differential cross section is set to the Klein-Nishina distribution \cite{klein1928scattering,KN}. To generate the data used in section \ref{MC}, we expand the single trial, single scatter site code to multiple scattering events along lines parallel to $x_1$ (e.g. $L$ as pictured in figure \ref{figmain}). 

\bibliography{referencesBragg}

\end{document}